\documentclass[11pt]{amsart} 

\usepackage{amssymb,tensor,enumitem,verbatim,theoremref,mathtools,stmaryrd}
\usepackage[bookmarksdepth=3]{hyperref}
\usepackage[margin=1in]{geometry}
\usepackage[all,cmtip]{xy}
\usepackage[sort,noadjust]{cite}
\usepackage{nicematrix}
\usepackage[color=yellow]{todonotes}

\newcommand{\sm}[1]{[ \begin{smallmatrix} #1 \end{smallmatrix} ]}

\newcommand{\subgrp}[1]{\langle #1 \rangle}

\newcommand{\set}[1]{\left\{ #1 \right\}}

\newcommand{\abs}[1]{| #1 |}

\newcommand{\wt}[1]{\widetilde{ #1}}

\DeclareMathOperator{\ConjClo}{ConjClo}

\DeclareMathOperator{\End}{End}
\DeclareMathOperator{\Grp}{Grp}
\DeclareMathOperator{\Hom}{Hom}
\DeclareMathOperator{\Id}{Id}

\DeclareMathOperator{\Ind}{Ind}
\DeclareMathOperator{\Inf}{Inf}

\DeclareMathOperator{\Lie}{Lie}

\DeclareMathOperator{\res}{res}

\DeclareMathOperator{\Res}{Res}
\DeclareMathOperator{\rk}{rk}

\DeclareMathOperator{\Span}{span}

\newcommand{\ve}{\varepsilon}

\newcommand{\C}{\mathbb{C}}

\newcommand{\Q}{\mathbb{Q}}
\newcommand{\R}{\mathbb{R}}
\newcommand{\T}{\mathbb{T}}
\newcommand{\Z}{\mathbb{Z}}

\newcommand{\CB}{\C \calB}
\newcommand{\CD}{\C \calD}
\newcommand{\CG}{\C G}

\newcommand{\CS}{\C \fS}
\newcommand{\CW}{\C W}

\newcommand{\RB}{\R \calB}

\newcommand{\calB}{\mathcal{B}}
\newcommand{\calD}{\mathcal{D}}
\newcommand{\calH}{\mathcal{H}}
\newcommand{\calL}{\mathcal{L}}

\newcommand{\calT}{\mathcal{T}}
\newcommand{\calW}{\mathcal{W}}
\newcommand{\calX}{\mathcal{X}}

\newcommand{\calBP}{\mathcal{BP}}

\newcommand{\fb}{\mathfrak{b}}

\newcommand{\fd}{\mathfrak{d}}

\newcommand{\g}{\mathfrak{g}}
\newcommand{\gl}{\mathfrak{gl}}
\newcommand{\fh}{\mathfrak{h}}
\newcommand{\osp}{\mathfrak{osp}}
\newcommand{\fr}{\mathfrak{r}}
\newcommand{\fs}{\mathfrak{s}}
\newcommand{\fS}{\mathfrak{S}}
\newcommand{\fsl}{\mathfrak{sl}}
\newcommand{\fso}{\mathfrak{so}}
\newcommand{\fsp}{\mathfrak{sp}}

\numberwithin{equation}{subsection}
\counterwithin{table}{section}

\newtheorem{theorem}{Theorem}[subsection]
\newtheorem*{theorem*}{Theorem}
\newtheorem{proposition}[theorem]{Proposition}
\newtheorem*{proposition*}{Proposition}

\newtheorem{corollary}[theorem]{Corollary}
\newtheorem{lemma}[theorem]{Lemma}

\theoremstyle{definition}

\newtheorem{definition}[theorem]{Definition}

\newtheorem{remark}[theorem]{Remark}
\newtheorem{convention-construction}[theorem]{Convention/Construction}

\newtheorem*{definition*}{Definition}

\title{Lie algebras generated by reflections in types BCD}

\author{Christopher M.\ Drupieski}
\address{Department of Mathematical Sciences,
		DePaul University,
		Chicago, IL 60614, USA}
\email{c.drupieski@depaul.edu}

\thanks{CMD was supported by a Summer Research Grant from the DePaul University College of Science and Health.}

\author{Jonathan R.\  Kujawa}
\address{Department of Mathematics \\
		Oregon State University \\
		Corvallis, OR 97331, USA}
\email{kujawaj@oregonstate.edu}

\thanks{JRK was supported in part by Simons Collaboration Grants for Mathematicians No.\ 525043 and No.\ 963912.}

\date{\today}


\allowdisplaybreaks

\begin{document}

\begin{abstract}
We consider the group algebra over the field of complex numbers of the Weyl group of type B (the hyperoctahedral group, or the group of signed permutations) and of the Weyl group of type D (the demihyperoctahedral group, or the group of even-signed permutations), viewed as Lie algebras via the commutator bracket, and determine the structure of the Lie sub\-algebras generated by the sets of reflections.
\end{abstract}



\maketitle

\section{Introduction}

\subsection{The Lie algebra of transpositions}

Let $W$ be a finite group, let $\CW$ be its group algebra over the field of complex numbers, and write $\Lie(\CW)$ for $\CW$ considered as a Lie algebra via the commutator bracket. It is an immediate corollary of Maschke's theorem and the Artin--Wedderburn theorem that $\Lie(\CW)$ is a direct sum of general linear Lie algebras. Namely, if $V_1,\ldots,V_m$ is a complete set of pairwise non-isomorphic simple $\CW$-modules, then
	\[
	\Lie(\CW) \cong \gl(V_1) \oplus \cdots \oplus \gl(V_m).
	\]
Now suppose $W$ is a finite reflection group (e.g., the Weyl group of an irreducible root system) with set of reflections $S$. A more subtle problem in this context is to determine the structure of the Lie subalgebra of $\Lie(\CW)$ generated by the set $S$.

Motivated by questions from the representation theory of the braid group, Marin \cite{Marin:2007} solved the preceding problem in the case $W = \fS_n$ of the symmetric group on $n$ letters (i.e., the Weyl group of type $A_{n-1}$) by computing the structure of the Lie algebra $\fs_n$ generated by the set of transpositions in $\fS_n$. Recall that the isomorphism classes of the simple $\CS_n$-modules are labeled by integer partitions of $n$. Marin showed that $\fs_n$ is a reductive Lie algebra whose semisimple derived subalgebra $\fs_n' = [\fs_n,\fs_n]$ is a direct sum of special linear Lie algebras, indexed by equivalence classes of partitions whose Young diagrams are non-symmetric and not hooks, and of orthogonal and symplectic Lie algebras, indexed by partitions whose Young diagrams are symmetric and not hooks. 
A more detailed summary of Marin's result is given below in Section \ref{subsubsec:improper-bipartitions}. The calculation of $\fs_n$ enabled Marin to determine the algebraic hull of the image of the braid group in each of the irreducible representations of the Iwahori--Hecke algebra of type A, and in the Iwahori--Hecke algebra itself. Earlier work by Marin \cite{Marin:2003} also showed that the computation of $\fs_n$ permits one to describe how tensor products of representations for the Iwahori--Hecke algebra decompose as representations over the braid group; see \cite{BMM:2017,Marin:2013} for related and additional results.

\subsection{The Lie algebra of reflections in types BCD}

This paper extends Marin's Lie algebra calculations to Weyl groups types $B_n$ and $D_n$. Let $\calB_n = W(B_n)$ be the Weyl group of type $B_n$ (the hyperoctahedral group, or the group of signed permutations), and let $\calD_n = W(D_n)$ be the Weyl group of type $D_n$ (the demihyperoctahedral group, or the group of even-signed permutations). Foundational results on the complex representation theories of $\calB_n$ and $\calD_n$ are recalled below in Sections \ref{sec:Type-B-classical} and \ref{sec:Type-D-classical}. In particular, the isomorphism classes of the simple $\CB_n$-modules are labeled by the set $\calBP(n)$ of bipartitions of $n$, i.e., by ordered pairs $(\lambda,\mu)$ of integer partitions such that $\abs{\lambda} + \abs{\mu} = n$. The isomorphism classes of the simple $\CD_n$-modules are labeled by the set $\calBP\{n\}$ of all unordered pairs $\set{\lambda,\mu}$ of integer partitions whose parts together sum to $n$, except that if $n$ is even and $\lambda \vdash n/2$, then the unordered pair $\set{\lambda,\lambda}$ gives rise to two simple $\CD_n$-modules labeled by $\set{\lambda,+}$ and $\set{\lambda,-}$. Theorems \ref{theorem:AW-type-B} and \ref{theorem:AW-type-D} record the descriptions of $\CB_n$ and $\CD_n$, and hence of $\Lie(\CB_{n})$ and $\Lie(\CD_{n})$, as direct sums of full matrix algebras indexed by these sets.

Once again, the more interesting problem is to determine the structure of the Lie subalgebras $\fb_{n} \subset \Lie (\CB_{n})$ and $\fd_{n} \subset \Lie (\CD_{n})$ generated by the sets of reflections in $\calB_n$ and $\calD_n$. The Lie algebras $\fb_n$ and $\fd_n$ turn out to be reductive, with centers spanned by sums of reflections in $\CB_n$ and $\CD_n$; see Section \ref{subsec:generalities} (where reductivity and the description of the center is established under more general hypotheses) and Section \ref{subsec:Lie-alg-gen-reflections} for details. The main problem is then to determine the structure of the semisimple derived subalgebras $\fb_n' = [\fb_n,\fb_n]$ and $\fd_n' = [\fd_n,\fd_n]$. Section \ref{S:LA-by-reflections} culminates in Proposition \ref{prop:restriction-of-AW-map}, which states that the restrictions to $\fb_{n}'$ and $\fd_{n}'$ of the Artin--Wedderburn maps for $\CB_n$ and $\CD_n$ have images in certain direct sums of Lie algebras indexed by subsets of $\calBP(n)$ and $\calBP\{n\}$; see Section \ref{subsec:factorizations} for an explanation of notation. 

\begin{proposition*}[\ref{prop:restriction-of-AW-map}]
For $n \geq 2$, the Artin--Wedderburn maps induce Lie algebra homomorphisms
	\begin{equation} \label{eq:bn-refined-factorization-intro}
	\fb_n' \to \calL_n \oplus \fb_{\beta} \oplus \fb_{\gamma} \oplus \Big[ \bigoplus_{(\lambda,\mu) \in E(n)/\sim} \fb_{(\lambda,\mu)} \Big] \oplus \Big[ \bigoplus_{(\lambda,\mu) \in F(n)} \fb_{(\lambda,\mu)} \Big],
	\end{equation}
and
	\begin{equation} \label{eq:dn-refined-factorization-intro}
	\begin{split}
	\fd_n' \to \calL_n \oplus \fd_{\beta} 
	&\oplus \Big[ \bigoplus_{\substack{\set{\lambda,\mu} \in E\{n\}/\sim \\ \lambda \neq \mu}} \fd_{\set{\lambda,\mu}} \Big] 
	\oplus \Big[ \bigoplus_{\set{\lambda,\lambda} \in E\{n\}/\sim} \fd_{\set{\lambda,+}} \oplus \fd_{\set{\lambda,-}} \Big] \\
	&\oplus \Big[ \bigoplus_{\substack{\set{\lambda,\mu} \in F\{n\} \\ \lambda \neq \mu}} \fd_{\set{\lambda,\mu}} \Big] 
	\oplus \Big[ \bigoplus_{\set{\lambda,\lambda} \in F\{n\}} \fd_{\set{\lambda,+}} \oplus \fd_{\set{\lambda,-}}^\star \Big],
	\end{split}
	\end{equation}
where $\calL_{n}$ is isomorphic to the Lie algebra $\fs_{n}'$ determined by Marin, and where the notation $\star$ means that the term $\fd_{\set{\lambda,-}}^\star$ is omitted if $n/2$ is odd.
\end{proposition*}

The main result of the paper is Theorem \ref{thm:main-theorem}, which states that the maps \eqref{eq:bn-refined-factorization-intro} and \eqref{eq:dn-refined-factorization-intro} are isomorphisms and describes the summands in the image of each map. In particular, each summand is a simple complex Lie algebra of classical type. 

\begin{theorem*}[\ref{thm:main-theorem}]
Let $n \geq 2$. Then the maps \eqref{eq:bn-refined-factorization-intro} and \eqref{eq:dn-refined-factorization-intro} are Lie algebra isomorphisms, and the following identifications hold:
	\begin{enumerate}
	\item $\calL_n = \fsl(\alpha) \oplus [ \bigoplus_{\lambda \in E_n/\sim} \fsl(\lambda) ] \oplus [ \bigoplus_{\lambda \in F_n} \osp(\lambda) ]$ as described in Section \ref{subsubsec:improper-bipartitions}.

	\item $\fd_\beta = \fb_\beta = \fsl(\beta)$ and $\fb_\gamma = \fsl(\gamma)$. For $0 < d < n$, one has $\fd_{\beta_d} = \fb_{\beta_d} \cong \fsl(\beta)$ and $\fb_{\gamma_d} \cong \fsl(\gamma)$.

	\item If $(\lambda,\mu) \in F(n)$, then $\fb_{(\lambda,\mu)} = \osp(\lambda,\mu)$ as described in Lemma \ref{lemma:b-osp-inclusion}.

	\item If $(\lambda,\mu) \in E(n)$, then $\fb_{(\lambda,\mu)} = \fsl(\lambda,\mu)$.

	\item If $\set{\lambda,\mu} \in F\{n\}$ and $\lambda \neq \mu$, then $\fd_{\set{\lambda,\mu}} = \osp\{\lambda,\mu\}$ as described in Lemma \ref{lemma:b-osp-inclusion} or \ref{lemma:d-osp(lam,mu)-inclusion}.

	\item If $\set{\lambda,\lambda} \in F\{n\}$ and $n/2$ is even, then $\fd_{\set{\lambda,\pm}} = \osp\{\lambda,\pm \}$ as described in Lemma \ref{lemma:d-osp(lam,pm)-inclusion}.

	\item If $\set{\lambda,\lambda} \in F\{n\}$ and $n/2$ is odd, then $\fd_{\set{\lambda,\pm}} = \fsl\{\lambda,\pm \}$.

	\item If $\set{\lambda,\mu} \in E\{n\}$ and $\lambda \neq \mu$, then $\fd_{\set{\lambda,\mu}} = \fsl\{\lambda,\mu\}$.

	\item If $\set{\lambda,\lambda} \in E\{n\}$, then $\fd_{\set{\lambda,\pm}} = \fsl\{\lambda,\pm \}$.
	\end{enumerate}
\end{theorem*}

The proof of Theorem \ref{thm:main-theorem} is by induction on $n$, with the values $n \leq 5$ handled via calculations in GAP \cite{GAP4}, and the inductive step for $n \geq 6$ treated in Sections \ref{subsec:item:beta-gamma}--\ref{subsec:complete-main-theorem}. The overall inductive strategy parallels that employed by Marin for type A, but the statements of results and the arguments required in types B and D can be substantially more involved. In particular, types B and D require an intricate, case-by-case analysis of the implications of the branching rules when one restricts simple modules from $\CB_{n}$ to $\CB_{n-1}$ and from $\CD_{n}$ to $\CD_{n-1}$.

\subsection{Further questions}

The results of this paper open the door to investigating analogues for types BCD of the applications pursued by Marin in type A. For example, one can ask for the algebraic hull of the Artin groups of types B and D inside the corresponding characteristic zero Iwahori--Hecke algebras. The image of these Artin groups inside the finite Iwahori--Hecke algebra (defined over a finite field) has been determined by Esterle \cite{Esterle:2018} in the semisimple case. There may also be applications to algebraic combinatorics: the referee pointed out to us that the representation theory of Lie superalgebras in type A has been used (for example) by Stanley to study unimodality and related properties of partition statistics and other combinatorial sequences \cite{Stanley:1985}.

More generally, let $W$ be an arbitrary Weyl group with set of reflections $S$, and let $\fs \subset \Lie(\CW)$ be the (reductive, by  Lemma \ref{lemma:g-reductive}) Lie subalgebra generated by $S$. For the sake of completeness, it would be nice to know the structure of $\fs$ beyond the classical types ABCD. In type $G_2$, the Weyl group $W(G_2)$ is a dihedral group of order $12$, and we leave it as an exercise for the reader to show in that case that $\fs'$ is a direct product of two copies of $\fsl(2)$. For type $F_4$, calculations in GAP show that the Lie algebra $\fs$ is $510$-dimensional, but the Lie algebra's precise structure is currently unknown to the authors.  The structure of $\fs$ for types $E_6$, $E_7$, and $E_8$ also remains open (and, currently, beyond our prowess in GAP to compute the dimensions). One can also ask for the structure of $\fs$ when $W$ is an arbitrary finite Coxeter group, but beyond that scope the generalities in Section \ref{subsec:generalities} begin to break down and it is no longer clear whether to expect nice answers.

If $W$ is a finite Coxeter group, then there exists a canonical supergroup structure on $W$ (i.e., a multi\-plicative grading by $\Z/2\Z$) defined by declaring the Coxeter generators to each be of odd super\-degree. This grading extends to the group algebra $\CW$, making it an associative super\-algebra. One can then consider the induced Lie superalgebra structure on $\CW$, with Lie bracket given by the graded commutator of elements, and ask for the structure of the Lie super\-algebra generated by the set of reflections. In \cite{DK:2023} the authors considered this problem in the case $W = \fS_n$, and showed that the Lie superalgebra generated by the transpositions is equal to the entire derived Lie superalgebra of $\CS_n$ plus the span of an odd central element. Generalizing that work to types BCD (and beyond) would be interesting and nontrivial.

\subsection{Conventions}

Unless indicated otherwise, all vector spaces are over the field $\C$ of complex numbers, all linear maps are $\C$-linear, and all tensor products are over $\C$. Except when indicated by a modifier (e.g., `Lie'), all algebras are assumed to be associative and unital. Given vector spaces $V$ and $W$, let $\dim(V) = \dim_{\C}(V)$ be the dimension of $V$ as a $\C$-vector space, let $\Hom(V,W) = \Hom_{\C}(V,W)$ be the space of all $\C$-linear maps $V \to W$, let $\End(V) = \Hom_{\C}(V,V)$, and let $V^* = \Hom(V,\C)$ be the $\C$-linear dual of $V$. Let $\Z_2 = \set{\pm 1}$ be the cyclic group of order $2$.


Let $\fS_n$ be the symmetric group on $n$ letters. Given a partition $\lambda \vdash n$, let $\lambda^* \vdash n$ be its conjugate (or transpose) partition, let $S^\lambda$ be the simple Specht module for $\CS_n$ labeled by $\lambda$, and let $f^\lambda = \dim(S^\lambda)$. Let $\ve'': \fS_n \to \set{\pm 1}$ be the sign character of $\fS_n$. The trivial $\CS_n$-module is $S^{[n]}$, and the sign module (i.e., the one-dimensional module afforded by $\ve''$) is $S^{[1^n]}$. By abuse of notation, we may also denote $S^{[1^n]}$ by $\ve''$. Then for each $\lambda \vdash n$, one has $S^\lambda \otimes \ve'' \cong S^{\lambda^*}$.

\subsection{Acknowledgements}

The second author thanks the  College of Engineering at Oregon State University for the use of their High Performance Computing Cluster. Calculations on the cluster helped guide the formulation of Theorem \ref{thm:main-theorem}. While the Lie algebra calculations for $n \leq 5$ could be completed on a laptop computer, computations for $n \geq 6$ required additional resources. The second author also thanks Nick Marshall for helpful conversations.

\section{Type B classical representation theory} \label{sec:Type-B-classical}

In Sections \ref{sec:Type-B-classical} and \ref{sec:Type-D-classical} we recall some classical results for Weyl groups of types $B_n$ and $D_n$, many of which can be found in Chapters 1 and 5 of \cite{Geck:2000}.

\subsection{Type B Weyl groups} \label{subsec:Type-B-definition}

Let $\calB_n \subseteq GL_n(\C)$ be the group of signed permutation matrices, i.e., matrices that have precisely one nonzero entry in each row and each column, equal to either $+1$ or $-1$. Let $H_n \subseteq \calB_n$ be the subgroup of permutation matrices in $\calB_n$ (i.e., matrices in $\calB_n$ whose nonzero entries are all $+1$), and let $N_n \subseteq \calB_n$ be the subgroup of diagonal matrices in $\calB_n$. Then $H_n$ identifies with the symmetric group $\fS_n$, $N_n$ identifies with the elementary abelian $2$-group $\Z_2^n$, and $\calB_n = N_n \rtimes H_n$. Henceforward we may make the identification $\calB_n = N_n \rtimes \fS_n$.

For $1 \leq i \leq n-1$, let $s_i \in H_n \subseteq \calB_n$ be the matrix obtained from the identity matrix by interchanging the $i$-th and $(i+1)$-th rows, so $s_i$ corresponds to the transposition $(i,i+1) \in \fS_n$. For $1 \leq i \leq n$, let $t_i \in N_n \subseteq \calB_n$ be the diagonal matrix whose $i$-th diagonal entry is equal to $-1$ and whose other diagonal entries are equal to $+1$, so $t_i$ corresponds to a multiplicative generator for the $i$-th factor in $\Z_2^n$. Then for $n \geq 2$, the pair $(\calB_n,\set{s_1,\ldots,s_{n-1},t_n})$ is a Coxeter system of type $B_n$. Each element of $\calB_n$ is uniquely of the form $t_1^{i_1} \cdots t_n^{i_n} \sigma$ with $\sigma \in H_n$ and $i_j \in \set{0,1}$ for each $j$. Identifying $H_n$ with $\fS_n$, the set of reflections in $\calB_n$ is
	\begin{equation} \label{eq:B-reflections}
	\set{ (i,j), \, t_i t_j (i,j) : 1 \leq i < j \leq n} \cup \set{ t_i : 1 \leq i \leq n}.
	\end{equation}
The two sets in this union are each conjugacy classes in $\calB_n$. Following \cite[\S2]{Okada:1990} or \cite[\S2]{Mishra:2016}, we may label these conjugacy classes as $C_{(21^{n-2},\emptyset)}$ and $C_{(1^{n-1},1)}$, respectively.

\subsection{Simple modules for Type B Weyl groups} \label{subsec:simples-B}

Let $\ve$, $\ve'$, and $\ve''$ be the group homomorphisms $\calB_n \to \set{\pm 1}$ that are defined on generators by
	\begin{equation} \label{eq:Bn-linear-characters}
	\begin{aligned}
	\ve(t_i) &= -1, \qquad & \ve'(t_i) &= -1, \qquad & \ve''(t_i) &= +1, \\
	\ve(s_j) &= -1, \qquad & \ve'(s_j) &= +1, \qquad & \ve''(s_j) &= -1.
	\end{aligned}
	\end{equation}
The map $\ve$ is the sign character of $\calB_n$, i.e., the character that evaluates to $-1$ on each of the Coxeter generators of $\calB_n$. The map $\ve''$ factors via the canonical quotient map $\calB_n = N_n \rtimes H_n \twoheadrightarrow H_n \cong \fS_n$ to the usual the sign character on $\fS_n$. By abuse of notation, we identify $\ve$, $\ve'$, and $\ve''$ with their corresponding one-dimensional $\CB_n$-modules. Then $\alpha = \beta \otimes \gamma$ whenever $\{ \alpha,\beta,\gamma \} = \{ \ve,\ve',\ve'' \}$.

The isomorphism classes of simple $\CB_n$-modules are labeled by bipartitions of $n$, i.e., ordered pairs of partitions $(\lambda,\mu)$ such that $\abs{\lambda}+\abs{\mu} = n$. Let $\calBP(n)$ denote the set of all bipartitions of $n$. Given $(\lambda,\mu) \in \calBP(n)$ with $\abs{\lambda} = a$ and $\abs{\mu} = b$, the corresponding simple $\CB_n$-module $S^{(\lambda,\mu)}$ can be constructed as follows; cf.\ \cite[\S5.5.4]{Geck:2000}. First, the diagonal embedding $GL_a(\C) \times GL_b(\C) \subseteq GL_n(\C)$ induces an embedding $\calB_a \times \calB_b \subseteq \calB_n$. Next, the Specht modules $S^\lambda$ and $S^\mu$ can be inflated to $\calB_a$ and $\calB_b$ via the canonical quotient maps $\calB_a \twoheadrightarrow \fS_a$ and $\calB_b \twoheadrightarrow \fS_b$. Then
	\begin{equation} \label{eq:type-B-module-construction}
	S^{(\lambda,\mu)} = \Ind_{\calB_a \times \calB_b}^{\calB_n}\left( S^\lambda \boxtimes (\ve' \otimes S^\mu) \right),
	\end{equation}
where $\Ind$ denotes the usual tensor induction functor and $\boxtimes$ denotes the external tensor product of modules. Then the trivial $\CB_n$-module is $S^{([n],\emptyset)}$, and identifying a linear character with its corresponding module, one has $\ve = S^{(\emptyset,[1^n])}$, $\ve' = S^{(\emptyset,[n])}$, and $\ve'' = S^{([1^n],\emptyset)}$. More generally, if $\lambda \vdash n$, then $S^{(\lambda,\emptyset)} = \Inf_{\fS_n}^{\calB_n}(S^\lambda)$, the inflation via the canonical quotient map $\calB_n \twoheadrightarrow \fS_n$ of the $\fS_n$-module $S^\lambda$, and $S^{(\emptyset,\lambda)} = \Inf_{\fS_n}^{\calB_n}(S^\lambda) \otimes \ve'$. As modules over the subalgebra $\CS_n \subset \CB_n$, one has $S^{(\lambda,\emptyset)} = S^\lambda = S^{(\emptyset,\lambda)}$. From \eqref{eq:type-B-module-construction}, one gets
	\begin{equation} \label{eq:dim-S(lambda,mu)}
	\dim(S^{(\lambda,\mu)}) = \frac{\abs{\calB_n}}{{\abs{\calB_a}} \cdot \abs{\calB_b}} \cdot \dim(S^\lambda) \cdot \dim(S^\mu) = \binom{n}{a} \cdot \dim(S^\lambda) \cdot \dim(S^\mu).
	\end{equation}

\begin{lemma} \label{lemma:Bn-subminimal-dimensions} \ 
\begin{enumerate}
\item The only one-dimensional $\CB_n$-modules are those labeled by the bipartitions $([n],\emptyset)$, $(\emptyset,[1^n])$, $(\emptyset,[n])$, and $([1^n],\emptyset)$, i.e., the trivial module and the modules afforded by $\ve$, $\ve'$, and $\ve''$.

\item For $n \geq 5$, if $(\lambda,\mu) \in \calBP(n)$ and $\dim(S^{(\lambda,\mu)}) \neq 1$, then $\dim(S^{(\lambda,\mu)}) \geq n-1$.

\item If $\lambda \vdash n/2$, then $\dim(S^{(\lambda,\lambda^*)}) = \dim(S^{(\lambda,\lambda)}) = \binom{n}{n/2} \cdot \dim(S^\lambda)^2 \geq \binom{n}{n/2}$.
\end{enumerate}
\end{lemma}

\begin{proof}
If $0 < k < n$, then $\binom{n}{k} \geq n$. Now apply \eqref{eq:dim-S(lambda,mu)} and \cite[Theorem 2.4.10]{James:1981}.
\end{proof}

For $(\lambda,\mu) \in \calBP(n)$, let $\rho_{(\lambda,\mu)}: \CB_n \to \End(S^{(\lambda,\mu)})$ denote the $\CB_n$-module structure map for $S^{(\lambda,\mu)}$. Then the Artin--Wedderburn Theorem implies:

\begin{theorem} \label{theorem:AW-type-B}
The map $\sigma \mapsto \bigoplus_{(\lambda,\mu) \in \calBP(n)} \rho_{(\lambda,\mu)}(\sigma)$ induces an algebra isomorphism
	\begin{equation} \label{eq:AW-type-B}
	\CB_n \cong \bigoplus_{(\lambda,\mu) \in \calBP(n)} \End(S^{(\lambda,\mu)}).
	\end{equation}
\end{theorem}

By \cite[Theorem 5.5.6(c)]{Geck:2000}, one has
	\begin{equation} \label{eq:B-tensor-with-character}
	\begin{aligned}
	S^{(\lambda,\mu)} \otimes \ve &\cong S^{(\mu^*,\lambda^*)}, \qquad &
	S^{(\lambda,\mu)} \otimes \ve' &\cong S^{(\mu,\lambda)}, \qquad &
	S^{(\lambda,\mu)} \otimes \ve'' &\cong S^{(\lambda^*,\mu^*)}.
	\end{aligned}
	\end{equation}
By \cite[Theorem 5.5.6(b)]{Geck:2000}, the field $\Q$ of rational numbers is a splitting field for $\calB_n$. We may write $S_{\R}^{(\lambda,\mu)}$ for the simple $\RB_n$-module whose scalar extension to $\C$ is $S^{(\lambda,\mu)}$. Since $\Q$ is a splitting field for $\calB_n$, it follows that all irreducible characters of $\calB_n$ are real valued, and hence all simple $\CB_n$-modules are self-dual.

\subsection{Weight space decompositions of Type B simple modules} \label{subsec:weight-space-type-B}

Our primary reference for this section is \cite{Mishra:2016}. The Young--Jucys--Murphy (YJM) elements $X_1,\ldots,X_n \in \CB_n$ are defined by
	\[
	X_i = \sum_{k=1}^{i-1} \left[ (k,i) + t_k t_i (k,i) \right].
	\]
The Gelfand--Tsetlin (GZ) algebra $GZ(\calB_n)$ is the commutative $\C$-sub\-algebra of $\CB_n$ generated by $t_1,\ldots,t_n$ and $X_1,\ldots,X_n$. Each simple $\CB_n$-module $V$ decomposes into one-dimen\-sional $GZ(\calB_n)$-submodules, called the GZ-subspaces of $V$. Each GZ-subspace $W$ is determined by the irreducible representation $\rho$ by which $\Z_2^n = \subgrp{t_1,\ldots,t_n}$ acts on $W$, and by the list $\chi = (\chi_1,\ldots,\chi_n)$ of eigen\-values by which the operators $X_1,\ldots,X_n$ act on $W$. We label the irreducible representations of $\Z_2^n$ by $\rho = (\rho_1,\ldots,\rho_n) \in \set{\pm 1}^n$, such that $\rho(t_i) = \rho_i$. The tuple
	\[
	\alpha = (\rho,\chi) = (\rho_1,\ldots,\rho_n,\chi_1,\ldots,\chi_n)
	\]
is called the weight of the GZ-subspace $W$. Given a simple $\CB_n$-module $V$, we may denote its decomposition into GZ-subspaces by $V = \bigoplus_{\alpha} V_{\alpha}$, where $\alpha$ is the weight of $V_{\alpha}$.

Given a partition $\lambda = (\lambda_1 \geq \lambda_2 \geq \cdots)$ of $n$, we draw the Young diagram of shape $\lambda$ via the ``English'' convention (see \cite{wiki-Young-tableau}), as an array of boxes with $\lambda_i$ boxes in the $i$-th row, the $(i+1)$-th row drawn below the $i$-th row, and the rows of boxes aligned on the left. A box in the $i$-th row and $j$-th column (from the left) of the diagram is said to have \emph{residue} $j-i$. A $\lambda$-tableau is a Young diagram of shape $\lambda$ whose boxes have been filled in some order with $n$ different positive integers. A \emph{standard} $\lambda$-tableau is a $\lambda$-tableau in which the values of the integers increase from top to bottom along columns, and from left to right along rows.

Given $(\lambda,\mu) \in \calBP(n)$, let $\T(\lambda,\mu)$ be the set of all standard bitableaux of shape $(\lambda,\mu)$, i.e., the set of all ordered pairs $T = (T_{+1},T_{-1})$ such that $T_{+1}$ and $T_{-1}$ are standard tableaux of shapes $\lambda$ and $\mu$, respectively, whose boxes have collectively been filled by the integers $1,2,\ldots,n$. The GZ-subspaces of $S^{(\lambda,\mu)}$ are indexed by the elements of $\T(\lambda,\mu)$. For $T = (T_{+1},T_{-1}) \in \T(\lambda,\mu)$, the corresponding weight $\alpha(T)$ is defined by
	\begin{equation} \label{eq:alpha(T)}
	\alpha(T) = \left(\rho_T(1),\ldots,\rho_T(n), 2 \cdot \res_T(1), \ldots, 2 \cdot \res_T(n) \right),
	\end{equation}
where $\rho_T(i) \in \set{\pm 1}$ is the subscript of the tableau ($T_{+1}$ or $T_{-1}$) in which the integer $i$ is located, and $\res_T(i)$ is the residue (or content, in the terminology of \cite{Mishra:2016}) of the box (in either $T_{+1}$ or $T_{-1}$) in which the integer $i$ is located. Then for each $(\lambda,\mu) \in \calBP(n)$, one has (cf.\ \cite[Theorem 6.5]{Mishra:2016})
	\[
	S^{(\lambda,\mu)} = \bigoplus_{T \in \T(\lambda,\mu)} S^{(\lambda,\mu)}_T = \bigoplus_{T \in \T(\lambda,\mu)} S^{(\lambda,\mu)}_{\alpha(T)}.
	\]
For $(\lambda,\mu) \in \calBP(n)$, set $\calW(\lambda,\mu) = \set{ \alpha(T) : T \in \T(\lambda,\mu)}$. Set $\T(n) = \bigcup_{(\lambda,\mu) \in \calBP(n)} \T(\lambda,\mu)$, and set $\calW(n) = \bigcup_{(\lambda,\mu) \in \calBP(n)} \calW(\lambda,\mu)$. Then the map $T \mapsto \alpha(T)$ is a bijection $\T(n) \to \calW(n)$.

Given a standard $\lambda$-tableau $T$, let $T^*$ be its transpose, which is then a standard $\lambda^*$-tableau. For $T = (T_{+1},T_{-1}) \in \T(\lambda,\mu)$, set $T^* = (T_{+1}^*,T_{-1}^*) \in \T(\lambda^*,\mu^*)$ and $T^\natural = (T_{-1},T_{+1}) \in \T(\mu,\lambda)$. Then if $T \in \T(\lambda,\mu)$ with $\alpha(T) = (\rho,\chi)$, one has $T^{*\natural} = T^{\natural*}$, and
	\begin{gather}
	\alpha(T^*) = (\rho,-\chi) \in \calW(\lambda^*,\mu^*), \label{eq:weights-transpose} \\
	\alpha(T^\natural) = (-\rho,\chi) \in \calW(\mu,\lambda), \label{eq:weights-swap} \\
	\alpha(T^{*\natural}) = -\alpha(T) \in \calW(\mu^*,\lambda^*). \label{eq:weights-swap-transpose}
	\end{gather}

\subsection{Orthogonal forms} \label{subsec:orthogonal-form-B}

The next theorem is a consequence of \cite[Theorem 6.12]{Mishra:2016}.

\begin{theorem} \label{thm:type-B-normal-form-action}
Let $(\lambda,\mu) \in \calBP(n)$. Then with respect to a Hermitian, conjugate-symmetric, $\calB_n$-invariant inner product on $S^{(\lambda,\mu)}$, there exists an orthonormal basis
	\begin{equation} \label{eq:B(lam,mu)}
	B^{(\lambda,\mu)} = \set{ c_T : T \in \T(\lambda,\mu)}
	\end{equation}
for $S^{(\lambda,\mu)}$, such that $c_T \in S_T^{(\lambda,\mu)}$, and such that for $1 \leq i < n$, the action of the Coxeter generator $s_i = (i,i+1) \in \fS_n \subset \calB_n$ is as follows:

Given $T \in \T(\lambda,\mu)$ with $\alpha(T) = (\rho,\chi)$, set $r_i = \frac{1}{2}(\chi_{i+1}-\chi_i) = \res_T(i+1) - \res_T(i)$.
	\begin{enumerate}
	\item Suppose $i$ and $i+1$ are not in the same tableau in $T$. Let $S = s_i \cdot T$. Then $S_T^{(\lambda,\mu)} \oplus S_S^{(\lambda,\mu)}$ is closed under the action of $s_i$, and the matrix of $s_i$ with respect to the basis $\set{c_T,c_S}$ is
		\[
		\begin{bmatrix}
		0 & 1 \\ 
		1 & 0 
		\end{bmatrix}.
		\]
	
	\item If $i$ and $i+1$ are in the same column of the same tableau in $T$, implying that $r_i = -1$, then $s_i$ acts on $S_T^{(\lambda,\mu)}$ as multiplication by $-1$.
	
	\item If $i$ and $i+1$ are in the same row of the same tableau in $T$, implying that $r_i = +1$, then $s_i$ acts on $S_T^{(\lambda,\mu)}$ as multiplication by $+1$.
	
	\item Suppose $i$ and $i+1$ are in the same tableau in $T$, but not in the same row or in the same column, implying that $\abs{r_i} \geq 2$. Let $S = s_i \cdot T$. Then $S_T^{(\lambda,\mu)} \oplus S_S^{(\lambda,\mu)}$ is closed under the action of $s_i$, and the matrix of $s_i$ with respect to the basis $\set{c_T,c_S}$ is
		\[
		\begin{bmatrix}
		r_i^{-1} & \sqrt{1-r_i^{-2}} \\
		\sqrt{1-r_i^{-2}} & -r_i^{-1}
		\end{bmatrix}.
		\]
	\end{enumerate}
\end{theorem}

\begin{remark} \label{remark:C-bilinear-form}
The simple $\RB_n$-module $S_{\R}^{(\lambda,\mu)}$ can be taken to be the $\R$-span in $S^{(\lambda,\mu)}$ of the set of vectors $\set{ c_T: T \in \T(\lambda,\mu)}$. The inner product on $S^{(\lambda,\mu)}$ then restricts to a symmetric $\calB_n$-invariant inner product $\subgrp{-,-}_{\R}$ on $S_{\R}^{(\lambda,\mu)}$ with respect to which \eqref{eq:B(lam,mu)} is an orthonormal basis. Extending scalars linearly in each variable, $\subgrp{-,-}_{\R}$ then extends to a symmetric $\calB_n$-invariant $\C$-bilinear form $\subgrp{-,-}_{(\lambda,\mu)}$ on $S^{(\lambda,\mu)}$ with respect to which \eqref{eq:B(lam,mu)} is an ortho\-normal basis.
\end{remark}

Recall the group homomorphism $\ve': \calB_n \to \set{\pm 1}$ defined in \eqref{eq:Bn-linear-characters}.

\begin{lemma} \label{lemma:epsilon'-intertwinor}
Given $(\lambda,\mu) \in \calBP(n)$, let $\phi_{\ve'}^{(\lambda,\mu)}: S^{(\lambda,\mu)} \to S^{(\mu,\lambda)}$ be the linear map defined in terms of the bases $B^{(\lambda,\mu)}$ and $B^{(\mu,\lambda)}$ of Theorem \ref{thm:type-B-normal-form-action} by
	\[
	\phi_{\ve'}^{(\lambda,\mu)}(c_T) = c_{T^{\natural}}.
	\]
Then $\phi_{\ve'}^{(\mu,\lambda)} \circ \phi_{\ve'}^{(\lambda,\mu)} = \Id_{S^{(\lambda,\mu)}}$, and for all $g \in \calB_n$ and all $v \in S^{(\lambda,\mu)}$, one has
	\begin{equation} \label{eq:epsilon'-intertwinor}
	\phi_{\ve'}^{(\lambda,\mu)}(g.v) = \ve'(g)g \cdot \phi_{\ve'}^{(\lambda,\mu)}(v).
	\end{equation}
\end{lemma}

\begin{proof}
The equality $\phi_{\ve'}^{(\mu,\lambda)} \circ \phi_{\ve'}^{(\lambda,\mu)} = \Id_{S^{(\lambda,\mu)}}$ is immediate. To prove \eqref{eq:epsilon'-intertwinor}, it suffices by linearity to assume that $v \in B^{(\lambda,\mu)}$, and then it suffices by multiplicativity to assume that $g$ is one of the generators $t_1,\ldots,t_n$ or $s_1,\ldots,s_{n-1}$ of $\calB_n$. The equality can then be checked using the description of weights in \eqref{eq:alpha(T)} and by cases using Theorem \ref{thm:type-B-normal-form-action}.
\end{proof}

Given $(\lambda,\mu) \in \calBP(n)$, let $R^{(\lambda,\mu)} = (R^{(\lambda,\mu)}_{+1},R^{(\lambda,\mu)}_{-1})$ be the ``row major'' standard bitableau of shape $(\lambda,\mu)$, in which $R^{(\lambda,\mu)}_{+1}$ (resp.\ $R^{(\lambda,\mu)}_{-1}$) is filled with the integers $1,\ldots,\abs{\lambda}$ (resp.\ $\abs{\lambda}+1,\ldots,n$) in row major order (i.e., $R^{(\lambda,\mu)}_{+1}$ contains $1,2,\ldots,\lambda_1$ in its first row, then $\lambda_1+1,\ldots,\lambda_2$ in its second row, etc., and similarly for $R^{(\lambda,\mu)}_{-1}$). Given $T \in \T(\lambda,\mu)$, there exists a unique permutation $\sigma_T \in \fS_n$ such that $\sigma_T \cdot R^{(\lambda,\mu)} = T$. Define the \emph{length} of $T$ by $\ell(T) = \ell(\sigma_T)$, the length of $\sigma_T$ as an element of the Coxeter group $\fS_n$. Then $(-1)^{\ell(T)} = (-1)^{\ell(\sigma_T)} = \ve''(\sigma_T)$. If $\tau \in \fS_n$, then
	\[
	(-1)^{\ell(\tau \cdot T)} = (-1)^{\ell(\tau \cdot \sigma_T)} = \ve''(\tau) \cdot \ve''(\sigma_T).
	\]

\begin{lemma} \label{lemma:epsilon''-intertwinor}
Given $(\lambda,\mu) \in \calBP(n)$, let $\phi_{\ve''}^{(\lambda,\mu)}: S^{(\lambda,\mu)} \to S^{(\lambda^*,\mu^*)}$ be the linear map defined in terms of the bases $B^{(\lambda,\mu)}$ and $B^{(\lambda^*,\mu^*)}$ of Theorem \ref{thm:type-B-normal-form-action} by
	\[
	\phi_{\ve''}^{(\lambda,\mu)}(c_T) = (-1)^{\ell(T)} \cdot c_{T^*}.
	\]
Then for all $g \in \calB_n$ and all $v \in S^{(\lambda,\mu)}$, one has
	\begin{equation} \label{eq:epsilon''-intertwinor}
	\phi_{\ve''}^{(\lambda,\mu)}(g.v) = \ve''(g)g \cdot \phi_{\ve''}^{(\lambda,\mu)}(v).
	\end{equation}
\end{lemma}

\begin{proof}
The strategy is the same as for the proof of Lemma \ref{lemma:epsilon'-intertwinor}.
\end{proof}

\begin{remark} \label{remark:associators-define-isos}
The linear maps $\phi_{\ve'}^{(\lambda,\mu)}$, $\phi_{\ve''}^{(\lambda,\mu)}$, and $\phi_{\ve}^{(\lambda,\mu)} := \phi_{\ve''}^{(\mu,\lambda)} \circ \phi_{\ve'}^{(\lambda,\mu)}$ can be interpreted as defining the isomorphisms $S^{(\lambda,\mu)} \otimes \ve' \cong S^{(\mu,\lambda)}$, $S^{(\lambda,\mu)} \otimes \ve'' \cong S^{(\lambda^*,\mu^*)}$, and $S^{(\lambda,\mu)} \otimes \ve \cong S^{(\mu^*,\lambda^*)}$.
\end{remark}

Say that the action of an adjacent transposition $s_i = (i,i+1)$ on a standard bitableau $T \in \T(\lambda,\mu)$ is \emph{admissible} if $s_i \cdot T$ is again a standard bitableau, i.e., if the action corresponds to the first or fourth cases of Lemma \ref{thm:type-B-normal-form-action}. The proof of \cite[Lemma 6.2]{Mishra:2016} shows that $T$ can be obtained from the row major standard bitableau $R^{(\lambda,\mu)}$ by applying a sequence of $\ell(T)$ admissible transpositions, i.e., such that at each intermediate step the resulting bitableau is standard. It follows for arbitrary $T,S \in \T(\lambda,\mu)$ that the permutation $\sigma_{T,S} \in \fS_n$ that maps $T$ to $S$ can be written as a product of admissible transpositions; let $\ell_{T,S}$ be the minimal number of terms in such a factorization.

\begin{lemma} \label{lemma:admissible-linear-combination}
Let $(\lambda,\mu) \in \calBP(n)$, let $T \in \T(\lambda,\mu)$, and let $\sigma \in \fS_n$ such that $\sigma \cdot T \in \T(\lambda,\mu)$. Then there exist scalars $\gamma_{T,S}^\sigma \in \C$ for $S \in \T(\lambda,\mu)$, with $\gamma_{T,\sigma \cdot T}^\sigma \neq 0$, such that
	\[
	\sigma \cdot c_T = \gamma_{T,\sigma \cdot T}^\sigma \cdot c_{\sigma \cdot T} +\sum_{\substack{S \in \T(\lambda,\mu) \\ \ell_{T,S} < \ell_{T,\sigma \cdot T}}} \gamma_{T,S}^\sigma \cdot c_S.
	\]
\end{lemma}

\begin{proof}
If $\ell_{T,\sigma \cdot T} \leq 1$, then the result is immediate from Theorem \ref{thm:type-B-normal-form-action}, so suppose that $\ell_{T,\sigma \cdot T} > 1$, let $\sigma = \sigma_1 \sigma_2 \cdots \sigma_\ell$ be a factorization of $\sigma$ into a product of $\ell = \ell_{T,\sigma \cdot T}$ admissible (adjacent) transpositions, and let $\tau = \sigma_2 \cdots \sigma_\ell$. Then $\ell_{T,\tau \cdot T} = \ell-1$, and by induction one has
	\[
	\tau \cdot c_T = \gamma_{T,\tau \cdot T}^\tau \cdot c_{\tau \cdot T} +\sum_{\substack{S \in \T(\lambda,\mu) \\ \ell_{T,S} < \ell_{T,\tau \cdot T}}} \gamma_{T,S}^\tau \cdot c_S.
	\]
with $\gamma_{T,\tau \cdot T}^\tau \neq 0$. Then
	\[
	\sigma \cdot c_T = \gamma_{T,\tau \cdot T}^\tau \cdot (\sigma_1 \cdot c_{\tau \cdot T}) +\sum_{\substack{S \in \T(\lambda,\mu) \\ \ell_{T,S} < \ell_{T,\tau \cdot T}}} \gamma_{T,S}^\tau \cdot (\sigma_1 \cdot c_S).
	\]
Now for $S \in \T(\lambda,\mu)$, either $\sigma_1 \cdot S \notin \T(\lambda,\mu)$, in which case $\sigma_1 \cdot c_S = \pm c_S$, or else $\sigma_1 \cdot S \in \T(\lambda,\mu)$ and $\ell_{T,\sigma_1 \cdot S} \leq \ell_{T,S}+1$. Then the result follows by applying Theorem \ref{thm:type-B-normal-form-action}.
\end{proof}

\subsection{Restriction of Type B simple modules} \label{subsec:branching-type-B}

Given partitions $\lambda \vdash n$ and $\mu \vdash (n-1)$, write $\mu \prec \lambda$ if the Young diagram of $\mu$ can be obtained by removing a box from the Young diagram of $\lambda$. In this case, let $\res(\lambda/\mu)$ denote the residue of the box that is removed from $\lambda$ to obtain $\mu$. Let $r_\lambda = \abs{\set{ \mu \vdash (n-1) : \mu \prec \lambda}}$ be the number of removable boxes in the Young diagram of $\lambda$. For $(\lambda,\mu) \in \calBP(n)$, write $(\nu,\tau) \prec (\lambda,\mu)$ if either $\nu \prec \lambda$ and $\tau = \mu$, or if $\nu = \lambda$ and $\tau \prec \nu$.

Identify $\calB_{n-1}$ with the subgroup of $\calB_n$ generated by $\set{s_1,\ldots,s_{n-2},t_{n-1}}$, and write $\Res_{\calB_{n-1}}^{\calB_n}(V)$ for the restriction of a $\CB_n$-module $V$ to the subalgebra $\CB_{n-1}$. Then for $(\lambda,\mu) \in \calBP(n)$, one has the multiplicity-free (hence canonical, by uniqueness of isotypical components) decomposition
	\begin{equation} \label{eq:restriction-Bn-toBn-1}
	\Res_{\calB_{n-1}}^{\calB_n}(S^{(\lambda,\mu)}) \cong \bigoplus_{(\nu,\tau) \prec (\lambda,\mu)} S^{(\nu,\tau)} = \Big[ \bigoplus_{\nu \prec \lambda} S^{(\nu,\mu)} \Big] \oplus \Big[ \bigoplus_{\tau \prec \mu} S^{(\lambda,\tau)} \Big];
	\end{equation}
see \cite[Theorem 4.1]{Okada:1990} or \cite[Theorem 6.6]{Mishra:2016}. In terms of weight spaces, \eqref{eq:restriction-Bn-toBn-1} takes the form:
	\begin{equation} \label{eq:B-to-Bn-1-weights}
	\Res_{\calB_{n-1}}^{\calB_n}(S^{(\lambda,\mu)}) \cong 
	\bigg[ \bigoplus_{\nu \prec \lambda} \bigoplus_{\substack{\alpha \in \calW(\lambda,\mu) \\ \rho_n = +1 \\ \chi_n = 2 \cdot \res(\lambda/\nu)}} S^{(\lambda,\mu)}_\alpha \bigg] \oplus
	\bigg[ \bigoplus_{\tau \prec \mu} \bigoplus_{\substack{\alpha \in \calW(\lambda,\mu) \\ \rho_n = -1 \\ \chi_n = 2 \cdot \res(\mu/\tau)}} S^{(\lambda,\mu)}_\alpha \bigg].
	\end{equation}
The summand indexed by $\nu$ (resp.\ $\tau$) is isomorphic as a $\CB_{n-1}$-module to $S^{(\nu,\mu)}$ (resp.\ $S^{(\lambda,\tau)}$).

\section{Type D classical representation theory} \label{sec:Type-D-classical}

\subsection{Type D Weyl groups}

Recall the linear characters $\ve'$ and $\ve''$ on $\calB_n$ defined in \eqref{eq:Bn-linear-characters}. Set $\calD_n = \ker(\ve')$ and $N_n' = N_n \cap \calD_n$. Let $\wt{s}_n = t_n s_{n-1} t_n = t_{n-1} t_n s_{n-1} \in \calD_n$. Then for $n \geq 4$, the pair $(\calD_n,\set{s_1,\ldots,s_{n-1},\wt{s}_n})$ is a Coxeter system of type $D_n$, and $\ve''$ restricts to the sign character of $\calD_n$. One has $\calD_n = N_n' \rtimes H_n$, and each element of $\calD_n$ is uniquely of the form $t_1^{i_1} \cdots t_n^{i_n} \sigma$ with $\sigma \in H_n$, $i_j \in \set{0,1}$, and $i_1 + \cdots + i_n$ even. For $1 \leq i \leq n-1$, let $u_i = t_1 t_{i+1}$. Then each element of $\calD_n$ is also uniquely of the form $u_1^{i_1} \cdots u_{n-1}^{i_{n-1}} \sigma$ with $\sigma \in H_n$ and $i_j \in \set{0,1}$ for each $j$. One has $N_n' \cong \Z_2^{n-1}$; via this identification, $u_i$ is a multiplicative generator for the $i$-th factor in $\Z_2^{n-1}$. Identifying $H_n$ with $\fS_n$, the set of reflections in $\calD_n$ is
	\begin{equation} \label{eq:D-reflections}
	\set{ (i,j), \, t_i t_j (i,j) : 1 \leq i < j \leq n}.
	\end{equation}
For $n \geq 3$, this set is a conjugacy class in $\calD_n$; for $n=2$, it is a union of two conjugacy classes.

\subsection{Simple modules for Type D Weyl groups}

Suppose $H$ is an index-$2$ subgroup of a finite group $G$. Given $t \in G \backslash H$ and a $\C H$-module $W$, let ${}^t W$ be the $\C H$-module obtained from $W$ by twisting the $H$-action via the automorphism $h \mapsto tht^{-1}$. Up to isomorphism, ${}^t W$ is independent of the choice of $t \in G \backslash H$. Recall that two $\C H$-modules $W$ and $W'$ are said to be \emph{conjugate} if $W' \cong {}^t W$ for some $t \in G \backslash H$. In particular, conjugate modules are of the same dimension.

For each $(\lambda,\mu) \in \calBP(n)$, the map $\phi_{\ve'}^{(\lambda,\mu)}: S^{(\lambda,\mu)} \to S^{(\mu,\lambda)}$ of Lemma \ref{lemma:epsilon'-intertwinor} defines a $\CD_n$-module isomorphism $\Res^{\calB_n}_{\calD_n}(S^{(\lambda,\mu)}) \cong \Res^{\calB_n}_{\calD_n}(S^{(\mu,\lambda)})$; in this way we identify these $\CD_n$-modules. Let
	\begin{equation} \label{eq:BP-set-n}
	\calBP\{n\} = \set{ \set{\lambda,\mu} : (\lambda,\mu) \in \calBP(n)}
	\end{equation}
be the set of all unordered pairs of partitions whose parts together sum to $n$. The simple $\CD_n$-modules are now obtained as follows; cf.\ \cite[Proposition 5.1]{Fulton:1991}.

\begin{lemma} \label{lemma:Dn-simples}
Up to isomorphism, each simple $\CD_n$-module arises uniquely via either:
	\begin{enumerate}
	\item If $\set{\lambda,\mu} \in \calBP\{n\}$ and $\lambda \neq \mu$, then $S^{\set{\lambda,\mu}} \coloneq \Res^{\calB_n}_{\calD_n}( S^{(\lambda,\mu)}) \cong \Res^{\calB_n}_{\calD_n}( S^{(\mu,\lambda)})$ is simple and isomorphic to its conjugate.
	
	\item If $\set{\lambda,\lambda} \in \calBP\{n\}$, then $S^{\set{\lambda,\lambda}} \coloneq \Res^{\calB_n}_{\calD_n}(S^{(\lambda,\lambda)})$ is canonically the direct sum of two simple, conjugate, non-isomorphic $\CD_n$-modules, which we denote $S^{\set{\lambda,+}}$ and $S^{\set{\lambda,-}}$.
	\end{enumerate}
\end{lemma}

By \cite[Corollary 5.6.4]{Geck:2000}, the field $\Q$ is a splitting field for $\calD_n$. Then all irreducible characters of $\calD_n$ are real valued, and hence all simple $\CD_n$-modules are self-dual.

\begin{remark} \label{remark:S(lam-lam)-decomposition}
Suppose $(\lambda,\lambda) \in \calBP(n)$. The $\CD_n$-module decomposition of $S^{(\lambda,\lambda)}$ is canonical, by the uniqueness of isotypical components. Using Lemma \ref{lemma:epsilon'-intertwinor}, one sees that the $\CD_n$-module summands of $S^{(\lambda,\lambda)}$ can be realized explicitly as
	\begin{align*}
	S^{\set{\lambda,+}} &= \Span\set{ c_{T} + c_{T^\natural} : T \in \T(\lambda,\lambda)}, \\
	S^{\set{\lambda,-}} &= \Span\set{ c_{T} - c_{T^\natural} : T \in \T(\lambda,\lambda)}.
	\end{align*}
For a given value of $n$, this labeling of $S^{\set{\lambda,\pm}}$ may be the opposite of that corresponding to the characters $\chi_{(\lambda,\pm)}$ defined in \cite[\S 5.6]{Geck:2000} or \cite[\S3]{Geck:2015}. For $T \in \T(\lambda,\lambda)$, set $S_T^{\set{\lambda,\pm}} = \Span\{ c_T \pm c_{T^\natural} \}$. Define subsets $\T(\lambda,\lambda)_{\pm 1}$ of $\T(\lambda,\lambda)$ by $\T(\lambda,\lambda)_{\pm 1} = \{ T \in \T(\lambda,\lambda) : \rho_T(1) = \pm 1 \}$. Then $\T(\lambda,\lambda)$ is the disjoint union of the sets $\T(\lambda,\lambda)_{+1}$ and $\T(\lambda,\lambda)_{-1}$, and the sets
	\begin{equation} \label{eq:(lam,pm)-bases}
	\set{ c_T + c_{T^\natural} : T \in \T(\lambda,\lambda)_{+ 1}} \quad \text{and} \quad \set{ c_T - c_{T^\natural} : T \in \T(\lambda,\lambda)_{+ 1}}
	\end{equation}
are bases for $S^{\set{\lambda,+}}$ and $S^{\set{\lambda,-}}$, respectively, with $\{ c_T \pm c_{T^\natural}\}$ a basis for $S^{\set{\lambda,\pm}}_T$. The bilinear form on $S^{(\lambda,\lambda)}$ described in Remark \ref{remark:C-bilinear-form} restricts to a bilinear form on $S^{\set{\lambda,\pm}}$, with respect to which the set $\set{ c_T \pm c_{T^\natural} : T \in \T(\lambda,\lambda)_{+ 1}}$ is an orthogonal basis.
\end{remark}

\begin{remark} \label{remark:D-tensor-with-sign}
The third isomorphism in \eqref{eq:B-tensor-with-character} restricts to a $\CD_n$-module isomorphism
	\[
	S^{\set{\lambda,\mu}} \otimes \ve'' \cong S^{\set{\lambda^*,\mu^*}}.
	\]
Suppose $(\lambda,\lambda) \in \calBP(n)$, so in particular, $n$ is even. For $T \in \T(\lambda,\lambda)$, one has $T^{\natural} = \tau \cdot T$, where $\tau$ is a product of $n/2$ disjoint transpositions. Then $(-1)^{\ell(T^\natural)} = (-1)^{n/2} \cdot (-1)^{\ell(T)}$, and identifying $S^{\set{\lambda,+}}$ and $S^{\set{\lambda,-}}$ as in Remark \ref{remark:S(lam-lam)-decomposition}, one sees that
	\[
	\phi_{\ve''}^{(\lambda,\lambda)}( S^{\set{\lambda,\pm}} ) = \begin{cases}
	S^{\set{\lambda^*,\pm}} & \text{if $n/2$ is even,} \\
	S^{\set{\lambda^*,\mp}} & \text{if $n/2$ is odd.}
	\end{cases}
	\]
The restriction of $\phi_{\ve''}^{(\lambda,\lambda)}$ to $S^{\set{\lambda,\pm}}$ can thus be interpreted as a $\CD_n$-module isomorphism
	\[
	S^{\set{\lambda,\pm}} \otimes \ve'' \cong \begin{cases}
	S^{\set{\lambda^*,\pm}} & \text{if $n/2$ is even,} \\
	S^{\set{\lambda^*,\mp}} & \text{if $n/2$ is odd;}
	\end{cases}
	\]
cf.\ \cite[Lemma 3.5]{Geck:2015} (where the labeling of characters may be opposite to our labeling).
\end{remark}

\begin{lemma} \label{lemma:Dn-subminimal-dimensions}
Let $n \geq 5$.
\begin{enumerate}
\item The only one-dimensional $\CD_n$-modules are those labeled by $\set{[n],\emptyset}$ and $\set{[1^n],\emptyset}$, i.e., the trivial $\CD_n$-module and the one-dimensional module afforded by the sign character $\ve''$.

\item If $\set{\lambda,\mu} \in \calBP\{n\}$ and $\dim(S^{\set{\lambda,\mu}}) \neq 1$, then $\dim(S^{\set{\lambda,\mu}}) \geq n-1$.

\item If $\set{\lambda,\lambda} \in \calBP\{n\}$, then $\dim(S^{\set{\lambda,\pm }}) = \frac{1}{2} \cdot \binom{n}{n/2} \cdot \dim(S^\lambda)^2 \geq \frac{1}{2} \cdot \binom{n}{n/2}$.
\end{enumerate}
\end{lemma}

\begin{proof}
These are immediate consequences of Lemma \ref{lemma:Bn-subminimal-dimensions}.
\end{proof}

For $\set{\lambda,\mu} \in \calBP\{n \}$, let $\rho_{\set{\lambda,\mu}}: \CD_n \to \End(S^{\set{\lambda,\mu}})$ denote the $\CD_n$-module structure map for $S^{\set{\lambda,\mu}}$. Similarly, for $\lambda \vdash n/2$, let $\rho_{\set{\lambda,\pm}}: \CD_n \to \End(S^{\set{\lambda,\pm}})$ be the $\CD_n$-module structure map for $S^{\set{\lambda,\pm}}$. Then the Artin--Wedderburn Theorem implies:

\begin{theorem} \label{theorem:AW-type-D}
The structure maps of the simple $\CD_n$-modules induce an algebra isomorphism
	\begin{equation} \label{eq:AW-type-D}
	\CD_n \cong \Big[ \bigoplus_{\substack{\set{\lambda,\mu} \in \calBP\{n\} \\ \lambda \neq \mu}} \End(S^{\set{\lambda,\mu}}) \Big] \oplus \Big[ \bigoplus_{\set{\lambda,\lambda} \in \calBP\{n\}} \End(S^{\set{\lambda,+}}) \oplus \End(S^{\set{\lambda,-}}) \Big]. 
	\end{equation}
\end{theorem}

\subsection{Weight space decompositions of Type D simple modules} \label{subsec:Type-D-weight-spaces}

Observe that the YJM elements $X_1,\ldots,X_n$ are elements of $\CD_n$. Recall also the elements $u_1,\ldots,u_{n-1} \in \calD_n$ defined by $u_i = t_1 t_{i+1}$. Define $GZ(\calD_n)$ to be the commutative subalgebra of $\CD_n$ generated by $u_1,\ldots,u_{n-1}$ and $X_1,\ldots,X_n$. We call $GZ(\calD_n)$ the GZ-subalgebra of $\CD_n$. By restriction, each simple $\CB_n$-module $V$ decomposes into a direct sum of one-dimensional $GZ(\calD_n)$-submodules, spanned for example by the weight vectors identified in Theorem \ref{thm:type-B-normal-form-action}.

If $V_\alpha$ is a GZ-subspace of $V$ of weight $\alpha = (\rho,\chi)$, then the action of $GZ(\calD_n)$ on $V_{\alpha}$ is determined by the restricted weight
	\[ 
	\wt{\alpha} \coloneq (\wt{\rho},\chi) = (\wt{\rho}_1,\ldots,\wt{\rho}_{n-1},\chi_1,\ldots,\chi_n),
	\] 
where $\wt{\rho}_i \coloneq \rho_1 \rho_{i+1}$ is the eigenvalue by which $u_i = t_1 t_{i+1}$ acts on $V_\alpha$. One immediately checks that if $\alpha,\beta \in \calW(n)$ with $\alpha = (\rho,\chi)$, then $\wt{\alpha} = \wt{\beta}$ if and only if $\beta = \alpha$ or $\beta = (-\rho,\chi)$. Combined with \eqref{eq:weights-swap}, this implies for $(\lambda,\mu) \in \calBP(n)$ that the sets $\calW(\lambda,\mu)$ and $\calW(\mu,\lambda)$ have the same restrictions to $GZ(\calD_n)$. Then we can unambiguously define $\calW\{\lambda,\mu\} \coloneq \set{\wt{\alpha} : \alpha \in \calW(\lambda,\mu)}$. One further deduces that the union $\calW\{n\} \coloneq \bigcup_{\set{\lambda,\mu} \in \calBP\{n\}} \calW\{\lambda,\mu\}$ is disjoint. One sees for $\lambda \neq \mu$ that the weights in the set $\calW(\lambda,\mu)$ restrict to distinct elements in $\calW\{\lambda,\mu\}$, while for $\lambda = \mu$, the weights in $\calW(\lambda,\lambda)$ occur in pairs $(\pm \rho,\chi)$ that have common restrictions to $GZ(\calD_n)$.

\begin{lemma} \label{lemma:type-D-simple-weight-spaces}
Each simple $\CD_n$-module decomposes into GZ-subspaces for the action of $GZ(\calD_n)$:
	\begin{enumerate}
	\item \label{item:lam-mu-weight-space-typeD} If $\set{\lambda,\mu} \in \calBP\{n\}$ and $\lambda \neq \mu$, then $S^{\set{\lambda,\mu}} = \bigoplus_{\wt{\alpha} \in \calW\{\lambda,\mu\}} S^{\set{\lambda,\mu}}_{\wt{\alpha}}$.
	
	\item \label{item:lam-lam-weight-space-typeD} If $\set{\lambda,\lambda} \in \calBP\{n\}$, then $S^{\set{\lambda,\pm }} = \bigoplus_{\wt{\alpha} \in \calW\{\lambda,\lambda\}} S^{\set{\lambda,\pm}}_{\wt{\alpha}}$.
	\end{enumerate}
\end{lemma}

\begin{proof}
The weight space decomposition in \eqref{item:lam-mu-weight-space-typeD} is true by definition, while the decomposition in \eqref{item:lam-lam-weight-space-typeD} can be seen from the bases for $S^{\set{\lambda,+}}$ and $S^{\set{\lambda,-}}$ specified in \eqref{eq:(lam,pm)-bases}.
\end{proof}

\subsection{Orthogonal forms}

The actions of the Coxeter generators $s_i \in \fS_n$ on weight bases for the simple $\CD_n$-modules can be deduced by restriction from Theorem \ref{thm:type-B-normal-form-action} and the identifications in Remark \ref{remark:S(lam-lam)-decomposition}. In the next lemma we record a special case of this action.

\begin{lemma} \label{lemma:S(lam-pm)-si-action}
Suppose $n$ is even, let $\lambda \vdash n/2$, and identify $S^{\set{\lambda,\pm}}$ as in Remark \ref{remark:S(lam-lam)-decomposition}. Given $T \in \T(\lambda,\lambda)_{+1}$ and $1 \leq i < n$, set $r_i = \res_T(i+1)-\res_T(i)$. Then for $1 \leq i < n$, the action of the Coxeter generator $s_i = (i,i+1) \in \fS_n$ on $S^{\set{\lambda,\pm}}$ is as follows:

	\begin{enumerate}
	\item \label{item:i-i+1-different-tableaux} Suppose $i$ and $i+1$ are not in the same tableau in $T$. Let $S = s_i \cdot T$. Then with respect to the basis $\set{c_T \pm c_{T^\natural},c_S \pm c_{S^\natural}}$, $s_i$ acts on $S_T^{\set{\lambda,\pm}} \oplus S_S^{\set{\lambda,\pm}}$ via the matrix $\sm{0 & 1 \\ 1 & 0}$.
	
	\item If $i$ and $i+1$ are in the same column of the same tableau in $T$, implying that $r_i = -1$, then $s_i$ acts on $S_T^{\set{\lambda,\pm}}$ as multiplication by $-1$.
	
	\item If $i$ and $i+1$ are in the same row of the same tableau in $T$, implying that $r_i = +1$, then $s_i$ acts on $S_T^{\set{\lambda,\pm}}$ as multiplication by $+1$.
	
	\item Suppose $i$ and $i+1$ are in the same tableau in $T$, but not in the same row or in the same column, implying that $\abs{r_i} \geq 2$. Let $S = s_i \cdot T$. Then $S_T^{\set{\lambda,\pm}} \oplus S_S^{\set{\lambda,\pm}}$ is closed under the action of $s_i$, and the matrix of $s_i$ with respect to the basis $\{c_T \pm c_{T^\natural},c_S \pm c_{S^\natural} \}$ is
		\[
		\begin{bmatrix}
		r_i^{-1} & \sqrt{1-r_i^{-2}} \\
		\sqrt{1-r_i^{-2}} & -r_i^{-1}
		\end{bmatrix}.
		\]
	\end{enumerate}
\end{lemma}

\begin{remark}
In the case $i=1$ of the lemma, if $1$ and $2$ are not in the same tableau in $T$ and if $S = s_1 \cdot T$, then $S \notin \T(\lambda,\lambda)_{+1}$ but rather $S^{\natural} \in \T(\lambda,\lambda)_{+1}$. Then $c_{S^{\natural}} - c_S = -(c_S - c_{S^{\natural}})$ is the preferred basis vector for $S_{S^{\natural}}^{\set{\lambda,-}}$ identified in \eqref{eq:(lam,pm)-bases}, and with respect to the basis $\set{c_T - c_{T^\natural},c_{S^{\natural}} - c_{S}}$, the element $s_1$ acts on $S_T^{\set{\lambda,-}} \oplus S_{S^\natural}^{\set{\lambda,-}}$ via the matrix $\sm{ \hfill 0 & - 1 \\ - 1 & \hfill 0}$.
\end{remark}

\subsection{Restriction of Type D simple modules} \label{subsec:branching-type-D}

The `folding subgroup' of $\calD_n$ is the subgroup $\calH_{n-1}$ generated by the set $\set{s_1,\ldots,s_{n-2},t_{n-1}t_n}$. Then $(\calH_{n-1},\set{s_1,\ldots,s_{n-2},t_{n-1}t_n})$ is a Coxeter system of type $B_{n-1}$. Making the identification $\calH_{n-1} = \calB_{n-1}$, one has:

\begin{lemma}[{\cite[Theorem and Corollary 2.1]{Tokuyama:1984}}] \label{lemma:restriction:D-to-folding}
Let $(\lambda,\mu) \in \calBP(n)$. Then
	\[
	\Res^{\calD_n}_{\calH_{n-1}}(S^{\set{\lambda,\mu}}) \cong \Big[ \bigoplus_{\nu \prec \lambda} S^{(\nu,\mu)} \Big] \oplus \Big[ \bigoplus_{\tau \prec \mu} S^{(\tau,\lambda)} \Big].
	\]
If $n$ is even and $\lambda \vdash n/2$, then $\Res^{\calD_n}_{\calH_{n-1}}(S^{\set{\lambda,\pm}}) \cong \bigoplus_{\nu \prec \lambda} S^{(\nu,\lambda)}$.
\end{lemma}

Identify $\calD_{n-1}$ with the subgroup of $\calD_n$ generated by $\set{s_1,\ldots,s_{n-2},\wt{s}_{n-1}}$, where
	\[
	\wt{s}_{n-1} = t_{n-1}s_{n-2}t_{n-1} = t_{n-2}t_{n-1} s_{n-2}.
	\]
One also has $\wt{s}_{n-1} = (t_{n-1}t_n)s_{n-2}(t_{n-1}t_n)$, so $\calD_{n-1}$ is a subgroup of $\calH_{n-1}$ in precisely the same manner that $\calD_n$ is a subgroup of $\calB_n$. For $\set{\lambda,\mu} \in \calBP\{n\}$, write $\set{\nu,\tau} \prec \set{\lambda,\mu}$ if either $(\nu,\tau) \prec (\lambda,\mu)$ or $(\nu,\tau) \prec (\mu,\lambda)$. Then Lemmas \ref{lemma:Dn-simples} and \ref{lemma:restriction:D-to-folding} imply:

\begin{lemma} \label{lemma:restriction-D-to-Dn-1}
Let $(\lambda,\mu) \in \calBP(n)$. Then
	\begin{equation} \label{eq:S(lam-mu)-restriction-typeD}
	\Res^{\calD_n}_{\calD_{n-1}}(S^{\set{\lambda,\mu}}) \cong \bigoplus_{\set{\nu,\tau} \prec \set{\lambda,\mu}} S^{\set{\nu,\tau}} = \Big[ \bigoplus_{\nu \prec \lambda} S^{\set{\nu,\mu}} \Big] \oplus \Big[ \bigoplus_{\tau \prec \mu} S^{\set{\lambda,\tau}} \Big].
	\end{equation}
If $n$ is even and $\lambda \vdash n/2$, then
	\begin{equation} \label{eq:S(lam-pm)-restriction}
	\Res^{\calD_n}_{\calD_{n-1}}(S^{\set{\lambda,\pm}}) \cong \bigoplus_{\nu \prec \lambda} S^{\set{\nu,\lambda}}.
	\end{equation}
\end{lemma}

\begin{remark}
Let $(\lambda,\mu) \in \calBP(n)$. If $\mu \prec \lambda$, then the summand $S^{\set{\mu,\mu}}$ in \eqref{eq:S(lam-mu)-restriction-typeD} is not simple, but decomposes into the sum of the simple modules $S^{\set{\mu,+}}$ and $S^{\set{\mu,-}}$, neither of which occur in any other summand of \eqref{eq:S(lam-mu)-restriction-typeD}. Similar remarks apply if $\lambda \prec \mu$. In this way, one sees that if $V$ is a simple $\CD_n$-module, then the module $\Res_{\calD_{n-1}}^{\calD_n}(V)$ is multiplicity-free.
\end{remark}

\begin{remark}
Suppose $n$ is even, let $\lambda \vdash n/2$, and identify $S^{\set{\lambda,\pm}}$ as in Remark \ref{remark:S(lam-lam)-decomposition}. Then in terms of weight spaces, the decomposition \eqref{eq:S(lam-pm)-restriction} takes the form
	\begin{equation} \label{eq:Dn-restriction-weights}
	S^{\set{\lambda,\pm}} = \bigoplus_{\nu \prec \lambda} \Big[ \bigoplus_{\substack{T \in \T(\lambda,\lambda)_{+1}\\ \res_T(n) = \res(\lambda/\nu)}} S^{\set{\lambda,\pm}}_T \Big],
	\end{equation}
with the summand indexed by $\nu$ isomorphic as a $\CD_{n-1}$-module to $S^{\set{\nu,\lambda}}$. To verify that \eqref{eq:Dn-restriction-weights} realizes \eqref{eq:S(lam-pm)-restriction}, one can first apply Lemma \ref{lemma:S(lam-pm)-si-action} to check for each $\nu \prec \lambda$ that the summand indexed by $\nu$ is a $\CD_{n-1}$-submodule of $S^{\set{\lambda,\pm}}$. Then one can check for each $\nu \prec \lambda$ that the summand indexed by $\nu$ has the same set of weights for the action of $GZ(\calD_{n-1})$ as $S^{\set{\nu,\lambda}}$.
\end{remark}

\section{Lie algebras generated by reflections}\label{S:LA-by-reflections}

We begin in Section \ref{subsec:generalities} by recording some straightforward generalizations of results originally stated in the context of symmetric groups by Marin; our primary sources are \cite[III.3]{Marin:2001} and \cite[\S2]{Marin:2007}. In Section \ref{subsec:Lie-alg-gen-reflections} we introduce the Lie algebras $\fb_n$ and $\fd_n$ that are the main objects of study in this paper, and then in Section \ref{subsec:factorizations} we begin to describe the images of $\fb_n$ and $\fd_n$ in the various factors on the right-hand sides of the Artin--Wedderburn isomorphisms \eqref{eq:AW-type-B} and \eqref{eq:AW-type-D}.

\subsection{Generalities} \label{subsec:generalities}

Let $G$ be a finite group and let $\emptyset \neq S \subseteq G$. Write $\Lie(\CG)$ for the group algebra $\CG$ considered as a Lie algebra via the commutator bracket $[x,y] = xy - yx$. Let $\g$ be the Lie subalgebra of $\Lie(\CG)$ generated by the elements of the set $S$, let $\g' = [\g,\g]$ be the derived subalgebra of $\g$, and let $Z(\g)$ be the center of $\g$. Throughout this section we make the following assumptions (though in practice \ref{item:S-union-CC} can be weakened; see Remark \ref{remark:smaller-generating-set}): 
	\begin{enumerate}[label={(G\arabic*)}]
	\item \label{item:S-generates-G} $G = \subgrp{S}$,
	\item \label{item:S-order-2} the set $S$ consists of elements of order $2$,
	\item \label{item:S-union-CC} the set $S$ is a union of conjugacy classes $C_1,\ldots,C_t$.
	\end{enumerate}
	
\begin{lemma} \label{lemma:Z(g)}
With notation and assumptions as above, $Z(\g) = \g \cap Z(\CG)$, and the projection map $p: \CG \to Z(\CG)$, defined by $p(z) = \frac{1}{\abs{G}} \sum_{g \in G} gzg^{-1}$, restricts to a map $p: \g \to Z(\g)$.
\end{lemma}

\begin{proof}
Let $z \in \g$. Since $\g$ is generated by the elements of $S$, one has $z \in Z(\g)$ if and only if $[z,s] = 0$ for all $s \in S$. And since $S$ generates $\CG$ as an associative algebra, it follows that $z \in Z(\CG)$ if and only if $[z,s] = 0$ for all $s \in S$. Thus for $z \in \g$, one has $z \in Z(\g)$ if and only if $z \in Z(\CG)$. Finally, since $S$ is a union of conjugacy classes, it follows that $\g$ is closed under the conjugation action of $G$, and hence $p$ maps $\g$ into $\g \cap Z(\CG) = Z(\g)$.
\end{proof}

\begin{lemma} \label{lemma:g-reductive}
With notation and assumptions as above, the Lie algebra $\g$ is reductive, its derived subalgebra $\g'$ is semisimple, and $\g = \g' \oplus Z(\g)$.
\end{lemma}

\begin{proof}
Let $V_1,V_2,\ldots,V_m$ be a complete set of pairwise non-isomorphic simple $\CG$-modules. Then $V := \bigoplus_{i=1}^m V$ is a faithful finite-dimensional $\CG$-module. Since $\g$ contains the set $S$ of associative algebra generators for $\CG$, it follows that each $V_i$ restricts to a simple $\g$-module, and hence $V$ restricts to a faithful, finite-dimensional, semisimple $\g$-module. Then $\g$ is reductive, $\g'$ is semi\-simple, and $\g = \g' \oplus Z(\g)$ by Proposition 5 of \cite[Chapter I, \S6, no.\ 4]{Bourbaki:1998}.
\end{proof}

\begin{lemma}
With notation and assumptions as above, one has $\ker(p) = \g'$.
\end{lemma}

\begin{proof}
Since $S$ is a union of conjugacy classes, it follows that first $\g$, and then $\g'$, are closed under conjugation by $G$. Then $p$ maps $\g'$ into $\g' \cap Z(\g)$, which is zero by Lemma \ref{lemma:g-reductive}.
\end{proof}

For $1 \leq i \leq t$, set $T_i = \sum_{c \in C_i} c$. Since $T_i$ is a class sum in $\CG$, then $T_i \in Z(\CG)$.

\begin{lemma} \label{lemma:Z(g)-class-sums}
With notation and assumptions as above, $Z(\g) = \Span\{ T_1,\ldots,T_t \}$.
\end{lemma}

\begin{proof}
For $1 \leq i \leq t$, one has $T_i \in Z(\g)$ by Lemma \ref{lemma:Z(g)}. Conversely, if $z \in Z(\g)$, then $z \in Z(\CG)$ by the lemma, and hence $z$ is a linear combination of class sums in $\CG$. Given a conjugacy class $C$ in $G$, let $\Psi_C: \CG \to \C$ be the (linear extension to $\CG$ of the) class function that evaluates to $1$ on the elements of $C$ and that evaluates to $0$ on the other elements of $G$. If $0 \neq z \in Z(\CG)$, then $\Psi_C(z) \neq 0$ for some conjugacy class $C$.

Given a sequence $E = (s_1,\ldots,s_r)$ of elements of the set $S$, set $\ell(E) = r$, set $s_E = s_1$ if $r = 1$, and set $s_E = [s_1,s_{E'}]$ if $r > 1$ and $E' = (s_2,\ldots,s_r)$. Then $\g$ is spanned by the elements $s_E$ as $E$ ranges over all finite sequences of elements in $S$, and $\g'$ is spanned by all $s_E$ with $\ell(E) > 1$. If $\ell(E) > 1$ and $s = s_1$, then $s_E = [s,s_{E'}] = s s_{E'} - s_{E'} s$ can be written as a linear combination of terms of the form $sg - gs$ with $g \in G$. Since the elements of $S$ are of order $2$, one has $gs = s(sg)s^{-1}$. Then $gs$ and $sg$ are conjugate in $G$, and hence $\Psi_C([s,g]) = \Psi_C(sg)-\Psi_C(gs) = 0$. It follows for each conjugacy class $C$ of $G$ that $\Psi_C$ vanishes on $\g'$, and $\Psi_C$ is nonzero on $\g$ only if $C = C_i$ for some $1 \leq i \leq t$. Thus if $0 \neq z \in Z(\g)$, then $z$ must be a linear combination of $T_1,\ldots,T_t$.
\end{proof}

Given a $\CG$-module $V$, we may write $\Res^{\CG}_{\g}(V)$ for the $\g$-module obtained by restricting the $\CG$-action to the Lie subalgebra $\g \subseteq \Lie(\CG)$. Then $\Res^{\CG}_{\g}(-)$ defines a functor from $\CG$-modules to $\g$-modules. This is not necessarily a tensor functor, since for $\CG$-modules $V$ and $V'$, the definitions of the $\g$-actions on $\Res^{\CG}_{\g}(V \otimes V')$ and $\Res^{\CG}_{\g}(V) \otimes \Res^{\CG}_{\g}(V')$ are different.

\begin{lemma} \label{lemma:simples-remain-simple}
With notation and assumptions as above, if $V$ is a simple $\CG$-module, then $\Res^{\CG}_{\g}(V)$ is simple as a $\g$-module, and as a module over the derived subalgebra $\g'$.
\end{lemma}

\begin{proof}
First, as pointed out previously, $V$ is simple as a $\g$-module because $\g$ contains a set of associative algebra generators for $\CG$. Next, by Lemma \ref{lemma:Z(g)-class-sums}, $Z(\g)$ is spanned by central elements of $\CG$, which by Schur's Lemma act as scalars on $V$. Since $\g = \g' \oplus Z(\g)$ by Lemma \ref{lemma:g-reductive}, this implies that $V$ must also be simple as a $\g'$-module.
\end{proof}

If $V$ is a $\CG$-module, then the dual space $V^* = \Hom(V,\C)$ is a $\CG$-module with action defined for $g \in G$, $\phi \in V^*$, and $v \in V$ by $(g.\phi)(v) = \phi(g^{-1}.v)$. On the other hand, $[\Res^{\CG}_{\g}(V)]^*$ is a $\g$-module via the contragredient action, defined for $x \in \g$, $\phi \in V^*$, and $v \in V$ by $(x.\phi)(v) = -\phi(x.v)$; we denote $V^*$ with this $\g$-module structure by $V^{*,\Lie}$. In general, the $\g$-modules $\Res^{\CG}_{\g}(V^*)$ and $V^{*,\Lie}$ need not be isomorphic.

\begin{lemma} \label{lemma:tensor-sign-dual}
Let $\kappa: G \to \set{\pm 1}$ be a linear character that evaluates to $-1$ on each element of $S$. Then with notation and assumptions as above, $\Res^{\CG}_{\g}(V^* \otimes \kappa) \cong V^{*,\Lie}$ as $\g$-modules. If $V \cong V^*$ as $\CG$-modules, then $\Res^{\CG}_{\g}(V \otimes \kappa) \cong V^{*,\Lie}$ as $\g$-modules.
\end{lemma}

\begin{proof}
Let $\phi \in V^*$. We make the evident vector space identification $V^* \otimes \kappa = V^*$. Since $s = s^{-1}$ for each $s \in S$, and $\kappa(s) = -1$, one sees for each $s \in S$ that the action of $s$ on $\phi$ is the same whether $\phi$ is regarded as an element of $\Res^{\CG}_{\g}(V^* \otimes \kappa)$ or $V^{*,\Lie}$. Since the elements of $S$ generate $\g$, it follows that $\Res^{\CG}_{\g}(V \otimes \kappa) \cong V^{*,\Lie}$ as $\g$-modules. The last statement of the lemma follows by the functoriality of $\Res^{\CG}_{\g}(-)$.
\end{proof}

\begin{remark} \label{remark:smaller-generating-set}
Given a subset $\emptyset \neq S \subseteq G$ of a finite group $G$, define the conjugation closure of $S$, denoted $\ConjClo(S)$, to be the unique smallest subset of $G$ such that
	\begin{itemize}
	\item $S \subseteq \ConjClo(S)$, and
	\item $ghg^{-1} \in \ConjClo(S)$ for all $g,h \in \ConjClo(S)$.
	\end{itemize}
Then $\ConjClo(S) = \bigcup_{i \geq 1} S_i$, where the sets $S_1 \subseteq S_2 \subseteq S_3 \subseteq \cdots$ are defined by $S_1 = S$, and for $i > 1$ by $S_i = S_{i-1} \cup \{ ghg^{-1} : g,h \in S_{i-1} \}$. Now if $s \in G$ is an element of order $2$ and if $x \in \Lie(\CG)$, then
	\[ \textstyle
	sxs = x - \frac{1}{2} \cdot [s,[s,x]] \in \Lie(\CG).
	\]
Using this, one sees that if $S$ consists of elements order $2$, then $S$ and $\ConjClo(S)$ generate the same Lie sub\-algebras of $\Lie(\CG)$, and the results of this section remain true if the assumption \ref{item:S-union-CC} is replaced by the following weaker assumption:
	\begin{enumerate}[label={(G\arabic*${}^\prime$)}]
	\setcounter{enumi}{2}
	\item the set $\ConjClo(S)$ is a union of conjugacy classes $C_1,\ldots,C_t$.
	\end{enumerate}
\end{remark}

\subsection{Lie algebras generated by reflections} \label{subsec:Lie-alg-gen-reflections}

We now define our main objects of study.

\begin{definition}
Let $n \geq 2$.
	\begin{enumerate}
	\item Let $\fs_n \subseteq \Lie(\CS_n)$ be the Lie subalgebra generated by the set of transpositions in $\fS_n$.
	\item Let $\fb_n \subseteq \Lie(\CB_n)$ be the Lie subalgebra generated by the set \eqref{eq:B-reflections} of reflections in $\calB_n$.
	\item Let $\fd_n \subseteq \Lie(\CD_n)$ be the Lie subalgebra generated by the set \eqref{eq:D-reflections} of reflections in $\calD_n$.
	\end{enumerate}
\end{definition}

\begin{remark}
It follows from Remark \ref{remark:smaller-generating-set} that $\fb_n$ is also generated by the set $\set{s_1,\ldots,s_{n-1},t_1}$, where $s_i = (i,i+1)$, and $\fd_n$ is generated by the set $\set{s_1,\ldots,s_{n-1},\wt{s}_n}$, where $\wt{s}_n = t_{n-1}t_n s_{n-1}$.
\end{remark}

The structure of $\fs_n$ was determined by Marin \cite{Marin:2007}. Let $\calX_n = \sum_{j=1}^n X_j \in \fd_n \subseteq \fb_n$ be the sum of the YJM elements, and let $\calT_n = \sum_{j=1}^n t_j \in \fb_n$. Then Lemmas \ref{lemma:g-reductive} and \ref{lemma:Z(g)-class-sums} give:

\begin{lemma} \label{lemma:bn-dn-center} \ 
	\begin{enumerate}
	\item For $n \geq 2$, one has $\fb_n = \fb_n' \oplus Z(\fb_n)$, with $\fb_n'$ semisimple and $Z(\fb_n) = \Span\{ \calX_n,\calT_n \}$.

	\item For $n \geq 3$, one has $\fd_n = \fd_n' \oplus Z(\fd_n)$, with $\fd_n'$ semisimple and $Z(\fd_n) = \Span\{ \calX_n \}$.
	\end{enumerate}
\end{lemma}

\begin{remark}
For $n=2$, the group $\calD_2$ is abelian. Then $\fd_2 = Z(\fd_2)$ and $\fd_2' = 0$.
\end{remark}

The next proposition extends \cite[Proposition 2]{Marin:2007} to types B and D.

\begin{proposition} \label{prop:iso-equivalent}
Let $n \geq 2$, let $(G,\g)$ be either the pair $(\calB_n,\fb_n)$ or $(\calD_n,\fd_n)$, and let $V_1$ and $V_2$ be two simple $\CG$-modules of dimension greater than $1$.\footnote{All simple $\CD_2$-modules are one-dimensional, so for $(G,\g) = (\calD_n,\fd_n)$ this implies that $n \geq 3$.} Consider the following statements:
	\begin{enumerate}
	\item \label{item:G-iso} $V_1$ and $V_2$ are isomorphic as $\CG$-modules.
	\item \label{item:g-iso} $V_1$ and $V_2$ are isomorphic as $\g$-modules.
	\item \label{item:D(g)-iso} $V_1$ and $V_2$ are isomorphic as $\g'$-modules.
	\end{enumerate}
Statements \eqref{item:G-iso} and \eqref{item:g-iso} are equivalent, and \eqref{item:g-iso} implies \eqref{item:D(g)-iso}.
	\begin{itemize}
	\item If $(G,\g) = (\calD_n,\fd_n)$, then \eqref{item:D(g)-iso} implies \eqref{item:G-iso}, and hence all three statements are equivalent.
	\item If $(G,\g) = (\calB_n,\fb_n)$, then \eqref{item:D(g)-iso} is equivalent to one of the following holding:
		\begin{itemize}
		\item $V_1 \cong V_2$ as $\CB_n$-modules, or
		\item up to reordering, $V_1 \cong S^{(\lambda,\emptyset)}$ and $V_2 \cong S^{(\emptyset,\lambda)}$ for some $\lambda \vdash n$.
		\end{itemize}
	\end{itemize}
\end{proposition}

\begin{proof}
Statements \eqref{item:G-iso} and \eqref{item:g-iso} are equivalent because $\g$ contains a set of associative algebra generators for $\CG$, and \eqref{item:g-iso} evidently implies \eqref{item:D(g)-iso}. Now suppose \eqref{item:D(g)-iso} holds. Let $\phi: V_1 \to V_2$ be a $\g'$-module isomorphism, and let $\rho_1: \CG \to \End(V_1)$ and $\rho_2: \CG \to \End(V_2)$ be the $\CG$-module structure maps for $V_1$ and $V_2$, respectively. Then for all $x \in \g'$, one has $\phi \circ \rho_1(x) = \rho_2(x) \circ \phi$. Let $C_{(21^{n-2},\emptyset)}$ be the first conjugacy class identified in \eqref{eq:B-reflections}. The elements of $C_{(21^{n-2},\emptyset)}$ generate the subgroup $\calD_n$ of $G$ (which may be equal to $G$). Then one can argue as in the proof of \cite[Proposition 2]{Marin:2007} (cf.\ also the proof of \cite[Proposition 4.6.3]{DK:2023}) to show that either:
	\begin{enumerate}[label={(Iso\arabic*)}]
	\item \label{item:transpositions-intertwinor} $\phi \circ \rho_1(s) = \rho_2(s) \circ \phi$ for all $s \in C_{(21^{n-2},\emptyset)}$, or
	\item \label{item:transpositions-act-by-scalars} there exist $a_1,a_2 \in \C$ such that $\rho_1(s) = a_1 \cdot \Id_{V_1}$ and $\rho_2(s) = a_2 \cdot \Id_{V_2}$ for all $s \in C_{(21^{n-2},\emptyset)}$.
	\end{enumerate}

First suppose \ref{item:transpositions-act-by-scalars} holds. Since the elements of $C_{(21^{n-2},\emptyset)}$ generate $\calD_n$, it follows by multiplicativity that $\rho_1(g) \in \C \cdot \Id_{V_1}$ for all $g \in \calD_n$. If $G = \calD_n$, this contradicts the assumption that $V_1$ is a simple module of dimension greater than $1$, so suppose that $G = \calB_n$. The group $\calB_n$ is generated by $\calD_n$ and the element $t_1$, and each simple $\CB_n$-module decomposes into eigenspaces for $t_1$. Since the elements of $\calD_n$ act on $V_1$ as scalar multiples of the identity, this again leads to a contradiction of the assumption that $V_1$ is simple of dimension greater than $1$. Thus \ref{item:transpositions-act-by-scalars} cannot hold.

Now since \ref{item:transpositions-intertwinor} holds, it follow by multiplicativity that $\phi \circ \rho_1(g) = \rho_2(g) \circ \phi$ for all $g \in \calD_n$. Then the proposition is true if $G = \calD_n$, so suppose that $G = \calB_n$. Let $C_{(1^{n-1},1)}$ be the second conjugacy class identified in \eqref{eq:B-reflections}. Then one can again argue as in the proof of \cite[Proposition 2]{Marin:2007} to show that either:
	\begin{enumerate}[resume,label={(Iso\arabic*)}]
	\item \label{item:signs-intertwinor} $\phi \circ \rho_1(t) = \rho_2(t) \circ \phi$ for all $t \in C_{(1^{n-1},1)}$, or
	\item \label{item:signs-act-by-scalars} there exist $b_1,b_2 \in \C$ such that $\rho_1(t) = b_1 \cdot \Id_{V_1}$ and $\rho_2(t) = b_2 \cdot \Id_{V_2}$ for all $t \in C_{(1^{n-1},1)}$.
	\end{enumerate}
If \ref{item:signs-intertwinor} holds, then one deduces by multiplicativity that $\phi \circ \rho_1(g) = \rho_2(g) \circ \phi$ for all $g \in \calB_n$, and hence $V_1 \cong V_2$ as $\CB_n$-modules. So suppose \ref{item:signs-act-by-scalars} holds. In this case the scalars $b_1,b_2$ must each be either $+1$ or $-1$, since each simple $\CB_n$-module decomposes into $\pm 1$-eigenspaces for the action of the elements $t \in C_{(1^{n-1},1)}$. Since each $t \in C_{(1^{n-1},1)}$ is acting on $V_1$ by the same scalar value, we get from the descriptions in Sections \ref{subsec:simples-B} and \ref{subsec:weight-space-type-B} that $V_1 \cong S^{(\lambda,\emptyset)}$ or $V_1 \cong S^{(\emptyset,\lambda)}$ for some partition $\lambda \vdash n$, and similarly $V_2 \cong S^{(\mu,\emptyset)}$ or $V_2 \cong S^{(\emptyset,\mu)}$ for some $\mu \vdash n$. Condition \ref{item:transpositions-intertwinor} implies that $V_1 \cong V_2$ as $\C\fS_n$-modules, so $\lambda = \mu$. Thus either $V_1 \cong V_2$, or $\{V_1,V_2 \} \cong \{ S^{(\lambda,\emptyset)}, S^{(\emptyset,\lambda)} \}$.

Finally, if $V_1 \cong S^{(\lambda,\emptyset)}$ and $V_2 \cong S^{(\emptyset,\lambda)}$ for some $\lambda \vdash n$, then one can show that any $\CS_n$-module isomorphism $S^{(\lambda,\emptyset)} \cong S^\lambda \cong S^{(\emptyset,\lambda)}$ lifts to a $\fb_n'$-module isomorphism $S^{(\lambda,\emptyset)} \cong S^{(\emptyset,\lambda)}$. This follows from the observation that if $s_E$ is an iterated Lie bracket as in the proof of Lemma \ref{lemma:Z(g)-class-sums} with $\ell(E) > 1$, and if the list $E$ contains one or more elements from the conjugacy class $C_{(1^{n-1},1)}$, then $s_E$ acts as zero on both $S^{(\lambda,\emptyset)}$ and $S^{(\emptyset,\lambda)}$.
\end{proof}

\subsection{Factorizations} \label{subsec:factorizations}

Given $(\lambda,\mu) \in \calBP(n)$, write $\gl(\lambda,\mu)$ for the space $\End(S^{(\lambda,\mu)})$ considered as a Lie algebra via the commutator, and write $\fsl(\lambda,\mu)$ for the Lie subalgebra of traceless operators. Similarly, we may write $\gl\{\lambda,\mu \} = \End(S^{\set{\lambda,\mu}})$, $\gl\{\lambda,\pm \}$, $\fsl\{\lambda,\mu\}$, etc. Write
	\begin{align*}
	\fb_{(\lambda,\mu)} &= \rho_{(\lambda,\mu)}(\fb_n'),
	& 
	\fd_{\set{\lambda,\mu}} &= \rho_{\set{\lambda,\mu}}(\fd_n'),
	&
	\fd_{\set{\lambda,\pm}} &= \rho_{\set{\lambda,\pm}}(\fd_n')
	\end{align*}
for the (semisimple) images of $\fb_n'$ and $\fd_n'$ under the corresponding module structure maps. Then the Artin--Wedderburn isomorphisms of Theorems \ref{theorem:AW-type-B} and \ref{theorem:AW-type-D} restrict to injective maps
	\begin{gather}
	\fb_n' \to \bigoplus_{(\lambda,\mu) \in \calBP(n)} \fb_{(\lambda,\mu)}, \label{eq:bn-AW-map} \\
	\fd_n' \to \Big[ \bigoplus_{\substack{\set{\lambda,\mu} \in \calBP\{n\} \\ \lambda \neq \mu}} \fd_{\{ \lambda,\mu \}} \Big] \oplus 
	\Big[ \bigoplus_{\set{\lambda,\lambda} \in \calBP\{n\}} \fd_{\{\lambda,+ \}} \oplus \fd_{\{\lambda,- \}} \Big]. \label{eq:dn-AW-map}
	\end{gather}
By restriction one has $\fd_{\set{\lambda,\mu}} \subseteq \fb_{(\lambda,\mu)}$.

\subsubsection{Dual pairs} \label{subsubsec:dual-pairs}

Let $(\lambda,\mu) \in \calBP(n)$. As remarked after Theorem \ref{theorem:AW-type-B}, all simple $\CB_n$-modules are self-dual. Then from \eqref{eq:B-tensor-with-character} and Lemma \ref{lemma:tensor-sign-dual}, one gets $[S^{(\lambda,\mu)}]^{*,\Lie} \cong S^{(\mu^*,\lambda^*)}$ as $\fb_n$-modules. Thus, with respect to an appropriate choice of basis, the map $\rho_{(\mu^*,\lambda^*)}: \fb_n \to \gl(\mu^*,\lambda^*)$ can be written in the form $X \mapsto -\rho_{(\lambda,\mu)}(X)^t$, where $A^t$ denotes the transpose of $A$. This map evidently factors through $\rho_{(\lambda,\mu)}$, and one gets the first commutative diagram in \eqref{eq:b-d-factorization-diagrams} below. Using Remark \ref{remark:D-tensor-with-sign}, one similarly sees that $[S^{\set{\lambda,\mu}}]^{*,\Lie} \cong S^{\set{\lambda^*,\mu^*}}$ as $\fd_n$-modules, and deduces that $\rho_{\set{\lambda^*,\mu^*}}$ factors through $\rho_{\set{\lambda,\mu}}$, yielding the second commutative diagram in \eqref{eq:b-d-factorization-diagrams}.
	\begin{equation} \label{eq:b-d-factorization-diagrams}
	\vcenter{\xymatrix{
	\fb_n' \ar@{->>}[r]^{\rho_{(\lambda,\mu)}} \ar@{->>}[dr]_{\rho_{(\mu^*,\lambda^*)}} & \fb_{(\lambda,\mu)} \ar@{->}[d]^{\cong} \\
	 & \fb_{(\mu^*,\lambda^*)}
	}} \qquad \qquad \qquad
	\vcenter{\xymatrix{
	\fd_n' \ar@{->>}[r]^{\rho_{\set{\lambda,\mu}}} \ar@{->>}[dr]_{\rho_{\set{\lambda^*,\mu^*}}} & \fd_{\set{\lambda,\mu}} \ar@{->}[d]^{\cong} \\
	 & \fd_{\set{\lambda^*,\mu^*}}
	}}
	\end{equation}
If $\lambda \vdash n/2$, then one gets $\fd_n$-module isomorphisms
	\[
	[S^{\set{\lambda,\pm}}]^{*,\Lie} \cong \begin{cases}
	S^{\set{\lambda^*,\pm}} & \text{if $n/2$ is even,} \\
	S^{\set{\lambda^*,\mp}} & \text{if $n/2$ is odd.}
	\end{cases}
	\]
As above, it follows that one gets commutative diagrams
	\begin{equation} \label{eq:d-factorization-diagrams-lam-pm}
	\vcenter{\xymatrix{
	\fd_n' \ar@{->>}[r]^{\rho_{\set{\lambda,\pm}}} \ar@{->>}[dr]_{\rho_{\set{\lambda^*,\pm}}} & \fd_{\set{\lambda,\pm}} \ar@{->}[d]^{\cong} \\
	 & \fd_{\set{\lambda^*,\pm}}
	}} \quad \text{if $n/2$ is even, and} \qquad 
	\vcenter{\xymatrix{
	\fd_n' \ar@{->>}[r]^{\rho_{\set{\lambda,\pm}}} \ar@{->>}[dr]_{\rho_{\set{\lambda^*,\mp}}} & \fd_{\set{\lambda,\pm}} \ar@{->}[d]^{\cong} \\
	 & \fd_{\set{\lambda^*,\mp}}
	}} \quad \text{if $n/2$ is odd.}
	\end{equation}

Let $\sim$ be the equivalence relation on $\calBP(n)$ generated by $(\lambda,\mu) \sim (\mu^*,\lambda^*)$. By abuse of notation, we also let $\sim$ denote the induced relation on $\calBP\{n\}$, generated by $\set{\lambda,\mu} \sim \set{\lambda^*,\mu^*}$. Then from the preceding discussion one gets that the Artin--Wedderburn maps \eqref{eq:bn-AW-map} and \eqref{eq:dn-AW-map} factor respectively through the maps
	\begin{gather}
	\fb_n' \to \bigoplus_{(\lambda,\mu) \in \calBP(n)/\sim} \fb_{(\lambda,\mu)}, \label{eq:bn-AW-map-mod-conjugates} \\
	\fd_n' \to \Big[ \bigoplus_{\substack{\set{\lambda,\mu} \in \calBP\{n\}/\sim \\ \lambda \neq \mu}} \fd_{\{\lambda,\mu \}} \Big] \oplus 
	\Big[ \bigoplus_{\set{\lambda,\lambda} \in \calBP\{n\}/\sim } \fd_{\{\lambda,+\}} \oplus \fd_{\{\lambda,-\}}^\star \Big], \label{eq:dn-AW-map-mod-conjugates}
	\end{gather}
where the notation $\star$ means that the term $\fd_{\{\lambda,-\}}^\star$ is omitted if $\lambda = \lambda^*$ and $n/2$ is odd (since by the second diagram in \eqref{eq:d-factorization-diagrams-lam-pm}, $\fd_{\set{\lambda,-}} \cong \fd_{\set{\lambda,+}}$ and hence $\fd_{\set{\lambda,-}}$ is redundant in that case). Here, the sums over elements of $\calBP(n)/\!\!\sim$ or $\calBP\{n\}/\!\!\sim$ mean that we take the sum over any complete set of representatives in $\calBP(n)$ or $\calBP\{n\}$, respectively, for the distinct equivalence classes under $\sim$.

\subsubsection{Self-dual Lie algebra representations}

The arguments in this section are based on those in \cite[III.3.2.3]{Marin:2001}. First let $(\lambda,\mu) \in \calBP(n)$ with $(\lambda,\mu) = (\mu^*,\lambda^*)$. Then $\mu = \lambda^*$, $n$ is even, and $\lambda$ and $\mu$ are partitions of $n/2$. Let $\subgrp{-,-}_{(\lambda,\mu)}$ be the bilinear form on $S^{(\lambda,\mu)}$ discussed in Remark \ref{remark:C-bilinear-form}, and let $(-|-)_{(\lambda,\mu)}$ be the nondegenerate bilinear form on $S^{(\lambda,\mu)}$ defined in terms of the linear isomorphism $\phi_{\ve} = \phi_{\ve}^{(\lambda,\mu)}: S^{(\lambda,\mu)} \to S^{(\lambda,\mu)}$ of Remark \ref{remark:associators-define-isos} by
	\begin{equation} \label{eq:b-(lam,mu)-bilinear-form}
	(u|v)_{(\lambda,\mu)} = \subgrp{u,\phi_{\ve}(v)}_{(\lambda,\mu)}.
	\end{equation}

\begin{lemma} \label{lemma:b-osp-inclusion}
Let $(\lambda,\mu) \in \calBP(n)$ with $(\lambda,\mu) = (\mu^*,\lambda^*)$, and let $(-|-) = (-|-)_{(\lambda,\mu)}$ be the bilinear form on $S^{(\lambda,\mu)}$ defined by \eqref{eq:b-(lam,mu)-bilinear-form}. Then
	\begin{equation} \label{eq:b-osp(lam,mu)}
	\fb_{(\lambda,\mu)} \subseteq \osp(\lambda,\mu) := \{ x \in \gl(\lambda,\mu) : (x.u|v) + (u|x.v) = 0 \text{ for all } u,v \in S^{(\lambda,\mu)} \}.
	\end{equation}
The bilinear form $(-|-)$ is symmetric if $n/2$ is even and is anti-symmetric if $n/2$ is odd. Thus,
	\[
	\fb_{(\lambda,\mu)} \subseteq \begin{cases} 
	\fso(S^{(\lambda,\mu)}) & \text{if $n/2$ is even,} \\
	\fsp(S^{(\lambda,\mu)}) & \text{if $n/2$ is odd.}
	\end{cases}
	\]
By restriction $\fd_{\set{\lambda,\mu}} \subseteq \fb_{(\lambda,\mu)}$, in which case we may write $\fd_{\set{\lambda,\mu}} \subseteq \osp\{ \lambda,\mu\}$.
\end{lemma}

\begin{proof}
Let $\subgrp{-,-} = \subgrp{-,-}_{(\lambda,\mu)}$. For $\sigma \in \calB_n$ and $u,v \in S^{(\lambda,\mu)}$, one has
	\[
	(\sigma^{-1}.u|v) = \subgrp{\sigma^{-1}.u,\phi_\ve(v)} = \subgrp{u,\sigma \cdot \phi_\ve(v)} = \ve(\sigma) \cdot \subgrp{u,\phi_\ve(\sigma.v)} = \ve(\sigma) \cdot (u|\sigma.v).
	\]
This implies first for $s$ in the set \eqref{eq:B-reflections} of generators for $\fb_n$, and then for all $x \in \fb_n$ by linearity, that $(x.u|v) + (u|x.v) = 0$. Next consider the basis \eqref{eq:B(lam,mu)} for $S^{(\lambda,\mu)}$. One has $\phi_{\ve}(c_T) = (-1)^{\ell(T^\natural)} \cdot c_{(T^\natural)^*}$, so $(c_S|c_T) \neq 0$ if and only if $S = (T^\natural)^*$. Let $\sigma \in \fS_n$ be the permutation that maps $T \in \T(\lambda,\mu)$ to $(T^\natural)^*$. Then $\sigma = \sigma^{-1}$ is a product of $n/2$ disjoint transpositions, independent of $T$. Since \eqref{eq:B(lam,mu)} is an orthonormal basis for $S^{(\lambda,\mu)}$ with respect to the form $\subgrp{-,-}$, then $\sigma.c_T - \subgrp{\sigma.c_T,c_{(T^\natural)^*}} \cdot c_{(T^\natural)^*}$ is a linear combination of vectors $c_S$ with $S \neq (T^\natural)^*$. Further, $\subgrp{\sigma.c_T,c_{(T^\natural)^*}} \neq 0$ as a consequence of Lemma \ref{lemma:admissible-linear-combination}. Then it follows that
	\[
	(c_{(T^\natural)^*}|c_T) = \frac{1}{\subgrp{\sigma.c_T,c_{(T^\natural)^*}}} \cdot (\sigma.c_T|c_T) = \frac{\ve(\sigma)}{\subgrp{\sigma.c_T,c_{(T^\natural)^*}}} \cdot (c_T|\sigma.c_T) = \ve(\sigma) \cdot (c_T|c_{(T^\natural)^*})
	\]
Since $\ve(\sigma) = (-1)^{n/2}$, this implies that the form $(-|-)$ is symmetric if $n/2$ is even, and is anti-symmetric if $n/2$ is odd. Finally, $(-|-)$ is non-degenerate because $\subgrp{-,-}$ is non-degenerate and $\phi_{\ve}$ is a linear isomorphism.
\end{proof}

Next let $\set{\lambda,\mu} \in \calBP\{n\}$ with $\set{\lambda,\mu} = \set{\lambda^*,\mu^*}$. If $(\lambda,\mu) = (\mu^*,\lambda^*)$, then $\fd_{\set{\lambda,\mu}} \subseteq \osp\{\lambda,\mu\}$ as defined in Lemma \ref{lemma:b-osp-inclusion}, so suppose that $(\lambda,\mu) = (\lambda^*,\mu^*)$. Let $(-|-)_{\set{\lambda,\mu}}$ be the nondegenerate bilinear form on $S^{\set{\lambda,\mu}} = \Res^{\calB_n}_{\calD_n}(S^{(\lambda,\mu)})$ defined via the map $\phi_{\ve''} = \phi_{\ve''}^{(\lambda,\mu)}$ of Lemma \ref{lemma:epsilon''-intertwinor} by
	\begin{equation} \label{eq:d-(lam,mu)-bilinear-form}
	(u|v)_{\set{\lambda,\mu}} = \subgrp{u,\phi_{\ve''}(v)}_{(\lambda,\mu)}.
	\end{equation}

\begin{lemma} \label{lemma:d-osp(lam,mu)-inclusion}
Let $(\lambda,\mu) \in \calBP(n)$ with $(\lambda,\mu) = (\lambda^*,\mu^*)$, and let $(-|-) = (-|-)_{\set{\lambda,\mu}}$ be the bilinear form on $S^{\set{\lambda,\mu}}$ defined by \eqref{eq:d-(lam,mu)-bilinear-form}. Then
	\begin{equation} \label{eq:d-osp(lam,mu)}
	\fd_{\set{\lambda,\mu}} \subseteq \osp\{ \lambda,\mu \} := \{ x \in \gl\{ \lambda,\mu \} : (x.u|v) + (u|x.v) = 0 \text{ for all } u,v \in S^{\set{\lambda,\mu}} \}.
	\end{equation}
Let $b(\lambda)$ and $b(\mu)$ be the lengths of the diagonals of the Young diagrams of $\lambda$ and $\mu$, respectively. Then $(-|-)$ is symmetric if $\frac{n-b(\lambda) - b(\mu)}{2}$ is even, and is anti-symmetric otherwise.\footnote{If $\lambda = \lambda^*$ and $\mu = \mu^*$, then $\abs{\lambda} - b(\lambda)$ and $\abs{\mu} - b(\mu)$ must each be even, hence their sum $n - b(\lambda,\mu)$ is even.} Thus
	\[
	\fd_{\set{\lambda,\mu}} \subseteq \begin{cases} 
	\fso(S^{\set{\lambda,\mu}}) & \text{if $\frac{n-b(\lambda) -b(\mu)}{2}$ is even,} \\
	\fsp(S^{\set{\lambda,\mu}}) & \text{if $\frac{n-b(\lambda) -b(\mu)}{2}$ is odd.}
	\end{cases}
	\]
\end{lemma}

\begin{proof}
The proof is parallel to that of Lemma \ref{lemma:b-osp-inclusion}. In determining whether $(-|-)$ is symmetric or anti-symmetric, let $\sigma \in \fS_n \subset \calD_n$ be the permutation that maps $T \in \T(\lambda,\mu)$ to $T^*$. Then $\sigma$ is a product of $\frac{\abs{\lambda}-b(\lambda)}{2} + \frac{\abs{\mu}-b(\mu)}{2}$ disjoint transpositions, independent of $T$, so $\ve''(\sigma) = 1$ if $\frac{n-b(\lambda)-b(\mu)}{2}$ is even, and $\ve''(\sigma) = -1$ if $\frac{n-b(\lambda)-b(\mu)}{2}$ is odd.
\end{proof}

Finally, suppose $\set{\lambda,\lambda} \in \calBP\{n\}$ with $\lambda = \lambda^*$ and $n/2$ even. It follows from Remarks \ref{remark:S(lam-lam)-decomposition} and \ref{remark:D-tensor-with-sign} that the bilinear form \eqref{eq:d-(lam,mu)-bilinear-form} on $S^{\set{\lambda,\lambda}}$ restricts to a nondegenerate bilinear form on $S^{\set{\lambda,\pm}}$. Then by restriction one gets:

\begin{lemma} \label{lemma:d-osp(lam,pm)-inclusion}
Let $\set{\lambda,\lambda} \in \calBP\{n\}$ with $\lambda = \lambda^*$ and $n/2$ even, and let $(-|-) = (-|-)_{\set{\lambda,\pm}}$ be the restriction to $S^{\set{\lambda,\pm}}$ of the bilinear form on $S^{\set{\lambda,\lambda}}$ defined by \eqref{eq:d-(lam,mu)-bilinear-form}. Then
	\begin{equation} \label{eq:d-osp(lam,pm)}
	\fd_{\set{\lambda,\pm}} \subseteq \osp\{ \lambda,\pm \} := \{ x \in \gl\{ \lambda,\pm \} : (x.u|v) + (u|x.v) = 0 \text{ for all } u,v \in S^{\set{\lambda,\pm}} \}.
	\end{equation}
Let $b(\lambda)$ be the length of the diagonal of the Young diagram of $\lambda$. Then $(-|-)$ is symmetric if $b(\lambda)$ is even, and is anti-symmetric if $b(\lambda)$ is odd. Thus
	\[
	\fd_{\set{\lambda,\pm}} \subseteq \begin{cases} 
	\fso(S^{\set{\lambda,\pm}}) & \text{if $b(\lambda)$ is even,} \\
	\fsp(S^{\set{\lambda,\pm}}) & \text{if $b(\lambda)$ is odd.}
	\end{cases}
	\]
\end{lemma}

\begin{remark}
If $\set{\lambda,\lambda} \in \calBP\{n\}$ and $\lambda = \lambda^*$, then \eqref{eq:b-(lam,mu)-bilinear-form} and \eqref{eq:d-(lam,mu)-bilinear-form} define two different bilinear forms on $S^{\set{\lambda,\lambda}} = \Res^{\calB_n}_{\calD_n}(S^{(\lambda,\lambda)})$. For $n/2$ even, one can check that with respect to either form, the subspaces $S^{\set{\lambda,+}}$ and $S^{\set{\lambda,-}}$ of $S^{\set{\lambda,\lambda}}$ are orthogonal, the restrictions of the two forms to $S^{\set{\lambda,+}}$ are the same, and the restrictions of the two forms to $S^{\set{\lambda,-}}$ differ by the scalar $-1$. Thus, for $n/2$ even, the ambiguity in definition does not affect the conclusions of Lemma \ref{lemma:d-osp(lam,pm)-inclusion}. For $n/2$ odd, one gets that with respect to either form, the spaces $S^{\set{\lambda,+}}$ and $S^{\set{\lambda,-}}$ are each totally isotropic, and the forms each define perfect pairings (differing by the scalar $-1$) between $S^{\set{\lambda,+}}$ and $S^{\set{\lambda,-}}$.
\end{remark}

\subsubsection{Arm and leg (A\&L) bipartitions} \label{subsubsec:arm-and-leg}

For $1 \leq d \leq n-1$, set
	\[
	\beta_d = \beta_{n,d} = ([n-d],[1^d]) \qquad \text{and} \qquad \gamma_d = \gamma_{n,d} = ([1^d],[n-d]),
	\]
We refer to bipartitions of these forms as arm and leg (A\&L) bipartitions. Similarly, we call $\set{\lambda,\mu}$ an arm and leg if $(\lambda,\mu)$ is an arm and leg. (The notation $\beta_d$ and $\gamma_d$ also makes sense for $d=0$ or $d=n$, but we do not consider these edge cases to be A\&L bipartitions.) Set $\beta = \beta_{n,1}$ and $\gamma = \gamma_{n,1}$. Then $S^\beta$ is the natural representation for $\calB_n$, i.e., via the defining embedding $\calB_n \subseteq GL_n(\C)$ of Section \ref{subsec:Type-B-definition}, $S^\beta$ is the restriction to $\calB_n$ of the natural module for $GL_n(\C)$. By abuse of notation, we may write $\fd_{\beta_d}$ for $\fd_{\{[n-d],[1^d] \}}$.

\begin{proposition}[{\cite[5.5.7]{Geck:2000}}] \label{prop:exterior-powers-natural}
For $0 \leq d \leq n$, one has $S^{\beta_d} \cong \Lambda^d(S^{\beta})$ as $\CB_n$-modules.
\end{proposition}

Let $G$, $S$, and $\g$ be as in Section \ref{subsec:generalities}, and let $V$ be a $\CG$-module. In general, the $\g$-modules $\Res^{\CG}_{\g}(\Lambda^d(V))$ and $\Lambda^d(\Res^{\CG}_{\g}(V))$ are not isomorphic. Let $\rho_{\Grp}^{\wedge^d}: \g \to \End(\Lambda^d(V))$ be the $\g$-module structure map for $\Res^{\CG}_{\g}(\Lambda^d(V))$, and let $\rho_{\Lie}^{\wedge^d}: \g \to \End(\Lambda^d(V))$ be the $\g$-module structure map for $\Lambda^d(\Res^{\CG}_{\g}(V))$. Then for $s \in S$ and $v_1,v_2,\ldots,v_d \in V$, one has
	\begin{gather*}
	\rho_{\Grp}^{\wedge^d}(s)(v_1 \wedge v_2 \wedge \cdots \wedge v_d) = (s.v_1) \wedge (s.v_2) \wedge \cdots \wedge (s.v_d), \\
	\textstyle \rho_{\Lie}^{\wedge^d}(s)(v_1 \wedge v_2 \wedge \cdots \wedge v_d) = \sum_{i=1}^d v_1 \wedge \cdots \wedge (s.v_i) \wedge \cdots \wedge v_d.
	\end{gather*}

\begin{lemma}
For $0 \leq d \leq n$, one has $\Res^{\CB_n}_{\fb_n'}(\Lambda^d(S^\beta)) \cong \Lambda^d(\Res^{\CB_n}_{\fb_n'}(S^\beta))$ as $\fb_n'$-modules, and hence the $\CB_n$-module isomorphism of Proposition \ref{prop:exterior-powers-natural} restricts to an isomorphism of $\fb_n'$-modules, with $\fb_n'$ acting on the exterior power $\Lambda^d(S^\beta)$ in the Lie-algebraic manner via derivations.
\end{lemma}

\begin{proof}
This is seen via the same line of reasoning as given in the proof of \cite[Lemme 12]{Marin:2007}, applied now to the natural representation $S^\beta$ for $\calB_n$. In this case one must also consider the action of the generator $t_1$, for which one checks that $\rho_{\Lie}^{\wedge^d}(t_1) = \rho_{\Grp}^{\wedge^d}(t_1) + (d-1) \cdot \Id$. This is the same relation as for the generators $s_1,\ldots,s_{n-1}$, so the reasoning in the proof of \cite[Lemme 12]{Marin:2007} applies. 
\end{proof}

Next consider the $\CB_n$-modules $\Lambda^d(S^\beta) \otimes \ve' \cong S^{\beta_d} \otimes \ve' \cong S^{\gamma_d}$ and $\Lambda^d(S^\beta \otimes \ve') \cong \Lambda^d(S^\gamma)$. 

\begin{lemma}
For $0 \leq d \leq n$, the canonical vector space identifications $\Lambda^d(S^\beta) \otimes \ve' \cong \Lambda^d(S^\beta) \cong \Lambda^d(S^\beta \otimes \ve')$ induce $\fb_n'$-module isomorphisms
	\[
	\Res^{\CB_n}_{\fb_n'}(\Lambda^d(S^\beta) \otimes \ve') \cong \Lambda^d(\Res^{\CB_n}_{\fb_n'}(S^\beta \otimes \ve')).
	\]
Then for $0 \leq d \leq n$, there exists a $\fb_n'$-module isomorphism $S^{\gamma_d} \cong \Lambda^d(S^\gamma)$, with $\fb_n'$ acting on $\Lambda^d(S^\gamma)$ in the Lie-algebraic manner via derivations.
\end{lemma}

\begin{proof}
Let $\rho_{\Grp}^{\wedge^d}: \fb_n \to \End(\Lambda^d(S^\beta))$ and $\rho_{\Lie}^{\wedge^d}: \fb_n \to \End(\Lambda^d(S^\beta))$ be the $\fb_n$-module structure maps for $\Res^{\CB_n}_{\fb_n}(\Lambda^d(S^\beta) \otimes \ve')$ and $\Lambda^d(\Res^{\CB_n}_{\fb_n}(S^\beta \otimes \ve'))$, respectively. The proof again mirrors that of \cite[Lemme 12]{Marin:2007}. In this case one sees that $\rho_{\Lie}^{\wedge^d}(t_1) = \rho_{\Grp}^{\wedge^d}(t_1) - (d-1) \cdot \Id$, while for $s_1,\ldots,s_{n-1}$ one again gets $\rho_{\Lie}^{\wedge^d}(s_i) = \rho_{\Grp}^{\wedge^d}(s_i) + (d-1) \cdot \Id$. Thus, for each element $s$ in a set of generators for $\fb_n$, one sees that $\rho_{\Lie}^{\wedge^d}(s)$ and $\rho_{\Grp}^{\wedge^d}(s)$ differ by a central element of $\End(\Lambda^d(S^\beta))$. This implies that the restrictions of $\rho_{\Lie}^{\wedge^d}$ and $\rho_{\Grp}^{\wedge^d}$ to $\fb_n'$ are the same.
\end{proof}

For $\eta \in \{ \beta,\gamma \}$, let $\Delta_{d}: \fsl(\eta) \to \fsl(\eta_{d}) \cong \fsl(\Lambda^d(S^\eta))$ denote the Lie algebra homomorphism corresponding to the action of $\fsl(S^{\eta})$ on $S^{\eta_d} \cong \Lambda^d(S^\eta)$ via derivations.

\begin{corollary} \label{cor:b-d-exterior-powers-factorization}
For $0 \leq d \leq n$, the module structure maps $\rho_{\beta_d}: \fb_n' \to \fsl(\beta_d)$, $\rho_{\gamma_d}: \fb_n' \to \gl(\gamma_d)$, and $\rho_{\beta_d}: \fd_n' \to \fsl(\beta_d)$ admit Lie algebra factorizations
	\begin{equation} \label{eq:linear-factorizations}
	\vcenter{\xymatrix{
	\fd_n' \ar@{^(->}[r] & \fb_n' \ar@{->}[r]^{\rho_{\beta}} \ar@{->}[dr]_{\rho_{\beta_d}} & \fsl(\beta) \ar@{->}[d]^{\Delta_d} \\ 
	&  & \fsl(\beta_d)
	}} \qquad \text{and} \qquad
	\vcenter{\xymatrix{
	\fb_n' \ar@{->}[r]^{\rho_{\gamma}} \ar@{->}[dr]_{\rho_{\gamma_d}} & \fsl(\gamma) \ar@{->}[d]^{\Delta_d} \\
	 & \fsl(\gamma_d)
	}}.
	\end{equation}
\end{corollary}

\subsubsection{Improper bipartitions} \label{subsubsec:improper-bipartitions}

We say that $(\lambda,\mu) \in \calBP(n)$ is \emph{proper} if $\lambda$ and $\mu$ are both nonempty, and that $(\lambda,\mu)$ is \emph{improper} otherwise. We may similarly refer to an element of $\calBP\{n\}$ as (im)proper. Let $\calBP(n)^\circ$ and $\calBP\{n\}^\circ$ be the subsets of proper elements in $\calBP(n)$ and $\calBP\{n\}$, respectively.

The split sequence of group maps $\fS_n \hookrightarrow \calB_n \twoheadrightarrow \fS_n$ gives rise to a split sequence of Lie algebras
	\begin{equation} \label{eq:split-sequence-group-algebras}
	\Lie(\CS_n) \hookrightarrow \Lie(\CB_n) \twoheadrightarrow \Lie(\CS_n).
	\end{equation}
In general, the quotient $\Lie(\CB_n) \twoheadrightarrow \Lie(\CS_n)$ does not restrict to a Lie algebra map from $\fb_n$ to $\fs_n$, because the generator $t_i \in \fb_n$ maps to the identity $1 \in \CS_n$, which in general is not an element of $\fs_n$. But considering iterated Lie brackets as in the proofs of Lemma \ref{lemma:Z(g)-class-sums} and Proposition \ref{prop:iso-equivalent}, one sees that any iterated Lie bracket $s_E$ with $\ell(E) > 1$ that is constructed using one or more elements from the conjugacy class $C_{(1^{n-1},1)}$ maps to zero under the quotient $\Lie(\CB_n) \twoheadrightarrow \Lie(\CS_n)$. This implies that \eqref{eq:split-sequence-group-algebras} restricts to a split sequence of derived algebras $\fs_n' \hookrightarrow \fb_n' \twoheadrightarrow \fs_n'$.

Now let $\lambda \vdash n$. As stated in Section \ref{subsec:simples-B}, one has $S^{(\lambda,\emptyset)} = \Inf_{\fS_n}^{\calB_n}(S^\lambda)$ and $S^{(\emptyset,\lambda)} = \Inf_{\fS_n}^{\calB_n}(S^\lambda) \otimes \ve'$, with $S^{(\lambda,\emptyset)} = S^\lambda = S^{(\emptyset,\lambda)}$ as modules for $\CS_n \subset \CB_n$. Set $\gl(\lambda) = \End(S^\lambda)$. Since the action of $\calB_n$ on $S^{(\lambda,\emptyset)}$ factors through the quotient map $\calB_n \twoheadrightarrow \fS_n$, it follows that the $\fb_n'$-module structure map $\rho_{(\lambda,\emptyset)}: \fb_n' \to \gl(\lambda)$ factors as $\fb_n' \twoheadrightarrow \fs_n' \xrightarrow{\rho_\lambda} \gl(\lambda)$, where $\rho_\lambda: \fs_n' \to \gl(\lambda)$ is the $\fs_n'$-module structure map for $S^\lambda$. By restriction, the $\fd_n'$-module structure map $\rho_{\set{\lambda,\emptyset}}$ also factors through $\rho_\lambda$. Set $\fs_\lambda = \rho_\lambda(\fs_n')$. Then for $\lambda \vdash n$, the Lie algebra maps $\fs_n' \hookrightarrow \fd_n' \hookrightarrow \fb_n' \twoheadrightarrow \fs_n'$ induce equalities
	\begin{equation} \label{eq:d-b-s-identifications}
	\fd_{\set{\lambda,\emptyset}} = \fb_{(\lambda,\emptyset)} = \fs_\lambda.
	\end{equation}

\begin{remark}
Let $\lambda \vdash n$. The reasoning in the last paragraph of the proof of  Proposition \ref{prop:iso-equivalent} implies that $\rho_{(\emptyset,\lambda)}: \fb_n' \to \gl(\lambda)$ also factors through the quotient $\fb_n' \twoheadrightarrow \fs_n'$, and hence $\fb_{(\emptyset,\lambda)}$ also identifies with $\fs_\lambda$. Thus for $\lambda \vdash n$, one has $\fs_\lambda = \fb_{(\emptyset,\lambda)} = \fb_{(\lambda,\emptyset)} \cong \fb_{(\emptyset,\lambda^*)} = \fb_{(\lambda^*,\emptyset)}$.
\end{remark}
	
The structure of $\fs_\lambda$ was determined by Marin in \cite{Marin:2007}. We summarize his results. Say that a partition $\lambda$ is a \emph{hook} if $\lambda = [n-d,1^d]$ for some $0 \leq d \leq n-1$.\footnote{Marin refers to partitions of this form as \emph{\'{e}querres}.} Let
	\begin{align*}
	E_n &= \set{ \lambda \vdash n : \text{$\lambda$ is not a hook and $\lambda \neq \lambda^*$}}, \\
	F_n &= \set{ \lambda \vdash n : \text{$\lambda$ is not a hook and $\lambda = \lambda^*$}}.
	\end{align*}

\begin{itemize}
\item If $\dim(S^\lambda) = 1$, then $\fs_\lambda = \fsl(\lambda) = 0$.

\item For $\lambda \vdash n$, one has $S^{\lambda^*} \cong [S^\lambda]^{*,\Lie}$ as $\fs_n$-modules.

\item If $\lambda \in E_n$, then $\fs_\lambda = \fsl(\lambda)$, the map $\rho_{\lambda^*}: \fs_n' \to \fs_{\lambda^*}$ factors through $\rho_\lambda$, and $\fs_\lambda \cong \fs_{\lambda^*}$.

\item Given $\lambda \in F_n$, let $b(\lambda)$ denote the length of the diagonal in the Young diagram of $\lambda$. Then there exists a nondegenerate bilinear form\footnote{Up to a nonzero scalar multiple, this is the form $(-|-)_{\set{\lambda,\emptyset}}$ on $S^{\set{\lambda,\emptyset}} = S^\lambda$ defined prior to Lemma \ref{lemma:d-osp(lam,mu)-inclusion}.} $(\;|\;)$ on $S^\lambda$, which is symmetric (resp.\ anti-symmetric) if $(n-b(\lambda))/2$ is even (resp.\ odd), and with respect to which
	\[
	\fs_\lambda = \osp(\lambda) := \{ x \in \gl(\lambda): (x.u|v)+ (u|x.v) = 0 \text{ for all } u,v \in S^\lambda \}.
	\]
Thus $\fs_\lambda$ is either a symplectic or an orthogonal Lie algebra.

\item For $0 \leq d \leq n-1$, set $\alpha_d = \alpha_{n,d} = [n-d,1^d]$, and set $\alpha = \alpha_{n,1}$. Then $S^{\alpha_d} \cong \Lambda^d(S^{\alpha})$ as $\fs_n'$-modules by \cite[Lemme 12]{Marin:2007}. The map $\rho_{\alpha_d}: \fs_n' \to \fsl(\alpha_d)$ factors through $\rho_{\alpha}: \fs_n' \to \fsl(\alpha)$, and $\fs_{\alpha} = \fsl(\alpha) \cong \fsl(n-1)$. Then $\fs_{\alpha_d} \cong \fs_\alpha$ if $1 \leq d \leq n-2$, and $\fs_{\alpha_d} = 0$ otherwise.
\end{itemize}

\begin{theorem}[{\cite[Theorem A]{Marin:2007}}] \label{thm:Marin-Theorem-A}
Let $\sim$ be the relation on $\set{ \lambda : \lambda \vdash n}$ generated by $\lambda \sim \lambda^*$. Then for $n \geq 2$, the Artin--Wedderburn map $\fs_n' \to \bigoplus_{\lambda \vdash n} \fs_\lambda$ factors through the map
	\begin{equation} \label{eq:Marin-Theorem-A}
	\fs_n' \to \calL_n := \fs_\alpha \oplus \Big[ \bigoplus_{\lambda \in E_n/\sim} \fs_\lambda \Big] \oplus \Big[ \bigoplus_{\lambda \in F_n} \fs_\lambda \Big] 
	= \fsl(\alpha) \oplus \Big[ \bigoplus_{\lambda \in E_n/\sim} \fsl(\lambda) \Big] \oplus \Big[ \bigoplus_{\lambda \in F_n} \osp(\lambda) \Big],
	\end{equation}
and this map is a Lie algebra isomorphism.\footnote{Marin states the theorem for $n \geq 3$. It is also vacuously true for $n=2$, since then all terms are zero.}
\end{theorem}

\subsubsection{Summary} \label{SS:AW-summary}

Recall that $\calBP(n)^\circ$ and $\calBP\{n\}^\circ$ denote the subsets of proper elements in $\calBP(n)$ and $\calBP\{n\}$, respectively. Define subsets of $\calBP(n)^\circ$ and $\calBP\{n\}^\circ$ by
	\begin{align*}
	E(n) &= \set{ (\lambda,\mu) \in \calBP(n)^\circ : \text{$(\lambda,\mu)$ is not an A\&L and $(\lambda,\mu) \neq (\mu^*,\lambda^*)$} }, \\
	F(n) &= \set{ (\lambda,\mu) \in \calBP(n)^\circ : \text{$(\lambda,\mu)$ is not an A\&L and $(\lambda,\mu) = (\mu^*,\lambda^*)$} }, \\
	E\{n\} &= \set{ \set{\lambda,\mu} \in \calBP\{n\}^\circ : \text{$\set{\lambda,\mu}$ is not an A\&L and $\set{\lambda,\mu} \neq \set{\lambda^*,\mu^*}$}}, \\
	F\{n\} &= \set{ \set{\lambda,\mu} \in \calBP\{n\}^\circ : \text{$\set{\lambda,\mu}$ is not an A\&L and $\set{\lambda,\mu} = \set{\lambda^*,\mu^*}$}}.
	\end{align*}
In particular, $F(n) \neq \emptyset$ only if $n$ is even. Now making the identifications in \eqref{eq:d-b-s-identifications}, the Lie algebra $\calL_n$ in \eqref{eq:Marin-Theorem-A} identifies with summands of the right-hand sides of \eqref{eq:bn-AW-map-mod-conjugates} and \eqref{eq:dn-AW-map-mod-conjugates} respectively. Combined with the factorizations provided by Corollary \ref{cor:b-d-exterior-powers-factorization}, it follows that \eqref{eq:bn-AW-map-mod-conjugates} and \eqref{eq:dn-AW-map-mod-conjugates} can be refined as follows:

\begin{proposition}\label{prop:restriction-of-AW-map}
For $n \geq 2$, the (injective) Artin--Wedderburn maps \eqref{eq:bn-AW-map} and \eqref{eq:dn-AW-map} factor respectively through the Lie algebra homomorphisms
	\begin{equation} \label{eq:bn-refined-factorization}
	\fb_n' \to \calL_n \oplus \fb_{\beta} \oplus \fb_{\gamma} \oplus \Big[ \bigoplus_{(\lambda,\mu) \in E(n)/\sim} \fb_{(\lambda,\mu)} \Big] \oplus \Big[ \bigoplus_{(\lambda,\mu) \in F(n)} \fb_{(\lambda,\mu)} \Big],
	\end{equation}
and
	\begin{equation} \label{eq:dn-refined-factorization}
	\begin{split}
	\fd_n' \to \calL_n \oplus \fd_{\beta} 
	&\oplus \Big[ \bigoplus_{\substack{\set{\lambda,\mu} \in E\{n\}/\sim \\ \lambda \neq \mu}} \fd_{\set{\lambda,\mu}} \Big] 
	\oplus \Big[ \bigoplus_{\set{\lambda,\lambda} \in E\{n\}/\sim} \fd_{\set{\lambda,+}} \oplus \fd_{\set{\lambda,-}} \Big] \\
	&\oplus \Big[ \bigoplus_{\substack{\set{\lambda,\mu} \in F\{n\} \\ \lambda \neq \mu}} \fd_{\set{\lambda,\mu}} \Big] 
	\oplus \Big[ \bigoplus_{\set{\lambda,\lambda} \in F\{n\}} \fd_{\set{\lambda,+}} \oplus \fd_{\set{\lambda,-}}^\star \Big],
	\end{split}
	\end{equation}
where the notation $\star$ means that the term $\fd_{\set{\lambda,-}}^\star$ is omitted if $n/2$ is odd.
\end{proposition}

\section{The main theorem}\label{S:main-theorem}

\subsection{Statement of the main theorem}

The main theorem of the paper is:

\begin{theorem} \label{thm:main-theorem}
Let $n \geq 2$. Then the maps \eqref{eq:bn-refined-factorization} and \eqref{eq:dn-refined-factorization} are Lie algebra isomorphisms (in particular, they are surjections), and the following identifications hold:
	\begin{enumerate}
	\item \label{item:Ln-decomp} $\calL_n = \fsl(\alpha) \oplus [ \bigoplus_{\lambda \in E_n/\sim} \fsl(\lambda) ] \oplus [ \bigoplus_{\lambda \in F_n} \osp(\lambda) ]$ as described in Section \ref{subsubsec:improper-bipartitions}.

	\item \label{item:beta-gamma} $\fd_\beta = \fb_\beta = \fsl(\beta)$ and $\fb_\gamma = \fsl(\gamma)$. For $0 < d < n$, one has $\fd_{\beta_d} = \fb_{\beta_d} \cong \fsl(\beta)$ and $\fb_{\gamma_d} \cong \fsl(\gamma)$.

	\item \label{item:F(n)-osp(lam,mu)} If $(\lambda,\mu) \in F(n)$, then $\fb_{(\lambda,\mu)} = \osp(\lambda,\mu)$ as described in Lemma \ref{lemma:b-osp-inclusion}.

	\item \label{item:E(n)-sl(lam,mu)} If $(\lambda,\mu) \in E(n)$, then $\fb_{(\lambda,\mu)} = \fsl(\lambda,\mu)$.

	\item \label{item:Fn-osp(lam,mu)} If $\set{\lambda,\mu} \in F\{n\}$ and $\lambda \neq \mu$, then $\fd_{\set{\lambda,\mu}} = \osp\{\lambda,\mu\}$ as described in Lemma \ref{lemma:b-osp-inclusion} or \ref{lemma:d-osp(lam,mu)-inclusion}.

	\item \label{item:Fn-osp(lam,pm)-d} If $\set{\lambda,\lambda} \in F\{n\}$ and $n/2$ is even, then $\fd_{\set{\lambda,\pm}} = \osp\{\lambda,\pm \}$ as described in Lemma \ref{lemma:d-osp(lam,pm)-inclusion}.

	\item \label{item:Fn-sl(lam,+)} If $\set{\lambda,\lambda} \in F\{n\}$ and $n/2$ is odd, then $\fd_{\set{\lambda,\pm}} = \fsl\{\lambda,\pm \}$.

	\item \label{item:En-sl(lam,mu)-d} If $\set{\lambda,\mu} \in E\{n\}$ and $\lambda \neq \mu$, then $\fd_{\set{\lambda,\mu}} = \fsl\{\lambda,\mu\}$.

	\item \label{item:En-sl(lam,pm)-d} If $\set{\lambda,\lambda} \in E\{n\}$, then $\fd_{\set{\lambda,\pm}} = \fsl\{\lambda,\pm \}$.
	\end{enumerate}
\end{theorem}

The proof of Theorem \ref{thm:main-theorem} is by induction on $n$. The base case of induction is handled by:

\begin{lemma} \label{lemma:main-theorem-base-case}
Theorem \ref{thm:main-theorem} is true for $n \in \set{2,3,4,5}$.
\end{lemma}

\begin{proof}
By Lemma \ref{lemma:bn-dn-center} and Proposition \ref{prop:restriction-of-AW-map}, the dimensions of $\fb_n' \subseteq \fb_n$ and $\fd_n' \subseteq \fd_n$ for $n \in \set{2,3,4,5}$ are less than or equal to the values specified in Table \ref{tab:upper-bounds}. One can then verify using GAP \cite{GAP4} that the dimensions are equal to these upper bounds, and hence the injective maps in Proposition \ref{prop:restriction-of-AW-map} must be isomorphisms by dimension comparison. A copy of the GAP code we used to carry out these calculations is included as an ancillary file with the version of this paper posted on the arXiv (\href{https://arxiv.org/abs/2506.01198}{arXiv:2506.01198}).
\end{proof}

\begin{center}
\begin{table}[tb]
\begin{tabular}{|l|l|l||l|l|}
\hline
$n$ & $\dim(\fb_n') \leq$ & $\dim(\fb_n) \leq$ & $\dim(\fd_n') \leq$ & $\dim(\fd_n) \leq $ \\
\hline
2 & 3 & 5 & 0 & 1 \\
\hline
3 & 19 & 21 & 11 & 12 \\
\hline
4 & 139 & 141 & 78 & 79 \\
\hline
5 & 1630 & 1632 & 840 & 841 \\
\hline
\end{tabular}

\ 

\caption{Upper bounds on dimensions of $\fb_n'$ and $\fb_n'$.} \label{tab:upper-bounds}
\end{table}
\end{center}

For the induction step we let $n \geq 6$ be fixed and we assume that Theorem \ref{thm:main-theorem} is true for the value $n-1$. Under this assumption we first show in Sections \ref{subsec:item:beta-gamma}--\ref{subsec:En-sl(lam,pm)-d} that parts \eqref{item:beta-gamma}--\eqref{item:En-sl(lam,pm)-d} of Theorem \ref{thm:main-theorem} are true for the value $n$; part \eqref{item:Ln-decomp} is true for all $n \geq 2$ by the identifications in \eqref{eq:d-b-s-identifications} and by Theorem \ref{thm:Marin-Theorem-A}. Then in Section \ref{subsec:complete-main-theorem} we show that the maps \eqref{eq:bn-refined-factorization} and \eqref{eq:dn-refined-factorization} are surjections. Since the maps are injections by construction, this finishes the proof.

In broad strokes, the proofs of parts \eqref{item:beta-gamma}--\eqref{item:En-sl(lam,pm)-d} are all the same: In addition to the given unknown Lie algebra $\g$ (either $\fd_\beta$, $\fb_\beta$, $\fb_\gamma$, $\fb_{(\lambda,\mu)}$, etc.), we consider a semisimple subalgebra $\fh \subseteq \g$ arising from the corresponding Weyl group of rank $n-1$. From the restriction formulas in Sections \ref{subsec:branching-type-B} and \ref{subsec:branching-type-D}, we first deduce that $\g$ is simple using \cite[Lemme 15]{Marin:2007}. Then in most cases, we see with the help of \cite[Lemme 13]{Marin:2007} and \cite[Lemme 14]{Marin:2007} that the rank of $\fh$ is sufficiently large to force $\g$ to have the claimed value. There is however one case where this strategy is not entirely sufficient to determine $\g$: For part \eqref{item:Fn-sl(lam,+)} we must apply more detailed information about the module $S^{\set{\lambda,\pm}}$ to rule out the possibility that $\g$ is either an orthogonal or symplectic Lie algebra.

\subsection{Proof of Theorem \ref{thm:main-theorem}(\ref{item:beta-gamma})} \label{subsec:item:beta-gamma}

\begin{lemma} \label{lemma:proof-of-item-beta-gamma}
Let $n \geq 6$, and suppose Theorem \ref{thm:main-theorem} is true for the value $n-1$. Then part \eqref{item:beta-gamma} of Theorem \ref{thm:main-theorem} is true for the value $n$, i.e.,
	\begin{center}
	$\fd_\beta = \fb_\beta = \fsl(\beta)$ and $\fb_\gamma = \fsl(\gamma)$. For $0 < d < n$, one has $\fd_{\beta_d} = \fb_{\beta_d} \cong \fsl(\beta)$ and $\fb_{\gamma_d} \cong \fsl(\gamma)$.
	\end{center}
\end{lemma}

\begin{proof}
First we show that the inclusions $\fd_\beta \subseteq \fb_\beta \subseteq \fsl(\beta)$ are equalities; the proof that the inclusion $\fb_\gamma \subseteq \fsl(\gamma)$ is an equality is entirely parallel. By Lemma \ref{lemma:restriction-D-to-Dn-1},
	\[
	\Res^{\calD_n}_{\calD_{n-1}}(S^\beta) \cong S^{\set{[n-2],[1]}} \oplus S^{\set{[n-1],\emptyset}}.
	\]
The dimensions of the modules are $\dim(S^\beta) = n$, $\dim(S^{\set{[n-2],[1]}}) = n-1$, and $\dim(S^{\set{[n-1],\emptyset}}) = 1$. Let $\fh = \rho_\beta(\fd_{n-1}')$. Then $\fh \subseteq \fd_{\beta}$ are two semisimple subalgebras of $\fsl(\beta)$, and from the induction hypothesis one gets $\fh = \fsl \{[n-2],[1] \}$. In particular, $\fh$ is semisimple of rank $\rk(\fh) = n-2$, while $\rk(\fd_\beta) \leq \rk(\fsl(\beta)) = n-1$, so $\rk(\fd_\beta) < 2 \cdot \rk(\fh)$. Then $\fd_\beta$ is simple by \cite[Lemme 15]{Marin:2007}. Now $2 \cdot \rk(\fd_\beta) \geq 2 \cdot \rk(\fh) = 2(n-2) > n = \dim(\beta)$, so $\fd_\beta = \fsl(\beta)$ by \cite[Lemme 13]{Marin:2007}.

Finally, commutativity of the diagrams in \eqref{eq:linear-factorizations} implies that $\fd_{\beta_d}$ is a Lie algebra quotient of the simple Lie algebra $\fd_\beta = \fsl(\beta)$. For $0 < d < n$, this quotient must be nonzero, because $S^{\beta_d}$ is a simple $\fd_n'$-module (by Lemma \ref{lemma:simples-remain-simple}) of dimension greater than $1$. Then $\fd_{\beta_d} \cong \fsl(\beta)$. By parallel reasoning, one gets $\fb_{\beta_d} \cong \fsl(\beta)$ and $\fb_{\gamma_d} \cong \fsl(\gamma)$.
\end{proof}

\subsection{Proof of Theorem \ref{thm:main-theorem}(\ref{item:F(n)-osp(lam,mu)})}

Elements of the set $F(n)$ are of the form $(\lambda,\lambda^*)$ with $n$ even and $\lambda \vdash n/2$, and one has $(\nu,\lambda^*) \prec (\lambda,\lambda^*)$ if and only if $(\lambda,\nu^*) \prec (\lambda,\lambda^*)$. One observes:

\begin{lemma} \label{lemma:(lam-lam*)-F(n)-preceding}
Let $n \geq 6$, let $(\lambda,\lambda^*) \in F(n)$, and let $\nu \prec \lambda$. Then $(\nu,\lambda^*)$ and $(\lambda,\nu^*)$ are proper, are not arm and leg bipartitions, and $(\nu,\lambda^*) \neq (\lambda,\nu^*)$. Thus $(\nu,\lambda^*) \in E(n-1)$.
\end{lemma}

\begin{lemma} \label{lemma:proof-of-item-F(n)-osp(lam,mu)}
Let $n \geq 6$, and suppose Theorem \ref{thm:main-theorem} is true for the value $n-1$. Then part \eqref{item:F(n)-osp(lam,mu)} of Theorem \ref{thm:main-theorem} is true for the value $n$, i.e.,
	\begin{center}
	If $(\lambda,\mu) \in F(n)$, then $\fb_{(\lambda,\mu)} = \osp(\lambda,\mu)$ as described in Lemma \ref{lemma:b-osp-inclusion}.
	\end{center}
\end{lemma}

\begin{proof}
Let $(\lambda,\lambda^*) \in F(n)$, and set $\fr = \fb_{n-1}$. Then by \eqref{eq:restriction-Bn-toBn-1}, the $\fb_n'$-module structure map $\rho_{(\lambda,\lambda^*)}: \fb_n' \twoheadrightarrow \fb_{(\lambda,\lambda^*)} \subseteq \osp(\lambda,\lambda^*) \subseteq \End(S^{(\lambda,\lambda^*)})$ restricts to a Lie algebra homomorphism
	\[
	\fr' \to \bigoplus_{\nu \prec \lambda} \Big[ \fr_{(\nu,\lambda^*)} \oplus \fr_{(\lambda,\nu^*)} \Big],
	\]
which by the discussion of Section \ref{subsubsec:dual-pairs} factors through the map
	\begin{equation} \label{eq:r-F(n)-sl(lam,lam*)-restricted}
	\rho: \fr' \to \bigoplus_{\nu \prec \lambda} \fr_{(\nu,\lambda^*)}.
	\end{equation}
It follows from Lemma \ref{lemma:(lam-lam*)-F(n)-preceding} that the right side of \eqref{eq:r-F(n)-sl(lam,lam*)-restricted} is a summand of the right side of \eqref{eq:bn-refined-factorization} for the value $n-1$. This implies by induction that \eqref{eq:r-F(n)-sl(lam,lam*)-restricted} is a surjection and $\fr_{(\nu,\lambda^*)} \cong \fsl(\nu,\lambda^*)$ for each $\nu \prec \lambda$. By the discussion of Section \ref{subsubsec:dual-pairs}, $S^{(\nu,\lambda^*)}$ and $S^{(\lambda,\nu^*)}$ are dual as $\fr'$-modules. Since $\dim(S^{(\nu,\lambda^*)}) \geq n-2 \geq 4$ by Lemmas \ref{lemma:(lam-lam*)-F(n)-preceding} and \ref{lemma:Bn-subminimal-dimensions}, and since the natural module for $\fsl(k)$ is self-dual only if $k=2$, it follows that $S^{(\nu,\lambda^*)} \not\cong S^{(\lambda,\nu^*)}$ as $\fr'$-modules.

Set $\fh = \rho_{(\lambda,\lambda^*)}(\fr')$. Then $\fh \cong \rho(\fr')$ is a semisimple subalgebra of $\fb_{(\lambda,\lambda^*)}$ such that the restriction of $S^{(\lambda,\lambda^*)}$ to each simple ideal of $\fh$ admits a simple factor of multiplicity one. Recall that $r_\lambda$ is the number of removable boxes in the Young diagram of $\lambda$. Then $r_\lambda < \abs{\lambda} = \frac{n}{2}$, and
	\[
	2 \cdot \rk(\fh) 
	= \sum_{\nu \prec \lambda} 2 \cdot \rk(\fsl(\nu,\lambda^*))
	=  \dim(S^{(\lambda,\lambda^*)}) - 2 \cdot r_\lambda > \dim(S^{(\lambda,\lambda^*)}) - n.
	\]
Recall that $f^\lambda = \dim(S^\lambda)$. 
Since $(\lambda,\lambda^*)$ is not an arm and leg by assumption, then $f^\lambda = f^{\lambda^*} \geq 2$. Now applying \eqref{eq:dim-S(lambda,mu)}, one gets
	\begin{equation} \label{eq:half-dim-S(lam,lam*)} \textstyle
	\frac{1}{2} \cdot \dim(S^{(\lambda,\lambda^*)}) 
	= \frac{1}{2} \cdot \binom{n}{n/2} \cdot f^\lambda \cdot f^{\lambda^*} 
	\geq 2 \cdot \binom{n}{n/2} 
	\geq 2n
	\end{equation}
and hence
	\begin{align*} 
	2 \cdot \rk(\fh) &> \textstyle \big( \tfrac{1}{2} \cdot \dim( S^{(\lambda,\lambda^*)}) \big) + \big( \tfrac{1}{2} \cdot \dim( S^{(\lambda,\lambda^*)}) \big) - n \\
	&\geq (\tfrac{1}{2} \cdot \dim( S^{(\lambda,\lambda^*)})) + n \\
	&> \tfrac{1}{2} \cdot \dim( S^{(\lambda,\lambda^*)}) = \rk(\osp(\lambda,\lambda^*)) \geq \rk(\fb_{(\lambda,\lambda^*)}).
	\end{align*}
Then $\fb_{(\lambda,\lambda^*)}$ is a simple Lie algebra by \cite[Lemme 15]{Marin:2007}. Since $\rk(\fh) \leq \rk(\fb_{(\lambda,\lambda^*)})$, the preceding inequalities also imply that $2 \cdot \rk(\fb_{(\lambda,\lambda^*)}) \leq \dim(S^{(\lambda,\lambda^*)}) < 4 \cdot \rk(\fb_{(\lambda,\lambda^*)})$. Now since $\dim(S^{(\lambda,\lambda^*)}) \geq \binom{n}{n/2}$, it follows from \cite[Lemme 14]{Marin:2007} that $\fb_{(\lambda,\lambda^*)} = \osp(\lambda,\lambda^*)$.
\end{proof}

\subsection{Proof of Theorem \ref{thm:main-theorem}(\ref{item:E(n)-sl(lam,mu)})} \label{subsec:E(n)-sl(lam,mu)}

The observations in the next lemma are straightforward.

\begin{lemma} \label{lemma:(lam-mu)-E(n)-preceding}
Let $n \geq 6$ and let $(\lambda,\mu) \in E(n)$.
	\begin{enumerate}
	\item Suppose $\abs{\lambda} \geq \abs{\mu}$. Then there exists an improper $\set{\nu,\tau} \prec \set{\lambda,\mu}$ if and only if $\mu = [1]$. If $\mu = [1]$, then $\set{\lambda,\emptyset}$ is the unique improper element such that $\set{\lambda,\emptyset} \prec \set{\lambda,\mu}$.

	\item Suppose $(\nu,\tau)$ is an arm and leg and $(\nu,\tau) \prec (\lambda,\mu)$. Then either
		\[
		\set{\lambda,\mu} = \{ [n-1-i,1],[1^i] \} \quad \text{or} \quad \set{\lambda,\mu} = \{ [n-1-i],[2,1^{i-1}] \}
		\]
	for some $1 \leq i \leq n-2$, and $(\nu,\tau)$ is the unique arm and leg such that $(\nu,\tau) \prec (\lambda,\mu)$.

	\item Suppose $(\nu,\tau) \prec (\lambda,\mu)$ and $(\tau^*,\nu^*) \prec (\lambda,\mu)$. Then either:
		\begin{enumerate}
		\item $(\nu,\tau) = (\mu^*,\mu)$ and $\mu^* \prec \lambda$, or
		\item $(\nu,\tau) = (\lambda,\lambda^*)$ and $\lambda^* \prec \mu$.
		\end{enumerate}
	Then $(\nu,\tau) = (\tau^*,\nu^*)$, and $(\nu,\tau)$ is unique such that $(\nu,\tau) \prec (\lambda,\mu)$ and $(\tau^*,\nu^*) \prec (\lambda,\mu)$.
	\end{enumerate}
\end{lemma}

\begin{lemma} \label{lemma:proof-of-item-E(n)-sl(lam,mu)}
Let $n \geq 6$, and suppose Theorem \ref{thm:main-theorem} is true for the value $n-1$. Then part \eqref{item:E(n)-sl(lam,mu)} of Theorem \ref{thm:main-theorem} is true for the value $n$, i.e.,
	\begin{center}
	If $(\lambda,\mu) \in E(n)$, then $\fb_{(\lambda,\mu)} = \fsl(\lambda,\mu)$.
	\end{center}
\end{lemma}

The proof of Lemma \ref{lemma:proof-of-item-E(n)-sl(lam,mu)} will occupy the rest of this section. Let $(\lambda,\mu) \in E(n)$, and let $\fr = \fb_{n-1}$. By \eqref{eq:restriction-Bn-toBn-1}, we have the $\fr$-module identification
	\begin{equation} \label{eq:S(lam,mu)-as-r-module}
	S^{(\lambda,\mu)} = \bigoplus_{(\nu,\tau) \prec (\lambda,\mu)} S^{(\nu,\tau)}.
	\end{equation}
Then $\rho_{(\lambda,\mu)}: \fb_n' \twoheadrightarrow \fb_{(\lambda,\mu)} \subseteq \fsl(\lambda,\mu)$ restricts to a Lie algebra homo\-morphism
	\begin{equation} \label{eq:r-E(n)-sl(lam,mu)}
	\rho: \fr' \to \bigoplus_{(\nu,\tau) \prec (\lambda,\mu)} \fr_{(\nu,\tau)}.
	\end{equation}
It follows from Lemma \ref{lemma:(lam-mu)-E(n)-preceding} that the right side of \eqref{eq:r-E(n)-sl(lam,mu)} is a summand of the right side of \eqref{eq:bn-refined-factorization} for the value $n-1$, except perhaps that a summand of the form $\fr_{\beta_{n-1,1}}$ or $\fr_{\gamma_{n-1,1}}$ may be replaced by an isomorphic summand of the form $\fr_{\beta_{n-1,d}}$ or $\fr_{\gamma_{n-1,d}}$ for some $1 \leq d \leq n-2$, and a summand of the form $\fs_{\alpha_{n-1,1}}$ may be replaced by an isomorphic summand $\fs_{\alpha_{n-1,d}}$ for some $1 \leq d \leq n-3$. This implies by the induction hypothesis that \eqref{eq:r-E(n)-sl(lam,mu)} is a surjection, and implies for each $(\nu,\tau) \prec (\lambda,\mu)$ that $\fr_{(\nu,\tau)}$ is a simple Lie algebra acting irreducibly on $S^{(\nu,\tau)}$.

Let $\fh = \rho_{(\lambda,\mu)}(\fr')$. Then $\fh \cong \rho(\fr')$ is a semisimple subalgebra of $\fb_{(\lambda,\mu)}$ such that the restriction of $S^{(\lambda,\mu)}$ to each simple ideal of $\fh$ admits a simple factor of multiplicity one. Since $\fh \subseteq \fb_{(\lambda,\mu)} \subseteq \fsl(\lambda,\mu)$, then $\rk(\fh) \leq \rk(\fb_{(\lambda,\mu)}) \leq \dim(S^{(\lambda,\mu)})$. Our goal now, as in the proof of Lemma \ref{lemma:proof-of-item-beta-gamma}, is to show that $2 \cdot \rk(\fh) > \dim(S^{(\lambda,\mu)})$. This will imply first by \cite[Lemme 15]{Marin:2007} that $\fb_{(\lambda,\mu)}$ is a simple Lie algebra, and then by \cite[Lemme 13]{Marin:2007} that $\fb_{(\lambda,\mu)} = \fsl(\lambda,\mu)$.

We break the remainder of the proof into the following cases:
	\begin{enumerate}
	\item \label{item:r(nu,tau)=sl(nu,tau)} $\fr_{(\nu,\tau)} \cong \fsl(\nu,\tau)$ for each $(\nu,\tau) \prec (\lambda,\mu)$.
	
	\item \label{item:(nu,tau)-improper-lambda-hook} $\mu = [1]$ and $\lambda$ is a hook partition.

	\item \label{item:(nu,tau)-improper-lambda-symmetric} $\mu = [1]$, $\lambda$ is not a hook, and $\lambda = \lambda^*$.
	
	\item \label{item:(nu,tau)-arm-and-leg} There exists $(\nu,\tau) \prec (\lambda,\mu)$ such that $(\nu,\tau)$ is an arm and leg.

	\item \label{item:(nu,tau)-self-dual} There exists $(\nu,\tau) \prec (\lambda,\mu)$ such that $(\nu,\tau) = (\tau^*,\nu^*)$, and $(\eta,\sigma)$ is not an arm and leg for all $(\eta,\sigma) \prec (\lambda,\mu)$.
	\end{enumerate}
For \eqref{item:(nu,tau)-improper-lambda-hook} and \eqref{item:(nu,tau)-improper-lambda-symmetric}, we may assume that $\mu = [1]$ by way of the isomorphism $\fb_{(\lambda,\mu)} \cong \fb_{(\mu^*,\lambda^*)}$.

\subsubsection{Completion of Lemma \ref{lemma:proof-of-item-E(n)-sl(lam,mu)} in case (\ref{item:r(nu,tau)=sl(nu,tau)})} \label{subsubsec:r(nu,tau)=sl(nu,tau)}

Suppose $\fr_{(\nu,\tau)} \cong \fsl(\nu,\tau)$ for all $(\nu,\tau) \prec (\lambda,\mu)$. Then
	\begin{align*}
	2 \cdot \rk(\fh) - \dim(S^{(\lambda,\mu)}) &= 2 \cdot \Big[ \sum_{(\nu,\tau) \prec (\lambda,\mu)}  (\dim(S^{(\nu,\tau)})-1) \Big] - \dim(S^{(\lambda,\mu)}) \\
	&= \sum_{(\nu,\tau) \prec (\lambda,\mu)} \big[ \dim(S^{(\nu,\tau)})-2 \big].
	\end{align*}
Since $(\lambda,\mu)$ is proper and is not an arm and leg, it follows from Lemma \ref{lemma:Bn-subminimal-dimensions} that $\dim(S^{(\nu,\tau)}) \geq 2$ for each $(\nu,\tau) \prec (\lambda,\mu)$, and there exists at least one $(\nu,\tau) \prec (\lambda,\mu)$ such that $\dim(S^{(\nu,\tau)}) > 2$. Then $2 \cdot \rk(\fh) > \dim(S^{(\lambda,\mu)})$, as desired.

\subsubsection{Completion of Lemma \ref{lemma:proof-of-item-E(n)-sl(lam,mu)} in case (\ref{item:(nu,tau)-improper-lambda-hook})} \label{subsubsec:(nu,tau)-improper-lambda-hook}

Suppose $\mu = [1]$ and $\lambda$ is a hook partition. Then $\lambda = \alpha_{n-1,d} = [n-1-d,1^d]$ for some $1 \leq d \leq n-3$; the values $d = 0$ and $d = n-2$ are excluded by the assumption that $(\lambda,\mu)$ is not an arm and leg. In this case, \eqref{eq:S(lam,mu)-as-r-module} becomes
	\[
	S^{(\lambda,\mu)} = S^{(\alpha_{n-2,d},[1])} \oplus S^{(\alpha_{n-2,d-1},[1])} \oplus S^{(\alpha_{n-1,d},\emptyset)},
	\]
By \eqref{eq:type-B-module-construction}, the dimensions of the modules are
	\begin{align*}
	\dim(S^{(\lambda,\mu)}) &= \textstyle n \cdot \binom{n-2}{d}, &
	\dim(S^{(\alpha_{n-2,d},[1])}) &= \textstyle (n-1) \cdot \binom{n-3}{d}, \\
	\dim(S^{(\alpha_{n-2,d-1},[1])}) &= \textstyle (n-1) \cdot \binom{n-3}{d-1}, &
	\dim(S^{(\alpha_{n-1,d},\emptyset)}) &= \textstyle \binom{n-2}{d}.
	\end{align*}
From Section \ref{subsubsec:improper-bipartitions}, $\fr_{(\alpha_{n-1,d},\emptyset)} \cong \fsl(\alpha_{n-1,1}) \cong \fsl(n-2)$. For $d \neq n-3$, the bipartition $(\alpha_{n-2,d},[1])$ is not an arm and leg, and hence $\fr_{(\alpha_{n-2,d},[1])} \cong \fsl(S^{(\alpha_{n-2,d},[1])})$, while for $d = n-3$, one has $(\alpha_{n-2,d},[1]) = ([1^{n-2}],[1]) = \gamma_{n-1,n-2}$. By induction, $\fr_{\gamma_{n-1,n-2}} \cong \fsl(\gamma_{n-1,1}) \cong \fsl(n-1)$, so by dimension comparison the isomorphism $\fr_{(\alpha_{n-2,d},[1])} \cong \fsl(S^{(\alpha_{n-2,d},[1])})$ is true for $d=n-3$. Similarly, one sees for all $1 \leq d \leq n-3$ that $\fr_{(\alpha_{n-2,d-1},[1])} \cong \fsl(S^{(\alpha_{n-2,d-1},[1])})$. Then
	\begin{align*}
	\rk(\fh) &= \textstyle \big[ (n-1) \cdot \binom{n-3}{d} - 1 \big] + \big[ (n-1) \cdot \binom{n-3}{d-1} - 1 \big] + \big[ (n-2) - 1 \big] \\
	&= \textstyle (n-1) \cdot \binom{n-2}{d} + (n-5)
	\end{align*}
and hence
	\[ \textstyle
	2 \cdot \rk(\fh) - \dim(S^{(\lambda,\mu)}) = (n-2) \cdot \binom{n-2}{d} + (2n-10) > 0.
	\]
Thus $2 \cdot \rk(\fh) > \dim(S^{(\lambda,\mu)})$, as desired.

\subsubsection{Completion of Lemma \ref{lemma:proof-of-item-E(n)-sl(lam,mu)} in case (\ref{item:(nu,tau)-improper-lambda-symmetric})}

Suppose $\mu = [1]$, $\lambda$ is not a hook partition, and $\lambda = \lambda^*$. Then \eqref{eq:S(lam,mu)-as-r-module} becomes $S^{(\lambda,[1])} = \left[ \bigoplus_{\nu \prec \lambda} S^{(\nu,[1])} \right] \oplus S^{(\lambda,\emptyset)}$, and the dimensions of the modules are $\dim(S^{(\lambda,[1])}) = n \cdot f^\lambda$, $\dim(S^{(\nu,[1])}) = (n-1) \cdot f^\nu$, and $\dim(S^{(\lambda,\emptyset)}) = f^\lambda$. For each $\nu \prec \lambda$, the bipartition $(\nu,[1])$ is proper and is not an arm and leg (because $\lambda$ is not a hook), so $\fr_{(\nu,[1])} \cong \fsl(S^{(\nu,[1])})$ by induction. Since $\lambda = \lambda^*$ and $\lambda$ is not a hook, then $\fr_{(\lambda,\emptyset)} \cong \osp(\lambda)$. Now
	\begin{align*}
	2 \cdot \rk(\fh) - \dim(S^{(\lambda,[1])}) &= 2 \cdot \Big[ \sum_{\nu \prec \lambda} (n-1) \cdot f^\nu - 1 \Big] + f^\lambda - n \cdot f^\lambda \\
	&=\sum_{\nu \prec \lambda} \left[ (n-1) \cdot f^\nu - 2 \right]
	\end{align*}
Since $(n-1) \cdot f^\nu \geq 5 \cdot 1$ for each $\nu \prec \lambda$, this implies that $2 \cdot \rk(\fh) > \dim(S^{(\lambda,[1])})$, as desired.

\subsubsection{Completion of Lemma \ref{lemma:proof-of-item-E(n)-sl(lam,mu)} in case (\ref{item:(nu,tau)-arm-and-leg})} \label{subsubsec:(nu,tau)-arm-and-leg}

Suppose there exists $(\nu,\tau) \prec (\lambda,\mu)$ such that $(\nu,\tau)$ is an arm and leg. First suppose $(\lambda,\mu) = ([n-1-i,1],[1^i])$ for some $1 \leq i \leq n-2$. Then
	\[
	S^{(\lambda,\mu)} = S^{([n-1-i],[1^i])} \oplus S^{([n-2-i,1],[1^i])} \oplus S^{([n-1-i,1],[1^{i-1}])},
	\]
where the second term on the right is zero if $i = n-2$. The dimensions of the modules are
	\begin{align*}
	\dim(S^{(\lambda,\mu)}) &= \textstyle \binom{n}{i} \cdot (n-1-i), &
	\dim(S^{([n-1-i],[1^i])}) &= \textstyle \binom{n-1}{i}, \\
	\dim(S^{([n-2-i,1],[1^i])}) &= \textstyle \binom{n-1}{i} \cdot (n-2-i), &
	\dim(S^{([n-1-i,1],[1^{i-1}])}) &= \textstyle \binom{n-1}{i-1} \cdot (n-1-i).
	\end{align*}
One has $([n-1-i],[1^i]) = \beta_{n-1,i}$, and $\fr_{\beta_{n-1,i}} \cong \fr_{\beta_{n-1,1}} \cong \fsl(n-1)$ by the induction hypothesis, while for the other $(\nu,\tau) \prec (\lambda,\mu)$ one has $\fr_{(\nu,\tau)} \cong \fsl(S^{(\nu,\tau)})$. For $1 \leq i \leq n-3$, this implies that
	\begin{align*}
	\rk(\fh) &= \textstyle (n-2) + \big[ \binom{n-1}{i} \cdot (n-2-i) - 1 \big] + \big[ \binom{n-1}{i-1} \cdot (n-1-i)-1 \big] \\
	&= \textstyle (n-4) + \binom{n}{i} \cdot (n-2-i) + \binom{n-1}{i-1},
	\end{align*}
and hence
	\[ \textstyle
	2 \cdot \rk(\fh) - \dim(S^{(\lambda,\mu)}) = (2n-8) + \binom{n}{i} \cdot (n-3-i) + 2 \cdot \binom{n-1}{i-1} > 0.
	\]
On the other hand, for $i = n-2$ one has
	\[ \textstyle
	\rk(\fh) = \big[ \binom{n-1}{1}-1 \big] + \big[ \binom{n-1}{n-3} - 1 \big] = \binom{n-1}{1} + \binom{n-1}{2} -2 = \binom{n}{2} - 2.
	\]
Then $2 \cdot \rk(\fh) - \dim(S^{(\lambda,\mu)}) = \binom{n}{2} - 4 > 0$. Thus, if $(\lambda,\mu) = ([n-1-i,1],[1^i])$ for some $1 \leq i \leq n-2$, then $2 \cdot \rk(\fh) > \dim(S^{(\lambda,\mu)})$, as desired.

Now by Lemma \ref{lemma:(lam-mu)-E(n)-preceding}, the are three other possibilities for $(\lambda,\mu)$. The proof in the case $(\lambda,\mu) = ([1^i],[n-1-i,1])$ is entirely similar to that given above, but with the parts of the bipartitions swapped and with $\beta_{n-1,i}$ replaced by $\gamma_{n-1,i}$; neither of these changes affect the dimension or rank calculations. Then the remaining cases are deduced via the isomorphism $\fb_{(\lambda,\mu)} \cong \fb_{(\mu^*,\lambda^*)}$.

\subsubsection{Completion of Lemma \ref{lemma:proof-of-item-E(n)-sl(lam,mu)} in case (\ref{item:(nu,tau)-self-dual})}

Suppose $(\eta,\sigma)$ is not an arm and leg for all $(\eta,\sigma) \prec (\lambda,\mu)$, and there exists $(\nu,\tau) \prec (\lambda,\mu)$ such that $(\nu,\tau) = (\tau^*,\nu^*)$. Then by Lemma \ref{lemma:(lam-mu)-E(n)-preceding}, either $\lambda^* \prec \mu$ or $\mu^* \prec \lambda$; we will assume that $\lambda^* \prec \mu$, the other case being entirely parallel. Then
	\[
	S^{(\lambda,\mu)} = \Big[ \bigoplus_{\nu \prec \lambda} S^{(\nu,\mu)} \Big] \oplus \Big[ \bigoplus_{\substack{\tau \prec \mu \\ \tau \neq \lambda^*}} S^{(\lambda,\tau)} \Big] \oplus S^{(\lambda,\lambda^*)}.
	\]
By induction, $\fr_{(\lambda,\lambda^*)} \cong \osp(\lambda,\lambda^*)$, and $\fr_{(\nu,\tau)} \cong \fsl(\nu,\tau)$ for all other $(\nu,\tau) \prec (\lambda,\mu)$. Then
	\[
	2 \cdot \rk(\fh) = 2 \cdot \Big[ \sum_{\nu \prec \lambda} \dim(S^{(\nu,\lambda)})-1 \Big] + 2 \cdot \Big[ \sum_{\substack{\tau \prec \mu \\ \tau \neq \lambda^*}} \dim(S^{(\lambda,\tau)})-1 \Big] + \dim(S^{(\lambda,\lambda^*)})
	\]
and hence
	\[
	2 \cdot \rk(\fh) - \dim(S^{(\lambda,\mu)}) = \Big[ \sum_{\nu \prec \lambda} \dim(S^{(\nu,\lambda)})-2 \Big] + \Big[ \sum_{\substack{\tau \prec \mu \\ \tau \neq \lambda^*}} \dim(S^{(\lambda,\tau)})-2 \Big].
	\]
As in Section \ref{subsubsec:r(nu,tau)=sl(nu,tau)}, one then deduces that $2 \cdot \rk(\fh) > \dim(S^{(\lambda,\mu)})$, as desired.

\subsection{Proof of Theorem \ref{thm:main-theorem}(\ref{item:Fn-osp(lam,mu)})} \label{subsec:Fn-osp(lam,mu)}

The observations in the next lemma are straightforward.

\begin{lemma} \label{lemma:(lam-mu)-F(n)-preceding}
Let $n \geq 6$, and let $\set{\lambda,\mu} \in F\{n\}$ with $\lambda \neq \mu$.
	\begin{enumerate}
	\item \label{item:d-F(n)-(lam-mu)=(mu*-lam*)} Suppose $(\lambda,\mu) = (\mu^*,\lambda^*)$, i.e., $\mu = \lambda^*$. Then for all $\nu \prec \lambda$, the elements $\set{\nu,\lambda^*}$ and $\set{\lambda,\nu^*}$ are proper and are not arm and legs, and hence are elements of $E\{n-1\}$.
	
	\item Suppose $(\lambda,\mu) = (\lambda^*,\mu^*)$ and $\abs{\lambda} \geq \abs{\mu}$. Then:
		\begin{enumerate}
		\item For all $\set{\nu,\tau} \prec \set{\lambda,\mu}$, the element $\set{\nu,\tau}$ is not an arm and leg.

		\item If there exists an improper $\set{\nu,\tau} \prec \set{\lambda,\mu}$, then $\mu = [1]$ and $\set{\nu,\tau} = \set{\lambda,\emptyset}$.
		
		\item If there exists $\set{\nu,\tau} \prec \set{\lambda,\mu}$ with $\nu = \tau$, then $\mu \prec \lambda$ and $\set{\nu,\tau} = \set{\mu,\mu}$.
		\end{enumerate}
	\end{enumerate}
\end{lemma}

\begin{lemma} \label{lemma:proof-of-item-Fn-osp(lam,mu)}
Let $n \geq 6$, and suppose Theorem \ref{thm:main-theorem} is true for the value $n-1$. Then part \eqref{item:Fn-osp(lam,mu)} of Theorem \ref{thm:main-theorem} is true for the value $n$, i.e.,
	\begin{center}
	If $\set{\lambda,\mu} \in F\{n\}$ and $\lambda \neq \mu$, then $\fd_{\set{\lambda,\mu}} = \osp\{\lambda,\mu\}$ as described in Lemma \ref{lemma:b-osp-inclusion} or \ref{lemma:d-osp(lam,mu)-inclusion}.
	\end{center}
\end{lemma}

The proof of Lemma \ref{lemma:proof-of-item-Fn-osp(lam,mu)} will occupy the rest of this section. If $(\lambda,\mu) = (\mu^*,\lambda^*)$, then the proof proceeds in an entirely parallel manner to the proof of Lemma \ref{lemma:proof-of-item-F(n)-osp(lam,mu)}, using Lemma \ref{lemma:(lam-mu)-F(n)-preceding}\eqref{item:d-F(n)-(lam-mu)=(mu*-lam*)} in lieu of Lemma \ref{lemma:(lam-lam*)-F(n)-preceding}, so suppose that $(\lambda,\mu) = (\lambda^*,\mu^*)$ and $\abs{\lambda} \geq \abs{\mu}$. Let $\fr = \fd_{n-1}$. We break the remainder of the proof into the following cases:
	\begin{enumerate}
	\item \label{item:mu=1-lam-hook} $\mu = [1]$ and $\lambda$ is a hook partition.
	\item \label{item:mu=1-lam-non-hook} $\mu = [1]$ and $\lambda$ is not a hook partition.
	\item \label{item:mu-neq-1-mu-prec-lam} $\mu \neq [1]$ and $\mu \prec \lambda$.
	\item \label{item:mu-neq-1-mu-nprec-lam} $\mu \neq [1]$ and $\mu \not\prec \lambda$.
	\end{enumerate}

\subsubsection{Completion of Lemma \ref{lemma:proof-of-item-Fn-osp(lam,mu)} in case (\ref{item:mu=1-lam-hook})}

Suppose $\mu = [1]$ and $\lambda$ is a hook partition. Since $\lambda = \lambda^*$, this implies that $n$ is even, say $n = 2m$ (so $m \geq 3$ because $n \geq 6$), and $\lambda = [m,1^{m-1}]$. Then as an $\fr$-module,
	\[
	S^{\{[m,1^{m-1}],[1] \}} = S^{\{[m-1,1^{m-1}],[1] \}} \oplus S^{\{[m,1^{m-2}],[1] \}} \oplus S^{\{ [m,1^{m-1}],\emptyset \}}.
	\]
The dimensions of the modules are 
	\begin{align*}
	\dim(S^{\{[m,1^{m-1}],[1] \}}) &= \textstyle n \cdot \binom{n-2}{m-1}, &
	\dim(S^{\{[m-1,1^{m-1}],[1] \}}) &= \textstyle (n-1) \cdot \binom{n-3}{m-1}, \\
	\dim(S^{\{[m,1^{m-2}],[1] \}}) &= \textstyle (n-1) \cdot \binom{n-3}{m-2}, &
	\dim(S^{\{ [m,1^{m-1}],\emptyset \}}) &= \textstyle \binom{n-2}{m-1}.
	\end{align*}
The module structure map $\rho_{\set{\lambda,\mu}}: \fd_n' \twoheadrightarrow \fd_{\set{\lambda,\mu}} \subseteq \osp\{\lambda,\mu\} \subseteq \End(S^{\set{\lambda,\mu}})$ restricts to a Lie algebra homomorphism $\fr' \to \fr_{\{[m-1,1^{m-1}],[1] \}} \oplus \fr_{\{[m,1^{m-2}],[1] \}} \oplus \fr_{\{ [m,1^{m-1}],\emptyset \}}$, which by the discussion of Section \ref{subsubsec:dual-pairs} factors through the map
	\begin{equation} \label{eq:r-F(n)-osp(lam,mu)-case-1-restricted}
	\rho: \fr' \to \fr_{\{[m-1,1^{m-1}],[1] \}} \oplus \fr_{\{ [m,1^{m-1}],\emptyset \}}.
	\end{equation}
The right side of \eqref{eq:r-F(n)-osp(lam,mu)-case-1-restricted} is a summand of the right side of \eqref{eq:dn-refined-factorization} for the value $n-1$. This implies by the induction hypothesis that \eqref{eq:r-F(n)-osp(lam,mu)-case-1-restricted} is a surjection. By induction, one has
	\[
	\fr_{\{[m-1,1^{m-1}],[1] \}} \cong \fsl\{[m-1,1^{m-1}],[1] \} \qquad \text{and} \qquad \fr_{\{ [m,1^{m-1}],\emptyset \}} \cong \fsl(n-2).
	\]

Let $\fh = \rho_{\set{\lambda,\mu}}(\fr') \cong \rho(\fr')$. Then $S^{\{[m-1,1^{m-1}],[1] \}}$ and $S^{\{[m,1^{m-2}],[1] \}}$ are dual but not isomorphic as $\fh$-modules, because the natural module for $\fsl(k)$ is self-dual only if $k=2$. Then $\fh$ is a semisimple subalgebra of $\fd_{\set{\lambda,\mu}}$ such that the restriction of $S^{\set{\lambda,\mu}}$ to each simple ideal of $\fh$ admits a simple factor of multiplicity one. One has
	\[
	2 \cdot \rk(\fh) = \textstyle \big[ 2 \cdot (n-1) \cdot \binom{n-3}{m-1} - 1 \big] + 2 \cdot (n-3) = (n-1) \cdot \binom{n-2}{m-1} + (2n-6).
	\]
Since $\fd_{\set{\lambda,\mu}} \subseteq \osp\{\lambda,\mu\}$, then $\rk(\fd_{\set{\lambda,\mu}}) \leq \frac{1}{2} \cdot \dim(S^{\set{\lambda,\mu}}) = \frac{n}{2} \cdot \binom{n-2}{m-1}$, and hence
	\[ \textstyle
	2 \cdot \rk(\fh) - \rk(\fd_{\set{\lambda,\mu}}) \geq (\frac{n}{2}-1) \cdot \binom{n-2}{m-1} + (2n-6) > 0.
	\]
Then $\fd_{\set{\lambda,\mu}}$ is a simple Lie algebra by \cite[Lemme 15]{Marin:2007}. Since $\rk(\fh) \leq \rk(\fd_{\set{\lambda,\mu}})$, the preceding inequalities also imply that $2 \cdot \rk(\fd_{\set{\lambda,\mu}}) \leq \dim(S^{\set{\lambda,\mu}}) < 4 \cdot \rk(\fd_{\set{\lambda,\mu}})$. Then since $n \geq 6$, it follows from \cite[Lemme 14]{Marin:2007} that $\fd_{\set{\lambda,\mu}} = \osp\{\lambda,\mu\}$.

\subsubsection{Completion of Lemma \ref{lemma:proof-of-item-Fn-osp(lam,mu)} in case (\ref{item:mu=1-lam-non-hook})}

Suppose $\mu = [1]$ and $\lambda$ is not a hook partition. Then as an $\fr$-module,
	\[
	S^{\set{\lambda,[1]}} = \Big[ \bigoplus_{\substack{\nu \prec \lambda \\ \res(\lambda/\nu) > 0}} S^{\set{\nu,[1]}} \oplus S^{\set{\nu^*,[1]}} \Big] \oplus \Big[ \bigoplus_{\substack{\nu \prec \lambda \\ \res(\lambda/\nu) = 0}} S^{\set{\nu,[1]}} \Big] \oplus S^{\set{\lambda,\emptyset}}.
	\]
There is at most one $\nu \prec \lambda$ such that $\res(\lambda/\nu) = 0$; if such a $\nu$ exists then $\set{\nu,[1]} \in F\{n-1\}$ with $\nu \neq [1]$. If $\nu \prec \lambda$ and $\res(\lambda/\nu) > 0$, then $\set{\nu,[1]} \in E\{n-1\}$ with $\nu \neq [1]$. Now the module structure map $\rho_{\set{\lambda,[1]}}: \fd_n' \twoheadrightarrow \fd_{\set{\lambda,[1]}} \subseteq \osp\{\lambda,[1]\} \subseteq \End(S^{\set{\lambda,[1]}})$ restricts to a Lie algebra homomorphism
	\[
	\fr' \to \Big[ \bigoplus_{\substack{\nu \prec \lambda \\ \res(\lambda/\nu) > 0}} \fr_{\set{\nu,[1]}} \oplus \fr_{\set{\nu^*,[1]}} \Big] \oplus \Big[ \bigoplus_{\substack{\nu \prec \lambda \\ \res(\lambda/\nu) = 0}} \fr_{\set{\nu,[1]}} \Big] \oplus \fr_{\set{\lambda,\emptyset}},
	\]
which by the discussion of Section \ref{subsubsec:dual-pairs} factors through the map
	\begin{equation} \label{eq:r-F(n)-osp(lam,mu)-case-2-restricted}
	\rho: \fr' \to \Big[ \bigoplus_{\substack{\nu \prec \lambda \\ \res(\lambda/\nu) > 0}} \fr_{\set{\nu,[1]}} \Big] \oplus \Big[ \bigoplus_{\substack{\nu \prec \lambda \\ \res(\lambda/\nu) = 0}} \fr_{\set{\nu,[1]}} \Big] \oplus \fr_{\set{\lambda,\emptyset}}.
	\end{equation}
The right side of \eqref{eq:r-F(n)-osp(lam,mu)-case-2-restricted} is a summand of the right side of \eqref{eq:dn-refined-factorization} for the value $n-1$. This implies by the induction hypothesis that \eqref{eq:r-F(n)-osp(lam,mu)-case-2-restricted} is a surjection. By induction, one has
	\begin{align*}
	\fr_{\set{\nu,[1]}} &\cong \fsl\{\nu,[1]\} & \text{if $\nu \prec \lambda$ and $\res(\lambda/\nu) > 0$,} \\
	\fr_{\set{\nu,[1]}} &\cong \osp\{\nu,[1]\} & \text{if $\nu \prec \lambda$ and $\res(\lambda/\nu) = 0$,} \\
	\fr_{\set{\lambda,\emptyset}} &\cong \osp(\lambda).
	\end{align*}

Let $\fh = \rho_{\set{\lambda,[1]}}(\fr') \cong \rho(\fr')$. Then for $\nu \prec \lambda$ with $\res(\lambda/\nu) > 0$, the modules $S^{\set{\nu,[1]}}$ and $S^{\set{\nu^*,[1]}}$ are dual but not isomorphic as $\fh$-modules, because $\dim(S^{\set{\nu,[1]}}) = (n-1) \cdot f^\nu \geq 5$, and the natural module for $\fsl(k)$ is self-dual only if $k=2$. Now $\fh$ is a semisimple subalgebra of $\fd_{\set{\lambda,[1]}}$ such that the restriction of $S^{\set{\lambda,[1]}}$ to each simple ideal of $\fh$ admits a simple factor of multiplicity one, and
	\begin{align*}
	2 \cdot \rk(\fh) &= \Big[ \sum_{\substack{\nu \prec \lambda \\ \res(\lambda/\nu) > 0}} 2 \cdot \dim(S^{\set{\nu,[1]}}) -1 \Big] + \Big[ \sum_{\substack{\nu \prec \lambda \\ \res(\lambda/\nu) = 0}} \dim(S^{\set{\nu,[1]}}) \Big] + \dim(S^{\set{\lambda,\emptyset}}) \\
	&\geq \dim(S^{\set{\lambda,[1]}}) - r_\lambda \\
	&>\dim(S^{\set{\lambda,[1]}}) - (n-1).
	\end{align*}
Since $\fd_{\set{\lambda,[1]}} \subseteq \osp\{\lambda,[1]\}$, then $\rk(\fd_{\set{\lambda,[1]}}) \leq \frac{1}{2} \cdot \dim(S^{\set{\lambda,[1]}}) = \frac{n}{2} \cdot f^\lambda$. Also, since $\lambda = \lambda^*$ and $n-1 \geq 5$, then $f^\lambda \geq (n-1)-1 \geq 2$ by \cite[Theorem 2.4.10]{James:1981}. Then
	\[ \textstyle
	2 \cdot \rk(\fh) - \rk(\fd_{\set{\lambda,[1]}}) > \frac{n}{2} \cdot f^\lambda - (n-1) \geq n - (n-1) = 1 > 0,
	\]
so $\fd_{\set{\lambda,[1]}}$ is simple by \cite[Lemme 15]{Marin:2007}. Further, $2 \cdot \rk(\fd_{\set{\lambda,[1]}}) \leq \dim(S^{\set{\lambda,[1]}}) < 4 \cdot \rk(\fd_{\set{\lambda,[1]}})$, and then one deduces using \cite[Lemme 14]{Marin:2007} that $\fd_{\set{\lambda,[1]}} = \osp\{\lambda,[1]\}$.

\subsubsection{Completion of Lemma \ref{lemma:proof-of-item-Fn-osp(lam,mu)} in case (\ref{item:mu-neq-1-mu-prec-lam})}

Suppose $\mu \neq [1]$ and $\mu \prec \lambda$. Since $\lambda$ and $\mu$ are both symmetric, then $\mu$ is the unique partition such that $\res(\lambda/\mu) = 0$. As an $\fr$-module,
	\begin{align*}
	S^{\set{\lambda,\mu}} = \Big[ &\bigoplus_{\substack{\nu \prec \lambda \\ \res(\lambda/\nu) > 0}} S^{\set{\nu,\mu}} \oplus S^{\set{\nu^*,\mu}} \Big] \oplus \Big[ S^{\set{\mu,+}} \oplus S^{\set{\mu,-}} \Big] \\
	&\oplus \Big[ \bigoplus_{\substack{\tau \prec \mu \\ \res(\mu/\tau) > 0}} S^{\set{\lambda,\tau}} \oplus S^{\set{\lambda,\tau^*}} \Big] \oplus \Big[ \bigoplus_{\substack{\tau \prec \mu \\ \res(\mu/\tau) = 0}} S^{\set{\lambda,\tau}} \Big].
	\end{align*}
There is at most one $\tau \prec \mu$ such that $\res(\mu/\tau) = 0$; if such a $\tau$ exists then $\set{\lambda,\tau} \in F\{n-1\}$ with $\lambda \neq \tau$. If $\tau \prec \mu$ and $\res(\mu/\tau) > 0$, then $\set{\lambda,\tau} \in E\{n-1\}$ with $\lambda \neq \tau$. Analogous statements hold for the partitions $\nu \prec \lambda$ with $\res(\lambda/\nu) > 0$. Now $\rho_{\set{\lambda,\mu}}: \fd_n' \twoheadrightarrow \fd_{\set{\lambda,\mu}} \subseteq \osp\{\lambda,\mu\} \subseteq \End(S^{\set{\lambda,\mu}})$ restricts to a Lie algebra homomorphism
	\begin{align*}
	\fr' &\to \Big[ \bigoplus_{\substack{\nu \prec \lambda \\ \res(\lambda/\nu) > 0}} \fr_{\set{\nu,\mu}} \oplus \fr_{\set{\nu^*,\mu}} \Big] \oplus \Big[ \fr_{\set{\mu,+}} \oplus \fr_{\set{\mu,-}} \Big] \\
	&\mathrel{\phantom{=}} \oplus \Big[ \bigoplus_{\substack{\tau \prec \mu \\ \res(\mu/\tau) > 0}} \fr_{\set{\lambda,\tau}} \oplus \fr_{\set{\lambda,\tau^*}} \Big] \oplus \Big[ \bigoplus_{\substack{\tau \prec \mu \\ \res(\mu/\tau) = 0}} \fr_{\set{\lambda,\tau}} \Big],
	\end{align*}
which by the discussion of Section \ref{subsubsec:dual-pairs} factors through the map
	\begin{multline} \label{eq:r-F(n)-osp(lam,mu)-case-3-restricted}
	\rho: \fr' \to \\
	\Big[ \bigoplus_{\substack{\nu \prec \lambda \\ \res(\lambda/\nu) > 0}} \fr_{\set{\nu,\mu}} \Big] 
	\oplus \Big[ \fr_{\set{\mu,+}} \oplus \fr_{\set{\mu,-}}^\star \Big] 
	\oplus \Big[ \bigoplus_{\substack{\tau \prec \mu \\ \res(\mu/\tau) > 0}} \fr_{\set{\lambda,\tau}} 
	\oplus \Big] \oplus \Big[ \bigoplus_{\substack{\tau \prec \mu \\ \res(\mu/\tau) = 0}} \fr_{\set{\lambda,\tau}} \Big],
	\end{multline}
where the $\star$ indicates that the term $\fr_{\set{\mu,-}}$ is omitted if $\abs{\mu}$ is odd. Then the right side of \eqref{eq:r-F(n)-osp(lam,mu)-case-3-restricted} is a summand of the right side of \eqref{eq:dn-refined-factorization} for the value $n-1$. This implies by the induction hypothesis that \eqref{eq:r-F(n)-osp(lam,mu)-case-3-restricted} is a surjection. By induction, one has
	\begin{align*}
	\fr_{\set{\nu,\mu}} &\cong \fsl\{\nu,\mu\} & \text{if $\res(\lambda/\nu) > 0$,} \\
	\fr_{\set{\mu,\pm}} &\cong \osp\{\mu,\pm\} & \text{if $\abs{\mu}$ is even,} \\
	\fr_{\set{\mu,+}} &\cong \fsl\{\mu,+\} & \text{if $\abs{\mu}$ is odd,} \\
	\fr_{\set{\lambda,\tau}} &\cong \fsl\{\lambda,\tau \} & \text{if $\res(\mu/\tau) > 0$,} \\
	\fr_{\set{\lambda,\tau}} &\cong \osp\{\lambda,\tau \} & \text{if $\res(\mu/\tau) = 0$.} 
	\end{align*}

Let $\fh = \rho_{\set{\lambda,\mu}}(\fr') \cong \rho(\fr')$. For $\nu \prec \lambda$ with $\res(\lambda/\nu) > 0$, the modules $S^{\set{\nu,\mu}}$ and $S^{\set{\nu^*,\mu}}$ are dual but not isomorphic as $\fh$-modules, because $\dim(S^{\set{\nu,\mu}}) = \binom{n-1}{\abs{\mu}} \cdot f^\nu \cdot f^\mu \geq n-1 \geq 5$, and the natural module for $\fsl(k)$ is self-dual only if $k=2$. Corresponding statements apply to $S^{\set{\lambda,\tau}}$ and $S^{\set{\lambda,\tau^*}}$ if $\res(\mu/\tau) > 0$, while if $\abs{\mu}$ is odd, then $S^{\set{\mu,+}}$ and $S^{\set{\mu,-}}$ are dual but not isomorphic as $\fh$-modules, because $\dim(S^{\set{\mu,\pm}}) = \frac{1}{2} \cdot \binom{n-1}{\abs{\mu}} \cdot f^\mu \cdot f^\mu \geq \frac{1}{2} \cdot (n-1) > 2$.

Now $\fh$ is a semisimple subalgebra of $\fd_{\set{\lambda,\mu}}$ such that the restriction of $S^{\set{\lambda,\mu}}$ to each simple ideal of $\fh$ admits a simple factor of multi\-plicity one. Noting that $\rk(\osp\{\mu,\pm\}) = \frac{1}{4} \cdot \dim(S^{\set{\mu,\mu}})$ and $\rk(\fsl\{\mu,+\}) = \frac{1}{2} \cdot \dim(S^{\set{\mu,\mu}})-1$, one has
	\begin{align*}
	2 \cdot \rk(\fh) &\geq \Big[ \sum_{\substack{\nu \prec \lambda \\ \res(\lambda/\nu) > 0}} 2 \cdot (\dim(S^{\set{\nu,\mu}})-1) \Big] 
	+ \Big[ \dim(S^{\set{\mu,\mu}})-2 \Big] \\
	&\mathrel{\phantom{=}} + \Big[ \sum_{\substack{\tau \prec \mu \\ \res(\mu/\tau) > 0}} 2 \cdot (\dim(S^{\set{\lambda,\tau}})-1) \Big] 
	+ \Big[ \sum_{\substack{\tau \prec \mu \\ \res(\mu/\tau) = 0}} \dim(S^{\set{\lambda,\tau}}) \Big] \\
	&\geq \dim(S^{\set{\lambda,\mu}}) - r_\lambda - r_\mu -1 \\
	&> \dim(S^{\set{\lambda,\mu}}) - n -1,
	\end{align*}
where at the second inequality we note that the additional $1$ is subtracted because there may not be any $\tau \prec \mu$ such that $\res(\mu/\tau) = 0$ (and hence the $-2$ term on the first line results in subtracting one more than $r_\lambda + r_\mu$). Since $\fd_{\set{\lambda,\mu}} \subseteq \osp\{\lambda,\mu\}$, then $\rk(\fd_{\set{\lambda,\mu}}) \leq \frac{1}{2} \cdot \dim(S^{\set{\lambda,\mu}})$, and hence
	\[ \textstyle
	2 \cdot \rk(\fh) - \rk(\fd_{\set{\lambda,\mu}}) > \frac{1}{2} \cdot \dim(S^{\set{\lambda,\mu}}) - (n+1).
	\]
Since $\lambda$ and $\mu$ are symmetric and not equal to $[1]$, then $f^\lambda \cdot f^\mu \geq 2 \cdot 2$, and hence
	\[ \textstyle
	\dim(S^{\set{\lambda,\mu}}) = \binom{n}{\abs{\mu}} \cdot f^\lambda \cdot f^\mu \geq n \cdot f^\lambda \cdot f^\mu \geq 4n.
	\]
Then $2 \cdot \rk(\fh) - \rk(\fd_{\set{\lambda,\mu}}) \geq 2n - (n+1) = n-1 > 0$, so $\fd_{\set{\lambda,\mu}}$ is simple by \cite[Lemme 15]{Marin:2007}. Finally, since $\rk(\fh) \leq \rk(\fd_{\set{\lambda,\mu}})$, one sees that $2 \cdot \rk(\fd_{\set{\lambda,\mu}}) \leq \dim(S^{\set{\lambda,\mu}}) < 4 \cdot \rk(\fd_{\set{\lambda,\mu}})$, and then deduces using \cite[Lemme 14]{Marin:2007} that $\fd_{\set{\lambda,\mu}} = \osp\{\lambda,\mu\}$.

\subsubsection{Completion of Lemma \ref{lemma:proof-of-item-Fn-osp(lam,mu)} in case (\ref{item:mu-neq-1-mu-nprec-lam})}

Suppose $\mu \neq [1]$ and $\mu \not\prec \lambda$. Then as an $\fr$-module,
	\begin{align*}
	S^{\set{\lambda,\mu}} = \Big[ &\bigoplus_{\substack{\nu \prec \lambda \\ \res(\lambda/\nu) > 0}} S^{\set{\nu,\mu}} \oplus S^{\set{\nu^*,\mu}} \Big] \oplus \Big[ \bigoplus_{\substack{\nu \prec \lambda \\ \res(\lambda/\nu) = 0}} S^{\set{\nu,\mu}} \Big] \\
	&\oplus \Big[ \bigoplus_{\substack{\tau \prec \mu \\ \res(\mu/\tau) > 0}} S^{\set{\lambda,\tau}} \oplus S^{\set{\lambda,\tau^*}} \Big] \oplus \Big[ \bigoplus_{\substack{\tau \prec \mu \\ \res(\mu/\tau) = 0}} S^{\set{\lambda,\tau}} \Big].
	\end{align*}
There is at most one $\nu \prec \lambda$ such that $\res(\lambda/\nu) = 0$; if such a $\nu$ exists then $\set{\nu,\mu} \in F\{n-1\}$ with $\nu \neq \mu$ because $\mu \not \prec \lambda$. If $\nu \prec \lambda$ and $\res(\lambda/\nu) > 0$, then $\set{\nu,\mu} \in E\{n-1\}$ with $\nu \neq \mu$. Analogous statements hold for the partitions $\tau \prec \mu$. Now $\rho_{\set{\lambda,\mu}}: \fd_n' \twoheadrightarrow \fd_{\set{\lambda,\mu}} \subseteq \osp\{\lambda,\mu\} \subseteq \End(S^{\set{\lambda,\mu}})$ restricts to a Lie algebra homomorphism
	\begin{align*}
	\fr' &\to \Big[ \bigoplus_{\substack{\nu \prec \lambda \\ \res(\lambda/\nu) > 0}} \fr_{\set{\nu,\mu}} \oplus \fr_{\set{\nu^*,\mu}} \Big] \oplus \Big[ \bigoplus_{\substack{\nu \prec \lambda \\ \res(\lambda/\nu) = 0}} \fr_{\set{\nu,\mu}} \Big] \\
	&\mathrel{\phantom{=}} \oplus \Big[ \bigoplus_{\substack{\tau \prec \mu \\ \res(\mu/\tau) > 0}} \fr_{\set{\lambda,\tau}} \oplus \fr_{\set{\lambda,\tau^*}} \Big] \oplus \Big[ \bigoplus_{\substack{\tau \prec \mu \\ \res(\mu/\tau) = 0}} \fr_{\set{\lambda,\tau}} \Big],
	\end{align*}
which by the discussion of Section \ref{subsubsec:dual-pairs} factors through the map
	\begin{equation} \label{eq:r-F(n)-osp(lam,mu)-case-4-restricted}
	\rho: \fr' \to \Big[ \bigoplus_{\substack{\nu \prec \lambda \\ \res(\lambda/\nu) > 0}} \fr_{\set{\nu,\mu}} \Big] \oplus \Big[ \bigoplus_{\substack{\nu \prec \lambda \\ \res(\lambda/\nu) = 0}} \fr_{\set{\nu,\mu}} \Big] \oplus \Big[ \bigoplus_{\substack{\tau \prec \mu \\ \res(\mu/\tau) > 0}} \fr_{\set{\lambda,\tau}} \oplus \Big] \oplus \Big[ \bigoplus_{\substack{\tau \prec \mu \\ \res(\mu/\tau) = 0}} \fr_{\set{\lambda,\tau}} \Big].
	\end{equation}
The right side of \eqref{eq:r-F(n)-osp(lam,mu)-case-4-restricted} is a summand of the right side of \eqref{eq:dn-refined-factorization} for the value $n-1$. This implies by the induction hypothesis that \eqref{eq:r-F(n)-osp(lam,mu)-case-4-restricted} is a surjection. By induction, one has
	\begin{align*}
	\fr_{\set{\nu,\mu}} &\cong \begin{cases} 
	\fsl\{\nu,\mu\} & \text{if $\res(\lambda/\nu) > 0$,} \\
	\osp\{\nu,\mu\} & \text{if $\res(\lambda/\nu) = 0$,}
	\end{cases} & 
	\fr_{\set{\lambda,\tau}} &\cong \begin{cases} 
	\fsl\{\lambda,\tau \} & \text{if $\res(\mu/\tau) > 0$,} \\
	\osp\{\lambda,\tau \} & \text{if $\res(\mu/\tau) = 0$.} 
	\end{cases}
	\end{align*}

Let $\fh = \rho_{\set{\lambda,\mu}}(\fr') \cong \rho(\fr')$. For $\nu \prec \lambda$ with $\res(\lambda/\nu) > 0$, the modules $S^{\set{\nu,\mu}}$ and $S^{\set{\nu^*,\mu}}$ are dual but not isomorphic as $\fh$-modules, because $\dim(S^{\set{\nu,\mu}}) = \binom{n-1}{\abs{\mu}} \cdot f^\nu \cdot f^\mu \geq n-1 \geq 5$, and the natural module for $\fsl(k)$ is self-dual only if $k=2$. Similarly, if $\res(\mu/\tau) > 0$, then $S^{\set{\lambda,\tau}}$ and $S^{\set{\lambda,\tau^*}}$ are dual but not isomorphic as $\fh$-modules. Now $\fh$ is a semisimple subalgebra of $\fd_{\set{\lambda,\mu}}$ such that the restriction of $S^{\set{\lambda,\mu}}$ to each simple ideal of $\fh$ admits a simple factor of multi\-plicity one, and
	\begin{align*}
	2 \cdot \rk(\fh) &= \Big[ \sum_{\substack{\nu \prec \lambda \\ \res(\lambda/\nu) > 0}} 2 \cdot (\dim(S^{\set{\nu,\mu}})-1) \Big] + \Big[ \sum_{\substack{\nu \prec \lambda \\ \res(\lambda/\nu) = 0}} \dim(S^{\set{\nu,\mu}}) \Big] \\
	&\mathrel{\phantom{=}} + \Big[ \sum_{\substack{\tau \prec \mu \\ \res(\mu/\tau) > 0}} 2 \cdot (\dim(S^{\set{\lambda,\tau}})-1) \Big] + \Big[ \sum_{\substack{\tau \prec \mu \\ \res(\mu/\tau) = 0}} \dim(S^{\set{\lambda,\tau}}) \Big] \\
	&\geq \dim(S^{\set{\lambda,\mu}}) - r_\lambda - r_\mu \\
	&> \dim(S^{\set{\lambda,\mu}}) - n.
	\end{align*}
Since $\fd_{\set{\lambda,\mu}} \subseteq \osp\{\lambda,\mu\}$, then $\rk(\fd_{\set{\lambda,\mu}}) \leq \frac{1}{2} \cdot \dim(S^{\set{\lambda,\mu}})$, and hence
	\[ \textstyle
	2 \cdot \rk(\fh) - \rk(\fd_{\set{\lambda,\mu}}) > \frac{1}{2} \cdot \dim(S^{\set{\lambda,\mu}}) - n.
	\]
Since $\lambda$ and $\mu$ are symmetric and not equal to $[1]$, then $f^\lambda \cdot f^\mu \geq 2 \cdot 2$, and hence
	\[ \textstyle
	\dim(S^{\set{\lambda,\mu}}) = \binom{n}{\abs{\mu}} \cdot f^\lambda \cdot f^\mu \geq n \cdot f^\lambda \cdot f^\mu \geq 4n.
	\]
Then $2 \cdot \rk(\fh) - \rk(\fd_{\set{\lambda,\mu}}) \geq n > 0$, so $\fd_{\set{\lambda,\mu}}$ is simple by \cite[Lemme 15]{Marin:2007}. Finally, since $\rk(\fh) \leq \rk(\fd_{\set{\lambda,\mu}})$, one sees that $2 \cdot \rk(\fd_{\set{\lambda,\mu}}) \leq \dim(S^{\set{\lambda,\mu}}) < 4 \cdot \rk(\fd_{\set{\lambda,\mu}})$, and then deduces using \cite[Lemme 14]{Marin:2007} that $\fd_{\set{\lambda,\mu}} = \osp\{\lambda,\mu\}$.

\subsection{Proof of Theorem \ref{thm:main-theorem}(\ref{item:Fn-osp(lam,pm)-d})} \label{subsec:item:Fn-osp(lam,pm)-d}

\begin{lemma} \label{lemma:proof-of-item-Fn-osp(lam,pm)-d}
Let $n \geq 6$, and suppose Theorem \ref{thm:main-theorem} is true for the value $n-1$. Then part \eqref{item:Fn-osp(lam,pm)-d} of Theorem \ref{thm:main-theorem} is true for the value $n$, i.e.,
	\begin{center}
	If $\set{\lambda,\lambda} \in F\{n\}$ and $n/2$ is even, then $\fd_{\set{\lambda,\pm}} = \osp\{\lambda,\pm \}$ as described in Lemma \ref{lemma:d-osp(lam,pm)-inclusion}.
	\end{center}
\end{lemma}

\begin{proof}
Suppose $\set{\lambda,\lambda} \in F\{n\}$. Let $\fr = \fd_{n-1}$. By \eqref{eq:S(lam-pm)-restriction}, one gets the $\fr$-module identification
	\[
	S^{\set{\lambda,\pm}} = \Big[ \bigoplus_{\substack{\nu \prec \lambda \\ \res(\lambda/\nu) > 0}} S^{\set{\nu,\lambda}} \oplus S^{\set{\nu^*,\lambda}} \Big] \oplus \Big[ \bigoplus_{\substack{\nu \prec \lambda \\ \res(\lambda/\nu) = 0}} S^{\set{\nu,\lambda}} \Big].
	\]
For all $\nu \prec \lambda$, $\set{\nu,\lambda}$ is proper and is not an arm and leg by Lemma \ref{lemma:(lam-mu)-E(n)-preceding}. Thus, if $\nu \prec \lambda$ and $\res(\lambda/\nu) = 0$, then $\set{\nu,\lambda} \in F\{n-1\}$, and hence $\fr_{\set{\nu,\lambda}} = \osp\{\nu,\lambda\}$ by induction. On the other hand, suppose $\nu \prec \lambda$ and $\res(\lambda/\nu) > 0$. Then $\nu \neq \nu^*$, implying that $\set{\nu,\lambda} \in E\{n-1\}$, and hence $\fr_{\set{\nu,\lambda}} \cong \fsl\{\nu,\lambda\}$. The modules $S^{\set{\nu,\lambda}}$ and $S^{\set{\nu^*,\lambda}}$ are dual as $\fr$-modules, but Lemma \ref{lemma:Dn-subminimal-dimensions} implies that $\dim(S^{\set{\nu,\lambda}}) = \dim(S^{\set{\nu^*,\lambda}}) \geq (n-1)-1 \geq 4$, and hence $S^{\set{\nu,\lambda}} \not\cong S^{\set{\nu^*,\lambda}}$ as $\fr$-modules because the natural representation of $\fsl(k)$ is self-dual only if $k=2$.

Now the module structure map $\rho_{\set{\lambda,\pm}}: \fd_n' \to \fd_{\set{\lambda,\pm}} \subseteq \osp\{\lambda,\pm\} \subseteq \End(S^{\set{\lambda,\pm}})$ restricts to a Lie algebra homomorphism
	\[
	\fr' \to \Big[ \bigoplus_{\substack{\nu \prec \lambda \\ \res(\lambda/\nu) > 0}} \fr_{\set{\nu,\lambda}} \oplus \fr_{\set{\nu^*,\lambda}} \Big] \oplus \Big[ \bigoplus_{\substack{\nu \prec \lambda \\ \res(\lambda/\nu) = 0}} \fr_{\set{\nu,\lambda}} \Big],
	\]
which in turn factors through the map
	\begin{equation} \label{eq:r-F(n)-sl(lam,pm)-restricted}
	\rho: \fr' \to \Big[ \bigoplus_{\substack{\nu \prec \lambda \\ \res(\lambda/\nu) > 0}} \fr_{\set{\nu,\lambda}} \Big] \oplus \Big[ \bigoplus_{\substack{\nu \prec \lambda \\ \res(\lambda/\nu) = 0}} \fr_{\set{\nu,\lambda}} \Big].
	\end{equation}
The right side of \eqref{eq:r-F(n)-sl(lam,pm)-restricted} is a summand of the right side of \eqref{eq:dn-refined-factorization} for the value $n-1$. This implies by the induction hypothesis that \eqref{eq:r-F(n)-sl(lam,pm)-restricted} is a surjection.

Let $\fh = \rho_{\set{\lambda,\pm}}(\fr')$. Then $\fh \cong \rho(\fr')$ is a semisimple subalgebra of $\fd_{\set{\lambda,\pm}}$ such that the restriction of $S^{\set{\lambda,\pm}}$ to each simple ideal $\fr_{\set{\nu,\lambda}}$ of $\fh$ admits a simple factor $S^{\set{\nu,\lambda}}$ of multiplicity one, and
	\begin{equation} \label{eq:d-(lam,lam)-F(n)-n/2-even}
	\begin{split}
	2 \cdot \rk(\fh) &= \Big[ \sum_{\substack{\nu \prec \lambda \\ \res(\lambda/\nu) > 0}} 2 \cdot (\dim(S^{\set{\nu,\lambda}})-1) \Big] + \Big[ \sum_{\substack{\nu \prec \lambda \\ \res(\lambda/\nu) = 0}} \dim(S^{\set{\nu,\lambda}}) \Big] \\
	&= \Big[ \sum_{\substack{\nu \prec \lambda \\ \res(\lambda/\nu) \neq 0}}\dim(S^{\set{\nu,\lambda}})-1 \Big] + \Big[ \sum_{\substack{\nu \prec \lambda \\ \res(\lambda/\nu) = 0}} \dim(S^{\set{\nu,\lambda}}) \Big] \\
	&\geq \dim(S^{\set{\lambda,\pm}}) - r_\lambda \\
	&> \dim(S^{\set{\lambda,\pm}}) - \tfrac{n}{2}.
	\end{split}
	\end{equation}
Since $\fd_{\set{\lambda,\pm}} \subseteq \osp\{\lambda,\pm\}$, then $\rk(\fd_{\set{\lambda,\pm}}) \leq \frac{1}{2} \cdot \dim(S^{\set{\lambda,\pm}})$, while by reasoning like that in \eqref{eq:half-dim-S(lam,lam*)}, one has $\dim(S^{\set{\lambda,\pm}}) = \frac{1}{2} \cdot \dim(S^{\set{\lambda,\lambda}}) \geq 2 \cdot \binom{n}{n/2} \geq 2n$. Then
	\begin{align*}
	2 \cdot \rk(\fh) - \rk(\fd_{\set{\lambda,\pm}}) &> \textstyle \big[ \dim(S^{\set{\lambda,\pm}}) - \tfrac{n}{2} \big] - \rk(\fd_{\set{\lambda,\pm}})\\
	&\geq \textstyle \frac{1}{2} \cdot \dim(S^{\set{\lambda,\pm}}) - \tfrac{n}{2} \\
	&\geq n - \tfrac{n}{2} > 0.
	\end{align*}
This implies by \cite[Lemme 15]{Marin:2007} that $\fd_{\set{\lambda,\pm}}$ is a simple Lie algebra. Finally, since $\rk(\fh) \leq \rk(\fd_{\set{\lambda,\pm}})$, one sees that $2 \cdot \rk(\fd_{\set{\lambda,\pm}}) \leq \dim(S^{\set{\lambda,\pm}}) < 4 \cdot \rk(\fd_{\set{\lambda,\pm}})$, and then deduces using \cite[Lemme 14]{Marin:2007} that $\fd_{\set{\lambda,\pm}} = \osp\{\lambda,\pm \}$.
\end{proof}

\subsection{Proof of Theorem \ref{thm:main-theorem}(\ref{item:Fn-sl(lam,+)})}

The proof of the next lemma requires extra work compared to the other parts of Theorem \ref{thm:main-theorem}, on account of the fact that the restriction of $S^{\set{\lambda,\pm}}$ to $\fd_{n-1}'$ is qualitatively the same for $n/2$ odd as for $n/2$ even, and hence a rank analysis like that given in Section \ref{subsec:item:Fn-osp(lam,pm)-d} is not sufficient to deduce that $\fd_{\set{\lambda,\pm}} = \fsl\{\lambda,\pm\}$.

\begin{lemma} \label{lemma:proof-of-item-Fn-sl(lam,+)}
Let $n \geq 6$, and suppose Theorem \ref{thm:main-theorem} is true for the value $n-1$. Then part \eqref{item:Fn-sl(lam,+)} of Theorem \ref{thm:main-theorem} is true for the value $n$, i.e.,
	\begin{center}
	If $\set{\lambda,\lambda} \in F\{n\}$ and $n/2$ is odd, then $\fd_{\set{\lambda,\pm}} = \fsl\{\lambda,\pm\}$.
	\end{center}
\end{lemma}

\begin{proof}
The proof begins with a word-for-word repetition of the proof of Lemma \ref{lemma:proof-of-item-Fn-osp(lam,pm)-d} through the conclusion in \eqref{eq:d-(lam,lam)-F(n)-n/2-even} that $2 \cdot \rk(\fh) \geq \dim(S^{\set{\lambda,\pm}}) - r_\lambda > \dim(S^{\set{\lambda,\pm}}) - \frac{n}{2}$. It then also follows, by the same reasoning as in the proof of Lemma \ref{lemma:proof-of-item-Fn-osp(lam,pm)-d}, that $\dim(S^{\set{\lambda,\pm}}) < 4 \cdot \rk(\fd_{\set{\lambda,\pm}})$.

Let $\g = \fd_{\set{\lambda,\pm}}$ and $V = S^{\set{\lambda,\pm}}$. First we will show that $\g$ is simple. Since $\g$ is semisimple, then $\g = \g_1 \oplus \cdots \oplus \g_t$ for some simple ideals $\g_1,\ldots,\g_t$. We want to show that $t = 1$. Since $\g$ acts irreducibly on $V$, it follows for each $1 \leq i \leq t$ that there exists a simple $\g_i$-module $V_i$ such that, via the identification $\g = \g_1 \oplus \cdots \oplus \g_t$, $V$ identifies with the outer (or external) tensor product of modules $V_1 \boxtimes \cdots \boxtimes V_t$. Since $\g \subseteq \End(S^{\set{\lambda,\pm}})$ acts faithfully on $V = S^{\set{\lambda,\pm}}$, it follows for each $i$ that $V_i$ cannot be the trivial $\g_i$ module (else $\g_i$ acts trivially on $V$), and hence $\dim(V_i) \geq 2$.

Let $\fh = \bigoplus_j \fh_j$ be the decomposition of $\fh$ into simple ideals. Then $V_i$ becomes an $\fh_j$-module via the canonical maps $\fh_j \hookrightarrow \g \twoheadrightarrow \g_i$. If $V_i$ were trivial as an $\fh_j$-module, it would follow that each simple $\fh_j$-module factor in $V$ would occur with multiplicity at least $\dim(V_i)$, contradicting the observation that $\fh_j$ has a simple factor in $V$ of multiplicity one. Then $V_i$ must contain a nontrivial simple factor for $\fh_j$. We know for each $j$ that either $\fh_j \cong \fsl\{\nu,\lambda\}$ or $\fh_j \cong \osp\{\nu,\lambda\}$ for some $\nu \prec \lambda$. In either case, each nontrivial simple $\fh_j$-module has dimension greater than or equal to
	\[ \textstyle
	\dim(S^{\set{\nu,\lambda}}) = \binom{n-1}{n/2} \cdot \dim(S^\nu) \cdot \dim(S^\lambda) \geq \binom{n-1}{n/2} \cdot \dim(S^\lambda).
	\]
Then $\dim(V_i) \geq \binom{n-1}{n/2} \cdot \dim(S^\lambda)$ for each $i$. Now suppose $t \geq 2$. Then
	\[ \textstyle
	\dim(V) \geq \dim(V_1) \cdot \dim(V_2) \geq \big[ \binom{n-1}{n/2} \cdot \dim(S^\lambda) \big]^2 \geq \binom{n-1}{n/2} \cdot \dim(S^{\set{\lambda,\pm}}) > \dim(V).
	\]
This is a contradiction, so we must have $t=1$, and hence $\g$ is simple.

Now if $2 \cdot \rk(\g) > \dim(V)$, then $\g \cong \fsl(V)$ by \cite[Lemme 13]{Marin:2007}, as desired. Otherwise one has
	\[
	2 \cdot \rk(\g) \leq \dim(V) < 4 \cdot \rk(\g).
	\]
One has $\dim(V) \geq 2 \cdot \binom{n}{n/2}$ as in \eqref{eq:half-dim-S(lam,lam*)}, and since $n \geq 6$, this implies that $\dim(V) \geq 40$. Then by \cite[Lemma 14]{Marin:2007}, exactly one of the following must be true:
	\begin{enumerate}
	\item $\g$ is of type $C$ and $V$ is isomorphic to the natural module for $\g$, or
	\item $\g$ is of type $D$ and $V$ is isomorphic to the natural module for $\g$.
	\end{enumerate}
This implies that there exists a nondegenerate $\g$-invariant bilinear form on $V$. 
We will show that this leads to a contradiction. 
We will assume that $V = S^{\set{\lambda,+}}$; the case $V = S^{\set{\lambda,-}}$ is similar.

Suppose there exists a nondegenerate $\g$-invariant symmetric bilinear form $(-,-)$ on $V$. Thus for all $u,v \in V$ and $x \in \g$, one has $(u,v) = (v,u)$ and $(x.u,v) = -(u,x.v)$. By abuse of notation, we will conflate an element $x \in \fd_n'$ with its image in $\g = \fd_{\set{\lambda,+}}$. Let $p: \fd_n \to Z(\fd_n)$ be the projection map of Lemma \ref{lemma:Z(g)}. One has $Z(\fd_n) = \Span\{\calX_n\}$ by Lemma \ref{lemma:bn-dn-center}, and since $\lambda = \lambda^*$, one sees that $\calX_n$ acts as zero on $S^{\set{\lambda,\lambda}}$, and hence as zero on $V$. Then for arbitrary $z \in \fd_n$, one has
	\[
	\big( z.u,v \big) = \big( [z-p(z)].u,v \big) = -\big( u,[z-p(z)].v \big) = 
	-\big( u,z.v \big)
	\]
and hence $(-,-)$ is $\fd_n$-invariant. This implies for $s$ in the set $\set{s_1,\ldots,s_{n-1},\wt{s}_n}$ of generators for $\fd_n$ that $(s.u,s.v) = -(u,v)$. Since these elements generate $\calD_n$, it follows for $g \in \calD_n$ that
	\begin{equation} \label{eq:Dn-sign-invariance}
	(g.u,g.v) = \ve''(g) \cdot (u,v).
	\end{equation}

For $T \in \T(\lambda,\lambda)_{+1}$, set $u_T = c_T + c_{T^\natural}$. Then $\set{u_T : T \in \T(\lambda,\lambda)_{+1}}$ is a weight basis for $V$ as in Lemma \ref{lemma:type-D-simple-weight-spaces}. In particular, the YJM element $X_i$ acts on $u_T$ as scalar multiplication by $2 \cdot \res_T(i)$. Then for $S,T \in \T(\lambda,\lambda)_{+1}$, the $\fd_n$-invariance of the form implies that either $(u_S, u_T) = 0$, or else $\res_S(i) = -\res_T(i)$ for all $1 \leq i \leq n$. Similarly, given indices $1 \leq i,j \leq n$, the condition \eqref{eq:Dn-sign-invariance} applied to the group element $t_it_j \in \calD_n$ implies that if $(u_S,u_T) \neq 0$, then $\rho_S(i) = \rho_S(j)$ if and only if $\rho_T(i) = \rho_T(j)$, i.e., the integers $i$ and $j$ are located in the same tableau of $S$ if and only if they are located in the same tableau of $T$. Now given $S \in \T(\lambda,\lambda)_{+1}$, there exists $T \in \T(\lambda,\lambda)_{+1}$ such that $(u_S,u_T) \neq 0$, by non-degeneracy of the form. The previous observations then imply that $S = T^*$, and one concludes in general for $S,T \in \T(\lambda,\lambda)_{+1}$ that $(u_S,u_T) \neq 0$ if and only if $S = T^*$.

Let $R = R^{(\lambda,\lambda)} \in \T(\lambda,\lambda)_{+1}$ be the row major bitableau defined prior to Lemma \ref{lemma:epsilon''-intertwinor}. Rescaling the form $(-,-)$ if necessary, we may assume that $(u_{R^*},u_R) = 1$. Let $T \in \T(\lambda,\lambda)_{+1}$. The proof of \cite[Lemma 6.2]{Mishra:2016} shows that $T$ can be obtained from $R$ by applying a sequence of $\ell(T)$ admissible adjacent transpositions, none of which move the integer $1$, so at each intermediate step one still has an element of $\T(\lambda,\lambda)_{+1}$. Then using \eqref{eq:Dn-sign-invariance}, one can show for arbitrary $T \in \T(\lambda,\lambda)_{+1}$ that
	\[
	(u_{T^*},u_T) = (-1)^{\ell(T)}.
	\]

Now let $T \in \T(\lambda,\lambda)_{+1}$ such that $\rho_T(2) = -1$, i.e., such that $1$ and $2$ are not in the same tableau in $T$. Let $s_1 = (1,2)$, and let $S = s_1 \cdot T$. Then $S^{\natural} \in \T(\lambda,\lambda)_{+1}$, and applying Lemma \ref{lemma:S(lam-pm)-si-action} one can check that $s_1 \cdot u_T = u_{S^\natural}$ and $s_1 \cdot u_{T^*} = u_{(S^\natural)^*}$. Then, using the fact that $n/2$ is odd,
	\[
	(s_1 \cdot u_{T^*}, s_1 \cdot u_T) = (u_{(S^\natural)^*},u_{S^\natural}) = (-1)^{\ell(S^\natural)} = (-1)^{n/2} (-1)^{\ell(S)} = -(-1)^{\ell(s_1 \cdot T)} = (-1)^{\ell(T)}.
	\]
On the other hand, \eqref{eq:Dn-sign-invariance} implies that
	\[
	(s_1 \cdot u_{T^*}, s_1 \cdot u_T) = \ve''(s_1) \cdot (u_{T^*},u_T) = -(-1)^{\ell(T)},
	\]
a contradiction. Thus, it must be that $2 \cdot \rk(\g) > \dim(V)$, and hence $\g \cong \fsl(V) = \fsl(S^{\set{\lambda,+}})$.
\end{proof}

\subsection{Proof of Theorem \ref{thm:main-theorem}(\ref{item:En-sl(lam,mu)-d})}

\begin{lemma} \label{lemma:En-sl(lam,mu)-d-preceding}
Let $n \geq 6$, let $\set{\lambda,\mu} \in E\{n\}$, so $\set{\lambda,\mu} \neq \set{\lambda^*,\mu^*}$, and suppose $\lambda \neq \mu$.
	\begin{enumerate}
	\item \label{item:En-sl(lam,mu)-d-preceding-improper} Suppose $\abs{\lambda} \geq \abs{\mu}$. Then there exists an improper element $\set{\nu,\tau} \prec \set{\lambda,\mu}$ if and only if $\mu = [1]$. If $\mu = [1]$, then $\set{\lambda,\emptyset}$ is the unique improper element such that $\set{\lambda,\emptyset} \prec \set{\lambda,\mu}$.

	\item \label{item:En-sl(lam,mu)-d-preceding-AL} Suppose $\set{\nu,\tau}$ is an arm and leg and $\set{\nu,\tau} \prec \set{\lambda,\mu}$. Then
		\[
		\set{\lambda,\mu} = \set{[n-1-i,1],[1^i]} \quad \text{or} \quad \set{\lambda,\mu} = \set{[n-1-i],[2,1^{i-1}]}
		\]
	for some $1 \leq i \leq n-2$, and $\set{\nu,\tau}$ is the unique arm and leg such that $\set{\nu,\tau} \prec \set{\lambda,\mu}$.

	\item \label{item:En-sl(lam,mu)-d-preceding-splits} Suppose $\set{\nu,\nu} \prec \set{\lambda,\mu}$. Then either $\nu = \mu$ and $\mu \prec \lambda$, or $\nu = \lambda$ and $\lambda \prec \mu$.

	\item \label{item:En-sl(lam,mu)-d-preceding-dual-pair} Suppose $\set{\nu,\tau} \prec \set{\lambda,\mu}$ and $\set{\nu^*,\tau^*} \prec \set{\lambda,\mu}$. Then either:
		\begin{enumerate}
		\item \label{item:d-dual-pair-mu*-mu} $\set{\nu,\tau} = \set{\mu^*,\mu}$ and $\mu^* \prec \lambda$, or

		\item \label{item:d-dual-pair-lam-lam*} $\set{\nu,\tau} = \set{\lambda,\lambda^*}$ and $\lambda^* \prec \mu$, or

		\item \label{item:d-dual-pair-nu-mu-symmetric} $\set{\nu,\tau} = \set{\nu,\mu}$ with $\nu^* = \nu \prec \lambda$ and $\mu = \mu^*$, or

		\item \label{item:d-dual-pair-lam-tau-symmetric} $\set{\nu,\tau} = \set{\lambda,\tau}$ with $\tau^* = \tau \prec \mu$ and $\lambda = \lambda^*$. 
		\end{enumerate}
		In particular, $\set{\nu,\tau} = \set{\nu^*,\tau^*}$.
	\end{enumerate}
\end{lemma}

\begin{proof}
Parts \eqref{item:En-sl(lam,mu)-d-preceding-improper} and \eqref{item:En-sl(lam,mu)-d-preceding-AL} were stated in Lemma \ref{lemma:(lam-mu)-E(n)-preceding}, and part \eqref{item:En-sl(lam,mu)-d-preceding-splits} is immediate. For part \eqref{item:En-sl(lam,mu)-d-preceding-dual-pair}, let $\nu \prec \lambda$, so that $\set{\nu,\mu} \prec \set{\lambda,\mu}$, and suppose that $\set{\nu^*,\mu^*} \prec \set{\lambda,\mu}$. We will show that either \eqref{item:d-dual-pair-mu*-mu} or \eqref{item:d-dual-pair-nu-mu-symmetric} must be true. Since $\set{\nu^*,\mu^*} \prec \set{\lambda,\mu}$, then either
	\begin{itemize}
	\item $\set{\nu^*,\mu^*} = \set{\lambda,\tau}$ for some $\tau \prec \mu$, or
	\item $\set{\nu^*,\mu^*} = \set{\eta,\mu}$ for some $\eta \prec \lambda$.
	\end{itemize}
If $\set{\nu^*,\mu^*} = \set{\lambda,\tau}$ for some $\tau \prec \mu$, then either $\nu^* = \lambda$ and $\mu^* = \tau$, a contradiction because $\abs{\nu^*} = \abs{\nu} < \abs{\lambda}$, or $\nu^* = \tau$ and $\mu^* = \lambda$, a contradiction because $\set{\lambda,\mu} \neq \set{\lambda^*,\mu^*}$. Then $\set{\nu^*,\mu^*} = \set{\eta,\mu}$ for some $\eta \prec \lambda$. Now either
	\begin{itemize}
	\item $\nu^* = \mu$ and $\mu^* = \eta$ (case \eqref{item:d-dual-pair-mu*-mu} of the lemma), or 
	\item $\nu^* = \eta$ and $\mu^* = \mu$.
	\end{itemize}
Suppose $\nu^* = \eta$ and $\mu = \mu^*$. If $\nu \neq \nu^*$, then the Young diagram of $\lambda$ must have distinct removable boxes that can be removed to produce $\nu$ and $\nu^*$ respectively. Then the union of the Young diagrams of $\nu$ and $\nu^*$ must be the entirety of the Young diagram of $\lambda$, and hence $\lambda = \lambda^*$. But then $\lambda = \lambda^*$ and $\mu = \mu^*$, contradicting the assumption that $\set{\lambda,\mu} \neq \set{\lambda^*,\mu^*}$. Thus it must be the case that $\nu = \nu^*$, and hence $\set{\nu,\mu}$ matches case \eqref{item:d-dual-pair-nu-mu-symmetric} of the lemma.

Similar reasoning to that given above yields that either \eqref{item:d-dual-pair-lam-lam*} or \eqref{item:d-dual-pair-lam-tau-symmetric} is true if one starts instead with some $\tau \prec \mu$ such that $\set{\lambda,\tau} \prec \set{\lambda,\mu}$ and $\set{\lambda^*,\tau^*} \prec \set{\lambda,\mu}$. 
\end{proof}

\begin{lemma} \label{lemma:proof-of-item-En-sl(lam,mu)-d}
Let $n \geq 6$, and suppose Theorem \ref{thm:main-theorem} is true for the value $n-1$. Then part \eqref{item:En-sl(lam,mu)-d} of Theorem \ref{thm:main-theorem} is true for the value $n$, i.e.,
	\begin{center}
	If $\set{\lambda,\mu} \in E\{n\}$ and $\lambda \neq \mu$, then $\fd_{\set{\lambda,\mu}} = \fsl\{\lambda,\mu\}$.
	\end{center}
\end{lemma}

The proof of Lemma \ref{lemma:proof-of-item-En-sl(lam,mu)-d} will occupy the rest of this section. Let $\set{\lambda,\mu} \in E\{n\}$ with $\lambda \neq \mu$, and let $\fr = \fd_{n-1}$. Then the $\fd_n'$-module structure map $\rho_{\set{\lambda,\mu}}: \fd_n' \twoheadrightarrow \fd_{\set{\lambda,\mu}} \subseteq \fsl\{\lambda,\mu\}$ restricts to a Lie algebra homomorphism
	\[
	\rho: \fr' \to \Big[ \bigoplus_{\substack{\set{\nu,\tau} \prec \set{\lambda,\mu} \\ \nu \neq \tau}} \fr_{\set{\nu,\tau}} \Big] \oplus \Big[ \bigoplus_{\set{\nu,\nu} \prec \set{\lambda,\mu}} \fr_{\set{\nu,+}} \oplus \fr_{\set{\nu,-}} \Big],
	\]
which in turn factors through the map
	\begin{equation} \label{eq:r-E(n)-sl(lam,mu)-restricted}
	\rho: \fr' \to \Big[ \bigoplus_{\substack{\set{\nu,\tau} \prec \set{\lambda,\mu} \\ \nu \neq \tau}} \fr_{\set{\nu,\tau}} \Big] \oplus \Big[ \bigoplus_{\set{\nu,\nu} \prec \set{\lambda,\mu}} \fr_{\set{\nu,+}} \oplus \fr_{\set{\nu,-}}^\star \Big].
	\end{equation}
where the $\star$ indicates that the term $\fr_{\set{\nu,-}}$ is omitted if $\nu = \nu^*$ and $\abs{\nu}$ is odd. It follows from Lemma \ref{lemma:En-sl(lam,mu)-d-preceding} that the right side of \eqref{eq:r-E(n)-sl(lam,mu)-restricted} is a summand of the right side of \eqref{eq:bn-refined-factorization} for the value $n-1$, except perhaps that a summand of the form $\fr_{\beta_{n-1,1}}$ may be replaced by an isomorphic summand of the form $\fr_{\beta_{n-1,d}}$ for some $1 \leq d \leq n-2$, and a summand of the form $\fs_{\alpha_{n-1,1}}$ may be replaced by an isomorphic summand $\fs_{\alpha_{n-1,d}}$ for some $1 \leq d \leq n-3$. This implies by induction that \eqref{eq:r-E(n)-sl(lam,mu)-restricted} is a surjection, and implies for each $\set{\nu,\tau} \prec \set{\lambda,\mu}$ with $\nu \neq \tau$ (resp.\ $\nu = \tau$) that $\fr_{\set{\nu,\tau}}$ is a simple Lie algebra acting irreducibly on $S^{\set{\nu,\tau}}$ (resp.\ that $\fr_{\set{\nu,\pm}}$ is a simple acting irreducibly on $S^{\set{\nu,\pm}}$).

Let $\fh = \rho_{\set{\lambda,\mu}}(\fr')$. If $\set{\nu,\nu} \prec \set{\lambda,\mu}$ with $\nu = \nu^*$ and $\abs{\nu}$ odd, then $\fr_{\set{\nu,\pm}} \cong \fsl\{\nu,\pm\}$ by induction, and it follows that $S^{\set{\nu,+}}$ and $S^{\set{\nu,-}}$ are dual but not isomorphic as $\fh$-modules because $\dim(S^{\set{\nu,\pm}}) = \frac{1}{2} \cdot \binom{n-1}{\abs{\nu}} \cdot f^\nu \cdot f^\nu \geq \frac{1}{2} \cdot (n-1) > 2$, and the natural module for $\fsl(k)$ is self-dual only if $k=2$. Now $\fh \cong \rho(\fr')$ is a semisimple subalgebra of $\fd_{\set{\lambda,\mu}}$ such that the restriction of $S^{\set{\lambda,\mu}}$ to each simple ideal of $\fh$ admits a simple factor of multiplicity one. Since $\fh \subseteq \fd_{\set{\lambda,\mu}} \subseteq \fsl\{\lambda,\mu\}$, then $\rk(\fh) \leq \rk(\fd_{\set{\lambda,\mu}}) \leq \dim(S^{\set{\lambda,\mu}})$. Our goal now is to show that $2 \cdot \rk(\fh) > \dim(S^{\set{\lambda,\mu}})$. This will imply first by \cite[Lemme 15]{Marin:2007} that $\fd_{\set{\lambda,\mu}}$ is a simple Lie algebra, and then by \cite[Lemme 13]{Marin:2007} that $\fd_{\set{\lambda,\mu}} = \fsl\{ \lambda,\mu \}$. We break the remainder of the proof into the following cases:
	\begin{enumerate}
	\item \label{item:r(nu,tau)=sl(nu,tau)-d} $\fr_{\set{\nu,\tau}} \cong \fsl\{\nu,\tau\}$ for all $\set{\nu,\tau} \prec \set{\lambda,\mu}$.

	\item \label{item:mu=1-lam-hook-d} $\mu = [1]$ and $\lambda$ is a hook partition.
	
	\item \label{item:mu=1-lam-non-hook-d} $\mu = [1]$ and $\lambda$ is not a hook partition.
	
	\item \label{item:(nu,tau)-arm-and-leg-d} There exists an arm and leg $\set{\nu,\tau} \prec \set{\lambda,\mu}$, and all $\set{\eta,\sigma} \prec \set{\lambda,\mu}$ are proper.
	
	\item Each $\set{\eta,\sigma} \prec \set{\lambda,\mu}$ is proper and is not not an arm and leg, and either:
		\begin{enumerate}
		\item \label{item:(nu,nu)-d-lam-symmetric} $\mu \prec \lambda$ and $\lambda = \lambda^*$.

		\item \label{item:(nu,nu)-d-lam-not-symmetric-mu-symmetric} $\mu \prec \lambda$, $\lambda \neq \lambda^*$, and $\mu = \mu^*$.

		\item \label{item:(nu,nu)-d-lam-not-symmetric-mu-not-symmetric} $\mu \prec \lambda$, $\lambda \neq \lambda^*$, and $\mu \neq \mu^*$.

		\item \label{item:(nu,tau)-self-dual-mu*-prec-lam} $\mu \not\prec \lambda$ and $\mu^* \prec \lambda$ (and hence $\lambda \not\prec \mu$).
		
		\item \label{item:(nu,tau)-self-dual-nu-prec-lam-mu} $\mu \not \prec \lambda$, $\lambda \not \prec \mu$, there exists $\kappa \prec \lambda$ such that $\kappa = \kappa^*$, and $\mu = \mu^*$.
		\end{enumerate}
	\end{enumerate}
The cases where the roles of $\lambda$ and $\mu$ are swapped are similar.

\subsubsection{Completion of Lemma \ref{lemma:proof-of-item-En-sl(lam,mu)-d} in case (\ref{item:r(nu,tau)=sl(nu,tau)-d})}

Suppose $\fr_{\set{\nu,\tau}} \cong \fsl\{\nu,\tau\}$ for all $\set{\nu,\tau} \prec \set{\lambda,\mu}$. Then the conclusion $2 \cdot \rk(\fh) > \dim(S^{\set{\lambda,\mu}})$ follows as in Section \ref{subsubsec:r(nu,tau)=sl(nu,tau)}.

\subsubsection{Completion of Lemma \ref{lemma:proof-of-item-En-sl(lam,mu)-d} in case (\ref{item:mu=1-lam-hook-d})}

Suppose $\mu = [1]$ and $\lambda$ is a hook partition. Since $\set{\lambda,\mu} \neq \set{\lambda^*,\mu^*}$, then $\lambda \neq \lambda^*$. If $\set{\nu,\tau} \neq \set{\nu^*,\tau^*}$ for all $\set{\nu,\tau} \prec \set{\lambda,\mu}$, i.e., if $\nu \neq \nu^*$ for all $\nu \prec \lambda$, then the conclusion  $2 \cdot \rk(\fh) > \dim(S^{\set{\lambda,\mu}})$ follows as in Section \ref{subsubsec:(nu,tau)-improper-lambda-hook}. So suppose $\nu = \nu^*$ for some $\nu \prec \lambda$. Then replacing $\set{\lambda,\mu}$ with $\set{\lambda^*,\mu^*}$ if necessary, we may assume that $n = 2d+1$ for some $d \geq 2$, and $\lambda = \alpha_{n-1,d} = [n-1-d,1^d]$. Then as an $\fr$-module,
	\[
	S^{\set{\lambda,\mu}} 
	= S^{\{ \alpha_{n-2,d},[1] \}} 
	\oplus S^{ \{ \alpha_{n-2,d-1},[1] \}} 
	\oplus S^{\{ \alpha_{n-1,d},\emptyset \}},
	\]
and $\{ \alpha_{n-2,d-1},[1] \} = \{ \alpha_{n-2,d-1}^*,[1]^* \}$. The dimensions of the modules are
	\begin{align*}
	\dim(S^{\set{\lambda,\mu}}) &= \textstyle n \cdot \binom{n-2}{d}, &
	\dim(S^{\{ \alpha_{n-2,d},[1] \}}) &= \textstyle (n-1) \cdot \binom{n-3}{d}, \\
	\dim(S^{ \{ \alpha_{n-2,d-1},[1] \}}) &= \textstyle (n-1) \cdot \binom{n-3}{d-1}, &
	\dim(S^{\{ \alpha_{n-1,d},\emptyset \}}) &= \textstyle \binom{n-2}{d}.
	\end{align*}
By induction, one has
	\begin{align*}
	\fr_{\{ \alpha_{n-2,d},[1] \}} &\cong \fsl\{ \alpha_{n-2,d},[1]\}, &
	\fr_{\{ \alpha_{n-2,d-1},[1] \}} &\cong \osp\{ \alpha_{n-2,d-1},[1] \}, &
	\fr_{\{ \alpha_{n-1,d},\emptyset \}} &\cong \fsl(n-2).
	\end{align*}
Then
	\begin{align*}
	2 \cdot \rk(\fh) - \dim(S^{\set{\lambda,\mu}} ) &= \textstyle 2 \cdot \big[ (n-1) \cdot \binom{n-3}{d} - 1 \big] +  (n-1) \cdot \binom{n-3}{d-1} + 2 \cdot \big[n-3\big] - n \cdot \binom{n-2}{d} \\
	&= \textstyle (n-1) \cdot \binom{n-3}{d} - \binom{n-2}{d} + (2n-8) \\
	&= \textstyle (2d) \cdot \binom{2d-2}{d} - \binom{2d-1}{d} + (2n-8) \\
	&= \textstyle \big[ (2d)(d-1)-(2d-1) \big] \cdot \frac{(2d-2)!}{d!(d-1)!} + (2n-8),
	\end{align*}
and this quantity is strictly positive for $d \geq 2$.

\subsubsection{Completion of Lemma \ref{lemma:proof-of-item-En-sl(lam,mu)-d} in case (\ref{item:mu=1-lam-non-hook-d})}

Suppose $\mu = [1]$ and $\lambda$ is not a hook partition. This implies by Lemma \ref{lemma:En-sl(lam,mu)-d-preceding} that $\set{\nu,\tau}$ is not an arm and leg for all $\set{\nu,\tau} \prec \set{\lambda,\mu}$. Since $\set{\lambda,\mu} \neq \set{\lambda^*,\mu^*}$, then $\lambda \neq \lambda^*$. If $\set{\nu,\tau} \neq \set{\nu^*,\tau^*}$ for all $\set{\nu,\tau} \prec \set{\lambda,\mu}$, i.e., if $\nu \neq \nu^*$ for all $\nu \prec \lambda$, then $\fr_{\set{\nu,\tau}} \cong \fsl\{\nu,\tau\}$ for all $\set{\nu,\tau} \prec \set{\lambda,\mu}$, and the conclusion  $2 \cdot \rk(\fh) > \dim(S^{\set{\lambda,\mu}})$ is true by case \eqref{item:r(nu,tau)=sl(nu,tau)-d}. So suppose $\nu = \nu^*$ for some (unique) $\nu \prec \lambda$. Then as an $\fr$-module,
	\[
	S^{\set{\lambda,[1]}} = \Big[ \bigoplus_{\substack{\tau \prec \lambda \\ \tau \neq \tau^*}} S^{\set{\tau,[1]}} \Big] \oplus S^{\set{\nu,[1]}} \oplus S^{\set{\lambda,\emptyset}}.
	\]
By induction, $\fr_{\{\lambda,\emptyset\}} \cong \fsl(\lambda)$, $\fr_{\{\nu,[1]\}} \cong \osp\{\nu,[1]\}$, and $\fr_{\{ \tau,[1] \}} \cong \fsl\{\tau,[1]\}$ for $\tau \neq \tau^*$. Then
	\[
	2 \cdot \rk(\fh) = 2 \cdot \Big[ \sum_{\substack{\tau \prec \lambda \\ \tau \neq \tau^*}} \dim(S^{\set{\tau,[1]}})-1 \Big] + \dim(S^{\set{\nu,[1]}}) + 2 \cdot \big[ \dim(S^{\set{\lambda,\emptyset}})-1 \big],
	\]
and hence
	\[
	2 \cdot \rk(\fh) - \dim(S^{\set{\lambda,[1]}}) = \Big[ \sum_{\substack{\tau \prec \lambda \\ \tau \neq \tau^*}} \dim(S^{\set{\tau,[1]}})-2 \Big] + \big[ \dim(S^{\set{\lambda,\emptyset}})-2 \big].
	\]
If $\tau \prec \lambda$, then $\dim(S^{\set{\tau,[1]}}) = (n-1) \cdot f^\tau \geq 5$. And since $\set{\lambda,[1]}$ is not an arm and leg, then $\lambda$ is not either $[n-1]$ or $[1^{n-1}]$, and hence $\dim(S^{\set{\lambda,\emptyset}}) = f^\lambda \geq (n-1)-1 \geq 4$. This implies that $2 \cdot \rk(\fh) > \dim(S^{\set{\lambda,[1]}})$, as desired.

\subsubsection{Completion of Lemma \ref{lemma:proof-of-item-En-sl(lam,mu)-d} in case (\ref{item:(nu,tau)-arm-and-leg-d})}

Suppose there exists an arm and leg $\set{\nu,\tau} \prec \set{\lambda,\mu}$, and all $\set{\eta,\sigma} \prec \set{\lambda,\mu}$ are proper. Replacing $\set{\lambda,\mu}$ with $\set{\lambda^*,\mu^*}$ if necessary, we may assume that $\set{\lambda,\mu} = \set{ [n-1-i,1], [1^i]}$ for some $2 \leq i \leq n-2$. If $\set{\nu,\tau} \neq \set{\nu^*,\tau^*}$ for all $\set{\nu,\tau} \prec \set{\lambda,\mu}$, then $2 \cdot \rk(\fh) - \dim(S^{\set{\lambda,\mu}}) > 0$ as in Section \ref{subsubsec:(nu,tau)-arm-and-leg}, so suppose $\set{\nu,\tau} = \set{\nu^*,\tau^*}$ for some $\set{\nu,\tau} \prec \set{\lambda,\mu}$. Then $n=2i+1$, $\set{\lambda,\mu} = \set{[i,1],[1^i]}$, and as $\fr$-modules,
	\[
	S^{\{[i,1],[1^i] \}} = S^{\{ [i],[1^i] \}} \oplus S^{\{ [i-1,1],[1^i] \}} \oplus S^{\{ [i,1],[1^{i-1}] \}},
	\]
with $\{[i],[1^i]\} = \{[i]^*,[1^i]^*\}$. The dimensions of the modules are
	\begin{align*}
	\dim(S^{\{[i,1],[1^i] \}}) &= \textstyle \binom{2i+1}{i} \cdot i, &
	\dim(S^{\{[i],[1^i] \}}) &= \textstyle \binom{2i}{i}, \\
	\dim(S^{\{[i-1,1],[1^i] \}}) &= \textstyle \binom{2i}{i} \cdot (i-1), &
	\dim(S^{\{[i,1],[1^{i-1}] \}}) &= \textstyle \binom{2i}{i-1} \cdot i.
	\end{align*}
By induction,
	\begin{align*}
	\fr_{\{ [i],[1^i] \}} &\cong \osp\{ [i],[1^i] \}, &
	\fr_{\{ [i-1,1],[1^i] \}} &\cong \fsl\{ [i-1,1],[1^i] \}, &
	\fr_{\{ [i,1],[1^{i-1}] \}} &\cong \fsl\{ [i,1],[1^{i-1}] \}.
	\end{align*}
Then
	\begin{align*}
	2 \cdot \rk(\fh) - \dim(S^{\set{\lambda,\mu}}) &= \textstyle \binom{2i}{i} + 2 \cdot \big[ \binom{2i}{i} \cdot (i-1) - 1 \big] + 2 \cdot \big[ \binom{2i}{i-1} \cdot i - 1 \big] - \binom{2i+1}{i} \cdot i \\
	&= \textstyle i \cdot \binom{2i+1}{i} - \binom{2i}{i} - 4 \\
	&\geq \textstyle 2 \cdot \binom{2i+1}{i} - \binom{2i}{i} - 4 \\
	&= \textstyle \binom{2i+1}{i} + \binom{2i}{i-1} - 4 \geq (4i+1) -4.
	\end{align*}
Thus $2 \cdot \rk(\fh) > \dim(S^{\set{\lambda,\mu}})$, as desired.

\subsubsection{Completion of Lemma \ref{lemma:proof-of-item-En-sl(lam,mu)-d} in case (\ref{item:(nu,nu)-d-lam-symmetric})}

Suppose each $\set{\eta,\sigma} \prec \set{\lambda,\mu}$ is proper and is not an arm and leg, and suppose $\mu \prec \lambda$ and $\lambda = \lambda^*$. Then $\mu \neq \mu^*$, because $\set{\lambda,\mu} \neq \set{\lambda^*,\mu^*}$, and hence $\res(\lambda/\mu) \neq 0$. We may assume that $\res(\lambda/\mu) > 0$. Then as an $\fr$-module,
	\begin{align*}
	S^{\set{\lambda,\mu}} = \Big[ &\bigoplus_{\substack{\mu \neq \nu \prec \lambda \\ \res(\lambda/\nu) > 0}} S^{\set{\nu,\mu}} \oplus S^{\set{\nu^*,\mu}} \Big] 
	\oplus \Big[ S^{\set{\mu,+}} \oplus S^{\set{\mu,-}} \Big] 
	\oplus S^{\set{\mu^*,\mu}} \\
	&\oplus \Big[ \bigoplus_{\substack{\nu \prec \lambda \\ \res(\lambda/\nu) = 0}} S^{\set{\nu,\mu}} \Big]
	\oplus \Big[ \bigoplus_{\substack{\tau \prec \mu \\ \tau \neq \tau^*}} S^{\set{\lambda,\tau}} \Big] 
	\oplus \Big[ \bigoplus_{\substack{\tau \prec \mu \\ \tau = \tau^*}} S^{\set{\lambda,\tau}} \Big]
	\end{align*}
The only labels $\set{\eta,\sigma}$ in the preceding sum such that $\set{\eta,\sigma} = \set{\eta^*,\sigma^*}$ are $\set{\mu^*,\mu}$ and $\set{\lambda,\tau}$ for $\tau = \tau^*$ (there is at most one $\tau$ with this property). Then by induction,
	\begin{align*}
	\fr_{\set{\nu,\mu}} &\cong \fsl\{\nu,\mu\} & \text{if $\nu \prec \lambda$ and $\nu \neq \mu$,} \\
	\fr_{\set{\mu,\pm}} &\cong \fsl\{\mu,\pm\}, \\
	\fr_{\set{\mu^*,\mu}} &\cong \osp\{\mu^*,\mu\}, \\
	\fr_{\set{\nu,\mu}} &\cong \fsl\{\nu,\mu\} & \text{if $\nu \prec \lambda$ and $\res(\lambda/\nu) = 0$,} \\
	\fr_{\set{\lambda,\tau}} &\cong \fsl\{\lambda,\tau\} & \text{if $\tau \prec \mu$ and $\tau \neq \tau^*$,} \\
	\fr_{\set{\lambda,\tau}} &\cong \osp\{\lambda,\tau\} & \text{if $\tau \prec \mu$ and $\tau = \tau^*$.}
	\end{align*}
Then
	\begin{align*}
	\rk(\fh) = \Big[ &\sum_{\substack{\mu \neq \nu \prec \lambda \\ \res(\lambda/\nu) > 0}} \dim(S^{\set{\nu,\mu}}) + \dim(S^{\set{\nu^*,\mu}})- 2 \Big] 
	+ \Big[ \dim(S^{\set{\mu,\mu}}) -2 \Big] 
	+ \tfrac{1}{2} \cdot \dim(S^{\set{\mu^*,\mu}}) \\
	&+ \Big[ \sum_{\substack{\nu \prec \lambda \\ \res(\lambda/\nu) = 0}} \dim(S^{\set{\nu,\mu}})-1 \Big]
	+ \Big[ \sum_{\substack{\tau \prec \mu \\ \tau \neq \tau^*}} \dim(S^{\set{\lambda,\tau}})-1 \Big] 
	+ \Big[ \sum_{\substack{\tau \prec \mu \\ \tau = \tau^*}} \tfrac{1}{2} \cdot \dim(S^{\set{\lambda,\tau}}) \Big],
	\end{align*}
and hence
	\begin{align*}
	2 \cdot \rk(\fh) - \dim(S^{\set{\lambda,\mu}}) = \Big[ &\sum_{\substack{\mu \neq \nu \prec \lambda \\ \res(\lambda/\nu) > 0}} \dim(S^{\set{\nu,\mu}}) + \dim(S^{\set{\nu^*,\mu}})- 4 \Big] 
	+ \Big[ \dim(S^{\set{\mu,\mu}}) -4 \Big]  \\
	&+ \Big[ \sum_{\substack{\nu \prec \lambda \\ \res(\lambda/\nu) = 0}} \dim(S^{\set{\nu,\mu}}) - 2 \Big]
	+ \Big[ \sum_{\substack{\tau \prec \mu \\ \tau \neq \tau^*}} \dim(S^{\set{\lambda,\tau}}) - 2 \Big].
	\end{align*}
Since each $\set{\eta,\sigma} \prec \set{\lambda,\mu}$ is proper by assumption, then $\dim(S^{\set{\eta,\sigma}}) = \binom{n}{\abs{\eta}} \cdot f^\eta \cdot f^\sigma \geq n \geq 6$. This implies that $2 \cdot \rk(\fh) > \dim(S^{\set{\lambda,\mu}})$, as desired.

\subsubsection{Completion of Lemma \ref{lemma:proof-of-item-En-sl(lam,mu)-d} in case (\ref{item:(nu,nu)-d-lam-not-symmetric-mu-symmetric})}

Suppose each $\set{\eta,\sigma} \prec \set{\lambda,\mu}$ is proper and is not an arm and leg, $\mu \prec \lambda$, $\lambda \neq \lambda^*$, and $\mu = \mu^*$. As an $\fr$-module,
	\[
	S^{\set{\lambda,\mu}} = \Big[ \bigoplus_{\substack{\nu \prec \lambda \\ \nu \neq \mu}} S^{\set{\nu,\mu}} \Big] 
	\oplus \Big[ S^{\set{\mu,+}} \oplus S^{\set{\mu,-}} \Big] 
	\oplus \Big[ \bigoplus_{\tau \prec \mu} S^{\set{\lambda,\tau}} \Big]
	\]
In this case, if $\set{\nu,\tau} \prec \set{\lambda,\mu}$ and $\set{\nu,\tau} \neq \set{\mu,\mu}$, then $\set{\nu,\tau} \neq \set{\nu^*,\tau^*}$. Now
	\[
	\fh \cong \Big[ \bigoplus_{\substack{\nu \prec \lambda \\ \nu \neq \mu}} \fr_{\set{\nu,\mu}} \Big] 
	\oplus \Big[ \fr_{\set{\mu,+}} \oplus \fr_{\set{\mu,-}}^\star \Big] 
	\oplus \Big[ \bigoplus_{\tau \prec \mu} \fr_{\set{\lambda,\tau}} \Big]
	\]
where the $\star$ indicates that $\fr_{\set{\mu,-}}$ is omitted if $\abs{\mu}$ is odd. By induction,
	\begin{align*}
	\fr_{\set{\nu,\mu}} &\cong \fsl\{\nu,\mu\} & \text{if $\nu \prec \lambda$ and $\nu \neq \mu$,} \\
	\fr_{\set{\mu,\pm}} &\cong \osp\{\mu,\pm\}, & \text{if $\abs{\mu}$ is even,} \\
	\fr_{\set{\mu,\pm}} &\cong \fsl\{\mu,\pm\}, & \text{if $\abs{\mu}$ is odd,} \\
	\fr_{\set{\lambda,\tau}} &\cong \fsl\{\lambda,\tau\} & \text{if $\tau \prec \mu$.}
	\end{align*}
Then it follows that
	\[
	\rk(\fh) \geq \Big[ \sum_{\substack{\nu \prec \lambda \\ \nu \neq \mu}} \dim(S^{\set{\nu,\mu}})-1 \Big] 
	+ \Big[ \dim(S^{\set{\mu,+}})-1 \Big] 
	+ \Big[ \sum_{\tau \prec \mu} \dim(S^{\set{\lambda,\tau}})-1 \Big],
	\]
and hence
	\[
	2 \cdot \rk(\fh) - \dim(S^{\set{\lambda,\mu}}) \geq \Big[ \sum_{\substack{\nu \prec \lambda \\ \nu \neq \mu}} \dim(S^{\set{\nu,\mu}})-2 \Big] 
	-2 
	+ \Big[ \sum_{\tau \prec \mu} \dim(S^{\set{\lambda,\tau}})-2 \Big].
	\]
Since each $\set{\eta,\sigma} \prec \set{\lambda,\mu}$ is proper by assumption, then $\dim(S^{\set{\eta,\sigma}}) \geq 6$ as in case \eqref{item:(nu,nu)-d-lam-symmetric}, and there is at least one $\tau \prec \mu$. This implies that $2 \cdot \rk(\fh) > \dim(S^{\set{\lambda,\mu}})$, as desired.

\subsubsection{Completion of Lemma \ref{lemma:proof-of-item-En-sl(lam,mu)-d} in case (\ref{item:(nu,nu)-d-lam-not-symmetric-mu-not-symmetric})}

Suppose each $\set{\eta,\sigma} \prec \set{\lambda,\mu}$ is proper and is not an arm and leg, and suppose that $\mu \prec \lambda$ and $\lambda \neq \lambda^*$ and $\mu \neq \mu^*$. Since $\lambda \neq \lambda^*$, then $\mu^* \not\prec \lambda$; cf.\ the reasoning in the proof of Lemma \ref{lemma:En-sl(lam,mu)-d-preceding}. Then by Lemma \ref{lemma:En-sl(lam,mu)-d-preceding}, there is no $\set{\eta,\sigma} \prec \set{\lambda,\mu}$ such that $\set{\eta,\sigma} = \set{\eta^*,\sigma^*}$. Now the reasoning in this case is parallel to that in \eqref{item:(nu,nu)-d-lam-not-symmetric-mu-symmetric}, except that the term $\fr_{\set{\mu,-}}$ is retained and $\fr_{\set{\mu,\pm}} \cong \fsl\{\mu,\pm\}$ regardless of the parity of $\abs{\mu}$. This results in
	\[
	2 \cdot \rk(\fh) - \dim(S^{\set{\lambda,\mu}}) = \Big[ \sum_{\substack{\nu \prec \lambda \\ \nu \neq \mu}} \dim(S^{\set{\nu,\mu}})-2 \Big] 
	+ \Big[ \dim(S^{\set{\mu,\mu}})-4 \Big] 
	+ \Big[ \sum_{\tau \prec \mu} \dim(S^{\set{\lambda,\tau}})-2 \Big].
	\]
By the same reasoning as in case \eqref{item:(nu,nu)-d-lam-symmetric}, one sees that this sum is strictly positive.

\subsubsection{Completion of Lemma \ref{lemma:proof-of-item-En-sl(lam,mu)-d} in case (\ref{item:(nu,tau)-self-dual-mu*-prec-lam})}

Suppose $\mu \not\prec \lambda$ and $\mu^* \prec \lambda$. Then $\mu \neq \mu^*$, and also $\lambda \neq \lambda^*$ (because if $\lambda$ were symmetric, it would be the case that $\nu \prec \lambda$ if and only if $\nu^* \prec \lambda$). Now as an $\fr$-module,
	\[
	S^{\set{\lambda,\mu}} = \Big[ \bigoplus_{\substack{\nu \prec \lambda \\ \nu \neq \mu^*}} S^{\set{\nu,\mu}} \Big] 
	\oplus S^{\set{\mu^*,\mu}} 
	\oplus \Big[ \bigoplus_{\tau \prec \mu} S^{\set{\lambda,\tau}} \Big],
	\]
and the only label $\set{\eta,\sigma} \prec \set{\lambda,\mu}$ such that $\set{\eta,\sigma} = \set{\eta^*,\sigma^*}$ is $\set{\mu^*,\mu}$. Then reasoning as in the previous cases, one gets
	\[
	2 \cdot \rk(\fh) - \dim(S^{\set{\lambda,\mu}}) = \Big[ \sum_{\substack{\nu \prec \lambda \\ \nu \neq \mu^*}} \dim(S^{\set{\nu,\mu}})-2 \Big]
	+ \Big[ \sum_{\tau \prec \mu} \dim(S^{\set{\lambda,\tau}})-2 \Big],
	\]
and deduces that this sum is strictly positive.

\subsubsection{Completion of Lemma \ref{lemma:proof-of-item-En-sl(lam,mu)-d} in case (\ref{item:(nu,tau)-self-dual-nu-prec-lam-mu})}

Suppose $\mu = \mu^*$, there exists a (unique) $\kappa \prec \lambda$ such that $\kappa = \kappa^*$, and neither $\mu \prec \lambda$ nor $\lambda \prec \mu$. Since $\set{\lambda,\mu} \neq \set{\lambda^*,\mu^*}$, then $\lambda \neq \lambda^*$. And since $\lambda \not \prec \mu$ and $\mu = \mu^*$, then also $\lambda^* \not\prec \mu$. Now as an $\fr$-module,
	\[
	S^{\set{\lambda,\mu}} = \Big[ \bigoplus_{\substack{\nu \prec \lambda \\ \nu \neq \nu^*}} S^{\set{\nu,\mu}} \Big] 
	\oplus S^{\set{\kappa,\mu}}
	\oplus \Big[ \bigoplus_{\tau \prec \mu} S^{\set{\lambda,\tau}} \Big],
	\]
and the only label $\set{\eta,\sigma} \prec \set{\lambda,\mu}$ such that $\set{\eta,\sigma} = \set{\eta^*,\sigma^*}$ is $\set{\kappa,\mu}$. Now reasoning as in case \eqref{item:(nu,tau)-self-dual-mu*-prec-lam}, but replacing $\set{\mu^*,\mu}$ with $\set{\kappa,\mu}$, one deduces that $2 \cdot \rk(\fh) > \dim(S^{\set{\lambda,\mu}})$.

\subsection{Proof of Theorem \ref{thm:main-theorem}(\ref{item:En-sl(lam,pm)-d})} \label{subsec:En-sl(lam,pm)-d}

\begin{lemma} \label{lemma:(lam,pm)-En-preceding}
Let $n \geq 6$, let $\set{\lambda,\lambda} \in E\{n\}$, so $\lambda \neq \lambda^*$, and let $\nu \prec \lambda$, so $\set{\nu,\lambda} \prec \set{\lambda,\lambda}$. Then $\set{\nu,\lambda}$ is proper, is not an arm and leg, and $\set{\nu,\lambda} \neq \set{\nu^*,\lambda^*}$. Thus, $\set{\nu,\lambda} \in E\{n-1\}$.

Further, $\set{\nu^*,\lambda^*} \not\prec \set{\lambda,\lambda}$.
\end{lemma}

\begin{lemma} \label{lemma:proof-of-item-En-sl(lam,pm)-d}
Let $n \geq 6$, and suppose Theorem \ref{thm:main-theorem} is true for the value $n-1$. Then part \eqref{item:En-sl(lam,pm)-d} of Theorem \ref{thm:main-theorem} is true for the value $n$, i.e.,
	\begin{center}
	If $\set{\lambda,\lambda} \in E\{n\}$, then $\fd_{\set{\lambda,\pm}} = \fsl\{\lambda,\pm \}$.
	\end{center}
\end{lemma}

\begin{proof}
Suppose $\set{\lambda,\lambda} \in E\{n\}$, and let $\fr = \fd_{n-1}$. Then the $\fd_n'$-module structure map $\rho_{\set{\lambda,\pm}}: \fd_n' \twoheadrightarrow \fd_{\set{\lambda,\pm}} \subseteq \fsl\{\lambda,\pm\}$ restricts to a Lie algebra homo\-morphism
	\begin{equation} \label{eq:r-En-sl(lam,pm)-d}
	\rho: \fr' \to \bigoplus_{\nu \prec \lambda} \fr_{\set{\nu,\lambda}}.
	\end{equation}
It follows from Lemma \ref{lemma:(lam,pm)-En-preceding} that the right side of \eqref{eq:r-En-sl(lam,pm)-d} is a summand of the right side of \eqref{eq:dn-refined-factorization} for the value $n-1$. This implies by induction that \eqref{eq:r-En-sl(lam,pm)-d} is a surjection, and implies for each $\nu \prec \lambda$ that $\fr_{\set{\nu,\lambda}} \cong \fsl\{\nu,\lambda\}$. Let $\fh = \rho_{\set{\lambda,\pm}}(\fr')$. Then $\fh \cong \rho(\fr')$ is a semisimple subalgebra of $\fd_{\set{\lambda,\pm}}$ such that the restriction of $S^{\set{\lambda,\pm}}$ to each simple ideal of $\fh$ admits a simple factor of multiplicity one. One has
	\[
	\rk(\fh) = \sum_{\nu \prec \lambda} \big[ \dim(S^{\set{\nu,\lambda}})-1 \big] = \dim(S^{\set{\lambda,\pm}}) - r_\lambda > \dim(S^{\set{\lambda,\pm}}) - \tfrac{n}{2},
	\]
and hence $2 \cdot \rk(\fh) - \dim(S^{\set{\lambda,\pm}}) > \dim(S^{\set{\lambda,\pm}}) - n$. As in \eqref{eq:half-dim-S(lam,lam*)} one sees that $\dim(S^{\set{\lambda,\pm}}) \geq 2n$, and hence $2 \cdot \rk(\fh) - \dim(S^{\set{\lambda,\pm}}) > 0$. This implies first that $\fd_{\set{\lambda,\pm}}$ is simple by \cite[Lemme 15]{Marin:2007}, and then that $\fd_{\set{\lambda,\pm}} \cong \fsl\{ \lambda,\pm \}$ by \cite[Lemme 13]{Marin:2007}.
\end{proof}

\subsection{Conclusion of the proof of Theorem \ref{thm:main-theorem}} \label{subsec:complete-main-theorem}

Given $(\lambda,\mu) \in \calBP(n)$, let $\ker(\rho_{(\lambda,\mu)})$ be the kernel of the $\fb_n'$-module structure map $\rho_{(\lambda,\mu)}: \fb_n' \twoheadrightarrow \fb_{(\lambda,\mu)} \subseteq \End(S^{(\lambda,\mu)})$, and let $\fb^{(\lambda,\mu)}$ be the orthogonal complement of $\ker(\rho_{(\lambda,\mu)})$ with respect to the Killing form on $\fb_n'$. Then $\fb^{(\lambda,\mu)}$ is an ideal in $\fb_n'$ such that $\fb_n' = \fb^{(\lambda,\mu)} \oplus \ker(\rho_{(\lambda,\mu)})$ as Lie algebras, and $\rho_{(\lambda,\mu)}$ induces a Lie algebra isomorphism $\fb^{(\lambda,\mu)} \cong \fb_{(\lambda,\mu)}$. We call $S^{(\lambda,\mu)}$ the module corresponding to $\fb^{(\lambda,\mu)}$. In an analogous fashion, one defines the ideals $\fd^{\set{\lambda,\mu}} \subseteq \fd_n'$ for $\set{\lambda,\mu} \in \calBP\{n\}$, and the ideals $\fd^{\set{\lambda,\pm}} \subseteq \fd_n'$ for $\set{\lambda,\lambda} \in \calBP\{n\}$. For $\lambda \vdash n$, set $\fb^\lambda = \fb^{(\lambda,\emptyset)}$ and $\fd^\lambda = \fd^{\set{\lambda,\emptyset}}$.

Recall the bipartitions $\beta$ and $\gamma$ defined in Section \ref{subsubsec:arm-and-leg}. Recall also the sets $E_n$ and $F_n$ and the partition $\alpha$ defined in Section \ref{subsubsec:improper-bipartitions}.

\begin{proposition} \label{prop:complete-simple-ideals-b}
Let $n \geq 6$ be fixed, suppose that Theorem \ref{thm:main-theorem} is true for the value $n-1$.
	\begin{enumerate}
	\item \label{item:complete-simple-ideals-b} The following collection consists of pairwise distinct simple ideals in $\fb_n'$:
		\[
		\{ \fb^\alpha,\fb^\beta,\fb^\gamma \} 
		\cup \{ \fb^\lambda : \lambda \in (E_n \cup F_n)/\!\!\sim \} 
		\cup \{ \fb^{(\lambda,\mu)} : (\lambda,\mu) \in [E(n) \cup F(n)]/\!\!\sim \}.
		\]
	
	\item \label{item:complete-simple-ideals-d} The following collection consists of pairwise distinct simple ideals in $\fd_n'$:
		\begin{align*}
		\{ \fd^\alpha,\fd^\beta \} 
		&\cup \{ \fd^\lambda : \lambda \in (E_n \cup F_n)/\!\!\sim \} \\
		&\mathrel{\cup} \{ \fd^{\set{\lambda,\mu}} : \set{\lambda,\mu} \in [E\{n\} \cup F\{n\}]/\!\!\sim \text{ and } \lambda \neq \mu \} \\
		&\mathrel{\cup} \{ \fd^{\set{\lambda,+}}, \fd^{\set{\lambda,-}}_\star : \set{\lambda,\lambda} \in [E\{n\} \cup F\{n\}]/\!\!\sim \},
		\end{align*}
	where the $\star$ indicates that $\fd^{\set{\lambda,-}}$ is omitted if $\lambda = \lambda^*$ and $\abs{\lambda}$ is odd.

	\end{enumerate}
\end{proposition}

\begin{proof}
We start with part \eqref{item:complete-simple-ideals-b}. It follows from Sections \ref{subsubsec:improper-bipartitions} and \ref{subsec:item:beta-gamma}--\ref{subsec:E(n)-sl(lam,mu)} that the collection consists of simple ideals of Lie types $A$, $C$, or $D$. (For ideals of type $\osp$, the corresponding module is even dimensional, so ideals of type $B$ do not occur.) Let $\fb^I$ be an ideal in the collection (so $I$ is one of the various superscript labels), and let $S^I$ be the corresponding module. Using Lemma \ref{lemma:Bn-subminimal-dimensions}, one observes that if $\fb^I$ is of type $A_r$, then $\dim(S^I) \geq n-1 \geq 5$, and hence $r \geq 4$, while if $\fb^I$ is of type $C_r$ or $D_r$, then $\dim(S^I) \geq \binom{n}{n/2} \geq \binom{6}{3} = 20$, and hence $r \geq 10$. Thus $\fb^I$ cannot be simultaneously labeled by two different Lie types.

Let $\fb^I$ and $\fb^J$ be two ideals from the collection. By simplicity, either $\fb^I \cap \fb^J = 0$ or $\fb^I = \fb^J$. Suppose $\fb^I = \fb^J$; we want to show that $I = J$. Let $\rho_I: \fb_n' \twoheadrightarrow \fb_I \subseteq \End(S^I)$ and $\rho_J: \fb_n' \twoheadrightarrow \fb_J \subseteq \End(S^J)$ be the structure maps for the corresponding modules $S^I$ and $S^J$. Since $\fb^I = \fb^J$, then $\dim(S^I) = \dim(S^J)$. We proceed based on the Lie type of $\fb^I$.

First suppose $\fb^I$ is of type $A_r$. Then fixing bases for $S^I$ and $S^J$, one has $\fb_I = \fsl(\C^{r+1}) = \fb_J$, and the composite induced map
	\[
	\fsl(\C^{r+1}) = \fb_I \xrightarrow{\rho_I^{-1}} \fb^I = \fb^J \xrightarrow{\rho_J} \fb_J = \fsl(\C^{r+1})
	\]
is a Lie algebra automorphism. Lie algebra automorphisms of $\fsl(\C^{r+1})$ come in two forms:
	\begin{itemize}
	\item $X \mapsto gXg^{-1}$ for some $g \in GL(\C^{r+1})$, or
	\item $X \mapsto -(gX^t g^{-1})$ for some $g \in GL(\C^{r+1})$, where $X^t$ denotes the transpose of $X$;
	\end{itemize}
see \cite[IX.5]{Jacobson:1979}. If $\rho_J \circ \rho_I^{-1}$ is of the first form, it follows that $S^I \cong S^J$ as $\fb_n'$-modules, implying by Proposition \ref{prop:iso-equivalent} and the range of possible values for $I$ and $J$ in \eqref{item:complete-simple-ideals-b} that $I = J$. If $\rho_J \circ \rho_I^{-1}$ is of the second form, it follows that $S^I \cong [S^J]^{*,\Lie}$ as $\fb_n'$-modules. However, since $[S^{(\lambda,\mu)}]^{*,\Lie} \cong S^{(\mu^*,\lambda^*)}$, and since the natural representation of $\fsl(k)$ is not self-dual for $k \geq 3$, it follows from Proposition \ref{prop:iso-equivalent} that $S^I \not\cong [S^J]^{*,\Lie}$ for any $I$ and $J$ in the range of possible values in \eqref{item:complete-simple-ideals-b}. Thus the automorphism must be of the first form, and $I = J$.

Now suppose $\fb^I$ is of type $C_r$ or $D_r$. Then $\dim(S^I) = \dim(S^J) = 2r$. Up to isomorphism, a simple Lie algebra of type $C_r$ or $D_r$ has a unique irreducible representation of dimension $2r$. Then $S^I \cong S^J$ as $\fb^I$-modules, and hence $S^I \cong S^J$ as $\fb_n'$-modules. Once again, this implies by Proposition \ref{prop:iso-equivalent} and the choice of possible values for $I$ and $J$ that $I = J$.

The proof of part \eqref{item:complete-simple-ideals-d} of the proposition proceeds along the same lines, using the inductive results from Sections \ref{subsubsec:improper-bipartitions}, \ref{subsec:item:beta-gamma}, and \ref{subsec:Fn-osp(lam,mu)}--\ref{subsec:En-sl(lam,pm)-d}. If an ideal $\fd^I$ is of type $C_r$ or $D_r$, then one observes using Lemma \ref{lemma:Dn-subminimal-dimensions} that $\dim(S^I) \geq \frac{1}{2} \cdot \binom{n}{n/2} \geq \frac{1}{2} \cdot \binom{6}{3} = 10$, and hence $r \geq 5$. The remainder of the argument is the same.
\end{proof}

\begin{proposition} \label{prop:main-theorem-surjectivity}
Let $n \geq 6$ be fixed, suppose that Theorem \ref{thm:main-theorem} is true for the value $n-1$. Then the (injective) maps \eqref{eq:bn-refined-factorization} and \eqref{eq:dn-refined-factorization} are surjections.
\end{proposition}

\begin{proof}
It follows from Proposition \ref{prop:complete-simple-ideals-b}\eqref{item:complete-simple-ideals-b} that the sum of ideals
	\[
	\fb^\alpha + \Big[ \sum_{\lambda \in E_n/\sim} \fb^\lambda \Big] + \Big[ \sum_{\lambda \in F_n} \fb^\lambda \Big] 
	+ \fb_{\beta} + \fb_{\gamma} 
	+ \Big[ \sum_{(\lambda,\mu) \in E(n)/\sim} \fb_{(\lambda,\mu)} \Big] 
	+ \Big[ \sum_{(\lambda,\mu) \in F(n)} \fb_{(\lambda,\mu)} \Big]
	\]
is a direct sum in $\fb_n'$. This implies by dimension comparison that \eqref{eq:bn-refined-factorization} is a surjection. Similarly, one sees that \eqref{eq:dn-refined-factorization} is a surjection by applying Proposition \ref{prop:complete-simple-ideals-b}\eqref{item:complete-simple-ideals-d}.
\end{proof}

\makeatletter
\renewcommand*{\@biblabel}[1]{\hfill#1.}
\makeatother

\bibliographystyle{eprintamsplain}
\bibliography{lie-algebras-reflections}

\end{document}